\documentclass{amsart}
\usepackage{amssymb}
\usepackage{amsmath}
\usepackage{lscape}

\usepackage[pdftex,colorlinks,backref,pagebackref,hypertexnames=false,
            linkcolor=blue,urlcolor=blue,citecolor=blue]{hyperref}

\newcommand{\N}{{\mathbb{N}}}
\newcommand{\Z}{{\mathbb{Z}}}
\newcommand{\Q}{{\mathbb{Q}}}
\newcommand{\R}{{\mathbb{R}}}
\newcommand{\C}{{\mathbb{C}}}
\newcommand{\F}{{\mathbb{F}}}
\newcommand{\G}{{\mathbf{G}}}
\newcommand{\T}{{\mathbf{T}}}
\newcommand{\W}{{\mathbf{W}}}

\newcommand{\LL}{{\mathbf{L}}}

\newcommand{\zetat}{\tilde{\zeta}}
\newcommand{\xit}{\tilde{\xi}}

\newcommand{\phit}{\tilde{\varphi}}
\newcommand{\epsilont}{\tilde{\varepsilon}}
\newcommand{\psit}{\tilde{\psi}}

\newcommand{\Irr}{\text{Irr}}
\newcommand{\Hom}{\operatorname{Hom}}
\newcommand{\Aut}{\operatorname{Aut}}
\newcommand{\Out}{\operatorname{Out}}
\newcommand{\Sz}{\operatorname{Sz}}
\newcommand{\Blk}{{\text{Blk}}}

\renewcommand{\epsilon}{\varepsilon}
\newcommand{\ra}{\rangle}
\newcommand{\la}{\langle}

\newtheorem{theorem}{Theorem}[section]
\newtheorem{lemma}[theorem]{Lemma}
\newtheorem{proposition}[theorem]{Proposition}

\begin{document}

\title[Dade's Invariant Conjecture for ${^2F_4(q^2)}$ in Defining 
  Characteristic]{Dade's Invariant Conjecture for the Ree
  Groups~$\mathbf{{^2F_4(q^2)}}$ in Defining Characteristic}      

\author{Frank Himstedt and Shih-chang Huang}

\address{F.H.: Technische Universit\"at M\"unchen, Zentrum Mathematik --
         M11, Boltzmannstr. 3, 85748 Garching, Germany}
\address{S.H.: Department of Mathematics, National Cheng Kung University,
  No. 1 Dasyue Rd, Tainan City 70101, Taiwan}

\email{F.H.: himstedt@ma.tum.de}
\email{S.H.: shua3@mail.ncku.edu.tw}


\begin{abstract}
We verify Dade's invariant conjecture for the simple Ree groups
${^2F}_4(2^{2n+1})$ for all $n>0$ in the defining characteristic,
i.e., in characteristic~2. Together with the results in \cite{An2F4},
this completes the proof of Dade's conjecture for the simple Ree
groups ${{}^2F}_4(2^{2n+1})$. 
\end{abstract}

\maketitle


\section{Introduction}

The representation theory of finite groups received great impetus from
several conjectures linking the representation theory of a finite
group $G$ to that of local subgroups, the normalizers of nontrivial
$p$-subgroups of $G$. A prominent family of such conjectures are
Dade's conjectures, published in \cite{Dade1}, \cite{Dade2} and 
\cite{Dade29}. These conjectures express the number of complex
irreducible characters with a fixed defect in a given $p$-block of a
finite group $G$ as an alternating sum of related values for
$p$-blocks of certain $p$-local subgroups. 

Dade's conjectures are closely related to several other
local-global conjectures: Dade's projective conjecture implies
Alperin's weight conjecture, the McKay conjecture and its blockwise
version due to Alperin, and one direction of Brauer's height-0
conjecture on abelian blocks. For details 
see~\cite{AlpMac}, \cite{Alpweight}, \cite{Dade1}, \cite{Dade2},
 \cite{KnoerrRob}, \cite{McKay} and also \cite{IsaacsNavarro},
 \cite{IsaacsMalleNavarro}. 

In this paper, we show that Dade's invariant conjecture is true for the
simple Ree groups ${^2F}_4(2^{2n+1})$, $n>0$, in the defining
characteristic, i.e., for $p=2$. Since the Ree groups
have a trivial Schur multiplier and a cyclic outer automorphism group
for $n>0$, our results imply that the most general version
of Dade's conjectures, the inductive conjecture, 
is also true for~${^2F}_4(2^{2n+1})$ and $p=2$. Together with An's
results~\cite{An2F4} and \cite{AnTits}, this completes the proof of
Dade's conjectures for the simple groups ${^2F}_4(2^{2n+1})'$ for all
$n \ge 0$. Note that this also completes the verification of Uno's
invariant conjecture \cite{Uno} for these groups in the defining
characteristic.  

The methods we use are similar to those in
\cite{HimstedtHuang3D4even}. By a corollary to the Borel-Tits theorem
\cite{BW}, the normalizers of radical $2$-chains are the parabolic
subgroups of ${{}^2F}_4(2^{2n+1})$. Using the character tables of
${{}^2F}_4(2^{2n+1})$ in the library of the MAPLE~\cite{Maple} part of
CHEVIE~\cite{CHEVIE} and the character tables of the parabolic
subgroups in \cite{HimstedtHuang2F4Borel} and
\cite{HimstedtHuang2F4MaxParab}, we count characters of the chain
normalizers which are fixed by certain outer automorphisms of the Ree
groups and then evaluate the alternating sums in Dade's invariant
conjecture. 

There are some differences between the verification of Dade's
conjecture for the Ree groups ${^2F}_4(2^{2n+1})$ and
for Steinberg's triality groups ${^3D_4}(2^n)$
in~\cite{HimstedtHuang3D4even}. Firstly, the parabolic subgroups of  
${^2F}_4(2^{2n+1})$ have more types of characters which have to be
treated separately. And secondly, when calculating the fixed points of
outer automorphisms on the irreducible characters of the parabolic
subgroups of ${^2F}_4(2^{2n+1})$ one has to deal with cyclotomic
polynomials over $\Q(\sqrt{2})$ instead of $\Q$ making
$\gcd$-computations more complicated; see for example
Lemma~\ref{la:gcds1}. 

This paper is organized as follows: In Section~2, we fix notation and
state Dade's invariant conjecture. In Section~3, we prove some lemmas
from elementary number theory which we use to count fixed points of
certain automorphisms of ${{}^2F}_4(2^{2n+1})$. In Section~4, we
compute the fixed points of these outer automorphisms on the
irreducible characters of parabolic subgroups. In Section~5, we verify
Dade's invariant conjecture for~${^2F}_4(2^{2n+1})$ in the defining
characteristic. Details on irreducible characters of the Ree groups
are summarized in tabular form in Appendix~A. In  Appendix B, we
describe the construction of several families of irreducible
characters of~${^2F}_4(2^{2n+1})$. Some of the data in Appendix B
might also be of independent interest and is used in~\cite{Himstedt2F4Decomp}. 

\section*{Acknowledgments}

\noindent A part of this work was done during visits of the first
author at the University of Auckland and at Chiba University. He
wishes to express his sincere thanks to all persons of 
the mathematics departments of these universities for their
hospitality, and also to the Marsden Fund of New Zealand and the Japan
Society for the Promotion of Science (JSPS) who supported his
visits. The second author acknowledges the support of his research
during the last years from the Foundation for Research, Science and
Technology of New Zealand, the JSPS and the National Science Council,
ROC. 


\section{Dade's Invariant Conjecture} \label{sec:dadeinvconj}

Let $R$ be a $p$-subgroup of a finite group $G$. Then $R$ is {\it
  radical} if $O_p(N(R))=R$, where $O_p(N(R))$ is the largest normal
  $p$-subgroup of the normalizer $N(R) := N_G(R)$. Denote by $\Irr(G)$
  the set of all irreducible ordinary characters of $G$, and by $\Blk(G)$ the set of
$p$-blocks. If $H \le G$, $\widetilde{B} \in \Blk(G)$, and $d$ is an integer,
we denote by $\Irr(H,\widetilde{B},d)$ the set of characters $\chi \in \Irr(H)$
satisfying $\text{d}(\chi) = d$ and $b(\chi)^G = \widetilde{B}$
  (in the sense of Brauer), where $\text{d}(\chi) = \log_p(|H|_p) -
  \log_p(\chi(1)_p)$ is the $p$-defect of $\chi$
  and $b(\chi)$ is the block of $H$ containing $\chi$.

Given a $p$-subgroup chain
$
C: P_0< P_1<\cdots< P_n
$
of $G$, define the length $|C| := n$, $C_k: P_0 < P_1 < \cdots < P_k$ and
\[
N(C) = N_G(C) := N_G(P_0)\cap N_G(P_1)\cap\cdots\cap N_G(P_n).
\]
The chain $C$ is said to be {\it radical} if it satisfies the
following two conditions:

\begin{enumerate}
\item[(a)] $P_0 = O_p(G)$ and
\item[(b)] $P_k = O_p(N(C_k))$ for $1\le k\le n$.
\end{enumerate}

\noindent Denote by $\mathcal{R} = \mathcal{R}(G)$ the set of all
radical $p$-chains of $G$.

\medskip
Suppose $1 \rightarrow G \rightarrow E \rightarrow \overline{E} \rightarrow 1$
 is an exact sequence, so that $E$ is an extension of $G$ by $\overline{E}$.
Then $E$ acts on $\mathcal{R}$ by conjugation.
Given $C \in \mathcal{R}$ and $\psi \in \Irr(N_G(C))$, let $N_E(C,\psi)$ be
the stabilizer of $(C,\psi)$ in  $E$, and
$
N_{\overline{E}}(C,\psi) := N_E(C,\psi)/N_G(C).
$
For $\widetilde{B} \in \Blk(G)$, an integer $d \ge 0$ and $U \le {\overline{E}}$, let
$k(N_G(C), \widetilde{B}, d, U)$ be the number of characters in the set
\[
\Irr(N_G(C),\widetilde{B},d,U) := \{\psi \in
\Irr(N_G(C),\widetilde{B},d) \, | \, N_{\overline{E}}(C,\psi)=U\}.
\]
Dade's invariant conjecture can be stated as follows:

\medskip

\noindent {\bf Dade's Invariant Conjecture} (\cite{Dade29}) \enspace
{\it If $O_p(G)=1$ and $\widetilde{B}\in\Blk(G)$ with defect group 
$D(\widetilde{B}) \neq 1$, then
\[
\sum_{C \in {\mathcal R}/G}(-1)^{|C|} k(N_{G}(C),\widetilde{B},d,U) = 0,
\]
where ${\mathcal R}/G$ is a set of representatives for the $G$-orbits
of $\mathcal R$.}

\smallskip

Let $\Aut(G)$ and $\Out(G)$ be the automorphism and outer automorphism groups 
of $G$, respectively. We may suppose ${\overline {E}}=\Out(G)$. If moreover, 
$\Out(G)$ is cyclic, then we write
\[
k(N_G(C), \widetilde{B}, d, |U|) := k(N_G(C), \widetilde{B}, d, U).
\]
For $G ={^2F}_4(2^{2n+1})$, $n>0$, $\Out(G)$ is cyclic and the  Schur
multiplier of $G$ is trivial. So the invariant conjecture for $G$ is
equivalent to the inductive conjecture.

\medskip


\section{Notation and Lemmas from Elementary Number Theory}
\label{sec:lemmas}

From now on, we assume that $p = 2$, $n$ is a 
positive integer and $q^2 = 2^{2n+1}$. We write~$\phi_i$ for the
$i$-th cyclotomic polynomial in $q$, for example: $\phi_1=q-1$,
$\phi_2=q+1$, $\phi_3=q^2+q+1$, $\phi_4=q^2+1$, $\phi_6=q^2-q+1$, $\phi_8=q^4+1$,
$\phi_{12}=q^4-q^2+1$, $\phi_{24}=q^8-q^4+1$. Furthermore, we set
$\phi_8' := q^2+\sqrt{2}q+1$, $\phi_8'' := q^2-\sqrt{2}q+1$, 
$\phi_{24}'=q^4+\sqrt{2}q^3+q^2+\sqrt{2}q+1$,
$\phi_{24}''=q^4-\sqrt{2}q^3+q^2-\sqrt{2}q+1$.
We denote by $\N = \{0, 1, 2, \dots \}$ the set of natural numbers
including zero. In the next section, we will use the following lemmas,
the first one is given by \cite[Lemma 3.1]{AnHimstedtHuang3D4odd}. 

\begin{lemma} \label{la:gcdfolk}
Suppose $m,n \in \Z$ with $m,n > 0$. Then 
$\gcd(2^m-1, 2^n-1) =|2^d-1|$ where $d := \gcd(m,n)$. 
\end{lemma}

\medskip

\begin{lemma} \label{la:gcds}
Let $f=2n+1$ and $t$ be a positive integer with $t\, |\, f$.  
Then the following hold.
\begin{enumerate}
\item[(i)] $\gcd(2^t-1, q^2-1) \hspace{0.15cm} = \, 2^t - 1$.

\item[(ii)] $\gcd(2^t-1, q^2+1) \hspace{0.15cm} = \, 1$.



\item[(iii)] $\gcd(2^t+1, q^2-1) \hspace{0.15cm} = \, 1$.

\item[(iv)] $\gcd(2^t+1, q^2+1) \hspace{0.15cm} = \, 2^t+1$.


\end{enumerate}
\end{lemma}

\noindent {\bf Proof:} (i) is clear by Lemma~\ref{la:gcdfolk}.

\smallskip

\noindent (ii) Suppose $d=\gcd(2^t-1, q^2+1)$. Since 
 $2^t-1 \, | \, q^2-1$ by (i), it follows that 
$d \mid \gcd(q^2-1, q^2+1)=1$.  


\smallskip

\noindent (iv) Since $t$ and $f$ are both odd, we have 
$2^t + 1 = -((-2)^t - 1)$ and $q^2 + 1 = -((-2)^f - 1)$. 
So Lemma~\ref{la:gcdfolk} implies
$\gcd(2^t+1, q^2+1) = \gcd((-2)^t - 1, (-2)^f - 1)
= |(-2)^t - 1| = 2^t + 1.$


\smallskip

\noindent (iii)  Suppose $d=\gcd(2^t+1, q^2-1)$. Since 
 $2^t+1 \, | \, q^2+1$ by (iv), it follows that 
$d \mid \gcd(q^2-1, q^2+1)=1$. \hfill $\Box$


\medskip

\begin{lemma} \label{la:gcds1}
Let $f=2n+1$ and $t$ be a positive integer with $t\, |\, f$. Then
there is a nonnegative integer $m$ such that $f-t=4tm$ or
$f-t=2t \cdot (2m+1)$ and for $\epsilon \in \{1,-1\}$ the following holds:
\begin{enumerate}
\item[(i)] If $f-t=4tm$, then
\begin{eqnarray*}
\gcd(2^{f-t}-1, q^2 + \epsilon \sqrt{2} q+1) & = & 2^t + \epsilon
\cdot (-1)^m \cdot 2^{(t+1)/2} + 1,\\
\gcd(2^{f-t}+1, q^2 + \epsilon \sqrt{2} q+1) & = & 1.
\end{eqnarray*}

\item[(ii)] If $f-t=2t \cdot (2m+1)$, then
\begin{eqnarray*}
\gcd(2^{f-t}-1, q^2 + \epsilon \sqrt{2}q+1) & = & 1,\\
\gcd(2^{f-t}+1, q^2 + \epsilon \sqrt{2}q+1) & = & 2^t + \epsilon \cdot
(-1)^{m+1} \cdot 2^{(t+1)/2} + 1.
\end{eqnarray*}
\end{enumerate}
\end{lemma}

\noindent {\bf Proof:} Since $t \, | \, f=2n+1$, the integer $t$ is
odd and the even integer $f-t$ is a multiple of $2t$. So, there is 
$m \in \N$ such that $f-t = 2t \cdot 2m$ or $f-t = 2t \cdot (2m+1)$.

\smallskip

\noindent (i) Suppose $f-t = 4tm$ and $\epsilon \in \{1,-1\}$. We
proceed in several steps:

\smallskip

\noindent Step 1: $\gcd(2^{f-t}-1, q^4+1) = 2^{2t}+1$.

\noindent Since $t \, | \, f$, there is an odd natural number $f_u$ such that
$f = t \cdot f_u$. So $f - t = t \cdot (f_u - 1) = 4t \cdot \frac{f_u - 1}{4}$. Hence 
$2^{f-t} - 1 = (-2^{2t})^{2 \cdot \frac{f_u-1}{4}} - 1$ and 
$q^4 + 1 = -((-2^{2t})^{f_u} - 1)$. So Lemma~\ref{la:gcdfolk}
implies  $
\gcd(2^{f-t}-1, q^4+1) = \gcd((-2^{2t})^{2 \cdot \frac{f_u-1}{4}} - 1, (-2^{2t})^{f_u} - 1)
= |(-2^{2t}) - 1| = 2^{2t} + 1.$

\smallskip

\noindent Step 2: The number $q^2 + \epsilon \sqrt{2} q+1$ is a
multiple of $2^t + \epsilon \cdot (-1)^m \cdot 2^{(t+1)/2} + 1$.

\noindent Since $f-t=2n+1-t=4tm$, we have $n = (4tm+t-1)/2$. 
Modulo $2^{2t}+1$, we compute:
\begin{eqnarray*}
q^2 + \epsilon \sqrt{2} q+1 & \equiv & 2^{2n+1} + \epsilon \cdot 2^{n+1} +
1 = 2^{4tm+t} + \epsilon \cdot 2^{(4tm+t+1)/2} + 1\\
& \equiv & (2^{2t})^{2m} \cdot 2^t + \epsilon \cdot (2^{2t})^m \cdot
2^{(t+1)/2} + 1\\
& \equiv & 2^t + \epsilon \cdot (-1)^m \cdot 2^{(t+1)/2} + 1.
\end{eqnarray*}
So, since $2^{2t}+1 = (2^t + \epsilon \cdot (-1)^m \cdot 2^{(t+1)/2} +
1) \cdot (2^t - \epsilon \cdot (-1)^m \cdot 2^{(t+1)/2} + 1)$, we can
conclude that $q^2 + \epsilon \sqrt{2} q+1$ is a multiple of 
$2^t + \epsilon \cdot (-1)^m \cdot 2^{(t+1)/2} + 1$.

\smallskip

\noindent Step 3: $\gcd(2^{f-t}-1, q^2 + \epsilon \sqrt{2} q+1) = 
2^t + \epsilon \cdot (-1)^m \cdot 2^{(t+1)/2} + 1$.

\noindent Since $q^4+1 = (q^2+\sqrt{2}q+1) \cdot (q^2-\sqrt{2}q+1)$
and $\gcd(q^2+\sqrt{2}q+1, q^2-\sqrt{2}q+1)=1$, we obtain from
Step~1:
\begin{eqnarray*}
&& \gcd(2^{f-t}-1, q^4+1)\\
& = & \gcd(2^{f-t}-1, q^2 + \epsilon \sqrt{2} q+1)
\cdot \gcd(2^{f-t}-1, q^2 - \epsilon \sqrt{2} q+1)\\
& = & 2^{2t} + 1 = (2^t + \epsilon \cdot (-1)^m \cdot 2^{(t+1)/2} +
1) \cdot (2^t - \epsilon \cdot (-1)^m \cdot 2^{(t+1)/2} + 1).
\end{eqnarray*}
So, Step~2 implies $\gcd(2^{f-t}-1, q^2 + \epsilon \sqrt{2} q+1) = 2^t
+ \epsilon \cdot (-1)^m \cdot 2^{(t+1)/2} + 1$. 

\smallskip

\noindent Step 4: $\gcd(2^{f-t}+1, q^2 + \epsilon \sqrt{2} q+1) = 1$.

\noindent Let $d = \gcd(2^{f-t}+1, q^4+1)$. As in Step~1, we have $f =
t \cdot f_u$ and $f - t = t \cdot (f_u - 1) = 4t
\cdot \frac{f_u - 1}{4}$ for an odd integer $f_u$. 
Since $d \, | \, 2^{f-t}+1$ and $d \, | \, 2^{2f}+1 = q^4+1$, we have 
$d \, | \, 2^{2(f-t)}-1$ and $d \, | \, 2^{4f}-1$. 
By Lemma~\ref{la:gcdfolk}, it follows that
$\gcd(2^{2(f-t)}-1, 2^{4f}-1) = \gcd(2^{4t \cdot 2 \cdot
  \frac{f_u-1}{4}}-1, 2^{4t \cdot   f_u}-1) = 2^{4t}-1$. 
Then $d \, | \, 2^{4t}-1$, so $2^{4t} \equiv 1$ mod $d$ and 
$0 \equiv 2^{f-t}+1 = 2^{4t \cdot \frac{f_u-1}{4}}+1 \equiv 2 
\mbox{ mod } d.$ So $d = 1$ as $d$ is odd. Since 
$q^2 + \epsilon \sqrt{2} q+1$ is a divisor of $q^4+1$, the claim
follows. 

\medskip

\noindent (ii) Suppose $f-t = 2t \cdot (2m+1)$ and 
$\epsilon \in \{1,-1\}$. Let $d = \gcd(2^{f-t}-1, q^4+1)$. Then
there is an odd $f_u \in \N$ such that $f = t \cdot f_u$, $f - t = t
\cdot (f_u - 1) = 2t \cdot \frac{f_u - 1}{2}$. Since $d \, | \, q^4+1 =
2^{2f}+1$, we have $d \, | \, 2^{4f}-1$.
By Lemma~\ref{la:gcdfolk}, it follows that $\gcd(2^{4f}-1, 2^{f-t}-1)
= \gcd(2^{2t \cdot 2f_u}-1, 2^{2t \cdot \frac{f_u-1}{2}}-1) = 2^{2t}-1$.
Then $d \, | \, 2^{2t}-1$, so $2^{2t} \equiv 1$ mod $d$ and
$0 \equiv q^4+1 = 2^{2f}+1 = 2^{2t \cdot f_u} + 1 \equiv 2
\mbox{ mod } d.$ So $d = 1$ as $d$ is odd. Since 
$q^2 + \epsilon \sqrt{2} q+1$ is a divisor of $q^4+1$, the first
statement of (ii) follows. The proof of the second statement of (ii) is
analogous to Steps 1-3 of (i). \hfill $\Box$

\medskip


\section{Action of Automorphisms on Irreducible Characters} 
\label{sec:action}

Let $G = {^2F_4(q^2)}$ be the simple Ree group with $q^2 = 2^{2n+1}$
and $n$ a positive integer. 
Let $O = \Out(G)$ and $A = \Aut(G)$. Then $O = \la \alpha \ra$
and $A = G \rtimes \la \alpha \ra$, where $\alpha$ is a field
automorphism of (odd) order $2n+1$. We fix a Borel subgroup $B$ and maximal
parabolic subgroups $P_a$ and $P_b$ of~$G$ containing $B$ as in
\cite{HimstedtHuang2F4Borel} and \cite{HimstedtHuang2F4MaxParab}. In
particular, $\alpha$ stabilizes $B$, $P_a$ and $P_b$.   

In this section, we determine the action of $O = \Out(G)$ on the irreducible
characters of the chain normalizers. Our notation for the parameter sets of
the irreducible characters of $G$, $B$, $P_a$ and $P_b$ is similar to the 
CHEVIE notation and is given in Table~A.1 in Appendix A. For a
definition and construction of the irreducible characters
$\chi_{42}(k), \dots, \chi_{51}(k,l)$, see Appendix B.

The first column of Table A.1 defines a name for the parameter
set which parameterizes those characters which are listed in the second column
of the table. The characters of $G$ are numbered according to the character
table of ${^2F}_4(q^2)$ in CHEVIE, and for the characters of 
$B$, $P_a$, $P_b$ we use the notation
from~\cite{HimstedtHuang2F4Borel} and~\cite{HimstedtHuang2F4MaxParab}. 
The list of parameters in the third column of Table~A.1 in Appendix A is of
the form  
\[
k = 0, \dots, n_1-1 \quad \quad \text{   or   } \quad \quad 
\begin{array}{lll}
k & = & 0, \dots, n_1-1\\
l & = & 0, \dots, n_2-1
\end{array}
\]
where the integers $n_j$ are given by polynomials in $q$ with
coefficients in $\Z[\sqrt{2}]$. In the first case, the parameter $k$ can be
substituted by an element of $\Z$, but two parameters which differ by
an element of $n_1 \Z$ yield the same character. In the second case,
the parameter vector $(k,l)$ can be substituted by an element of 
$\Z \times \Z$, but two parameter vectors which differ by an element of
$n_1 \Z \times n_2 \Z$ yield the same character. In other words, $k$
can be taken to be an element of $\Z_{n_1}$ and $(k,l)$ can be taken
to be an element of $\Z_{n_1} \times \Z_{n_2}$. The groups $\Z_{n_1}$
and $\Z_{n_1} \times \Z_{n_2}$ are also called \textsl{character
parameter groups} (see Section 3.7 of the CHEVIE \cite{CHEVIE}
manual). The next lines of Table~A.1 list elements which have to be
excluded from the character parameter group. The remaining parameters
are called \emph{admissible} in the following. Different values of
admissible parameters may give the same character. The fourth column
of Table~A.1 defines an equivalence relation on the set of admissible
parameters. If no equivalence relation is listed we mean the identity
relation. The parameter set is defined to be the set of these 
equivalence classes. Finally, the last column of Table~A.1 gives the
cardinality of the parameter set.   

We consider the example $_{P_b}I_5$. The character parameter group is
$\Z_{q^2-1} \times \Z_{q^2-1}$. The parameter vectors $(k,l)$ and $(k,-l)$
yield the same character and the equivalence class of $(k,l)$ is 
$\{(k,l), (k,-l)\}$. Hence, the characters ${_{P_b}\chi}_5(k,l)$ are
parameterized by the set 
\[
_{P_b}I_5 = \{\{(k,l), (k,-l)\} \, | \, (k,l) \in \Z_{q^2-1} \times \Z_{q^2-1} \,
, \, l\neq 0\}.
\]
If we want to emphasize the dependence of a parameter set, say
$_{P_b}I_5$, from $q$ we write $_{P_b}I_5(q)$. Table~A.1 does not give any detailed
information about the parameter sets $_GI_{34}$, $_GI_{35}$, $_GI_{36}$, $\ldots$,
$_GI_{41}$, $_GI_{42}$, $_GI_{46}$, $_GI_{50}$, $_GI_{51}$, since we
will not need an explicit knowledge of these sets 
(note that the sets $_GI_{1}$, $_GI_{22}$, $_GI_{26}$, $_GI_{28}$,
$_GI_{29}$, $_GI_{32}$, $_GI_{34}$, $_GI_{35}$, $\ldots$, $_GI_{41}$,
$_GI_{42}$, $_GI_{46}$, $_GI_{50}$, $_GI_{51}$ parameterize the semisimple
irreducible characters of $G$). The data in Table~A.1 is taken from
CHEVIE library, the Appendix of \cite{HimstedtHuang2F4Borel} and
\cite{HimstedtHuang2F4MaxParab} and Appendix~B of this paper.

\smallskip

The action of $O = \Out(G)$ on the conjugacy classes of elements of
$G$, $B$, $P_a$ and $P_b$ induces an action of $O$ on the sets 
$\Irr(G)$, $\Irr(B)$, $\Irr(P_a)$ and 
$\Irr(P_b)$ and then an action on the parameter sets. Using the
values of the irreducible characters of $G$, $B$, $P_a$ and $P_b$ on the
classes listed in the last column of Tables~A.2-A.15 we can describe the
action of $O$ on the parameter sets.

For an $O$-set $I$ and each subgroup $H \le O$ let $C_I(H)$
denote the set of fixed points of $I$ under the action of $H$.
In the following proposition we determine $|C_I(H)|$ where~$I$ runs through
all (disjoint) unions of parameter sets which are listed in Table~A.16 except
for ${_GI^{ss}}$, where ${_GI^{ss}}$ is the disjoint union of
parameter sets parameterizing the semisimple irreducible characters of
$G$. The number of fixed points of ${_GI^{ss}}$ under the action of
$H$ was computed in \cite[Lemma~6.2]{HimstedtHuang2F4Borel}.   

\begin{proposition} \label{prop:nrfixedpoints}
Let $t \, | \, 2n+1$ and $I \neq {_GI^{ss}}$ 
be one of the (disjoint) unions of parameter sets listed in Table~A.16. If
$H = \la \alpha^t \ra$ is a subgroup of $O$, then the second column of
Table~A.16 shows the number of fixed points $|C_I(H)|$ of 
$I$ under the action of $H$.  
\end{proposition}

\noindent {\bf Proof:} We have to consider the following parameter sets $I$.

\medskip

\noindent First let \\ 
\noindent $I \in \{{_GI_{2}} \cup {_GI_{3}}\, , \, {_GI_{4}} \cup {_GI_{30}}\, ,\, {_GI_{5}} \cup {_GI_{6}} \cup {_GI_{8}} \cup {_GI_{9}} \cup {_GI_{11}} \cup {_GI_{12}} \cup {_GI_{13}} \cup {_GI_{14}}\, ,\\ 
\, 
\hspace*{0.8cm} {_GI_{7}} \cup {_GI_{10}}\, , \, {_GI_{17}}\, ,  \,
        {_GI_{18}}\, , \, {_GI_{19}} \cup {_GI_{20}}\, ,  \,
        {_GI_{31}}\, , {_BI_{6}} \cup {_BI_{7}}\, ,  \,{_BI_{8}}\, ,
        \, {_BI_{10}}\, , {_BI_{13}} \cup {_BI_{14}} \, ,\\
\hspace*{0.8cm}{_BI_{15}} \cup {_BI_{16}} \cup {_BI_{17}} \cup {_BI_{18}} \cup
             {_BI_{19}} \cup {_BI_{20}} \cup {_BI_{21}} \cup {_BI_{22}}\, , \, {_BI_{28}}\cup {_BI_{29}}\, , \, {_BI_{43}}\, ,\\
\hspace*{0.8cm} {_BI_{53}} \cup {_BI_{54}}\, , \, {_BI_{55}} \cup {_BI_{56}}\cup {_BI_{57}} \cup {_BI_{58}} \, , \,{_{P_a}I_{6}}\, , \,{_{P_a}I_{9}} \cup {_{P_a}I_{10}}\, , \,{_{P_a}I_{11}}\, , \,{_{P_a}I_{13}} \cup {_{P_a}I_{25}}\, , \\
\hspace*{0.8cm} {_{P_a}I_{14}} \cup {_{P_a}I_{15}} \cup {_{P_a}I_{16}}
\cup {_{P_a}I_{17}} \cup {_{P_a}I_{21}} \cup {_{P_a}I_{22}} \cup
     {_{P_a}I_{23}} \cup {_{P_a}I_{24}}\, , \,{_{P_a}I_{20}}\, , \,
     {_{P_a}I_{29}}\, , \,{_{P_a}I_{30}}\, , \\
\hspace*{0.8cm}{_{P_a}I_{31}}\, , \,{_{P_b}I_{9}} \cup {_{P_b}I_{10}}
\, , \, {_{P_b}I_{11}} \, , \, {_{P_b}I_{13}}\, , \, {_{P_b}I_{14}}
\cup {_{P_b}I_{15}} \, , \, {_{P_b}I_{17}} \cup {_{P_b}I_{18}} \cup
     {_{P_b}I_{19}} \cup {_{P_b}I_{20}} \\
\hspace*{0.8cm} \cup {_{P_b}I_{22}} \cup {_{P_b}I_{23}} \cup
        {_{P_b}I_{24}} \cup {_{P_b}I_{25}}\, , \, {_{P_b}I_{29}} \cup
        {_{P_b}I_{30}}\, , \, {_{P_b}I_{43}} \cup {_{P_b}I_{50}} \, ,
        \, {_{P_b}I_{44}}\, \cup {_{P_b}I_{45}} \\
\hspace*{0.8cm} \cup {_{P_b}I_{51}} \cup {_{P_b}I_{52}}\, , \, {_{P_b}I_{46}} \cup {_{P_b}I_{53}}\}$. 

\noindent The degrees and character values on the conjugacy classes listed in
Tables~A.2-A.15 show $C_I(H) = I$ and hence $|C_I(H)| = |I|$.  We demonstrate this
for the parameter set $I = {_{P_b}I_{43}} \cup {_{P_b}I_{50}}$. The degrees in Table~A.5 show that ${_{P_b}\chi_{43}}$ and ${_{P_b}\chi_{50}}$
are the only irreducible characters of $P_b$ of degree
$\frac{q^9}{\sqrt{2}}(q^2-1)$.
Hence, ${_{P_b}\chi_{43}}^\alpha \in \{{_{P_b}\chi_{43}}, {_{P_b}\chi_{50}}\}$. The class representatives
in Table~A.7 in \cite{HimstedtHuang2F4MaxParab} show that the conjugacy class
$c_{1,9}$ is fixed by $\alpha$ and we can see from the character table A.9 of
$P_b$ in \cite{HimstedtHuang2F4MaxParab} that the values of ${_{P_b}\chi_{43}}$ and ${_{P_b}\chi_{50}}$
on $c_{1,9}$ are different. So, ${_{P_b}\chi_{i}}^\alpha = {_{P_b}\chi_{i}}$ for $i =
43$, $50$ and $|C_I(H)| = |I|$.

\noindent In each of the following cases, we have that the action of $\alpha$ on $I$ is given by $x^\alpha = 2x$ for all $x\in I$ using the character values
on the classes listed in the last column of Tables A.5-A.15.
We demonstrate this for the parameter set $I = {_{P_a}I_{3}} \cup
{_{P_a}I_{4}}$. The degrees in Table~A.15
show that the ${_{P_a}\chi_{3}}(k,l)$'s are the only irreducible characters
of $P_a$ of degree $q^2+1$, so 
${_{P_a}\chi_{3}}(k,l)^\alpha = {_{P_a}\chi_{3}}(k',l')$ for some $\{(k',l'),
\dots \} \in {_{P_a}I_{3}}$. We see from the class representatives in
\cite{HimstedtHuang2F4MaxParab} that $\alpha$ acts on the semisimple conjugacy
classes of $P_a$ like the $2$nd power map which implies that the values
of ${_{P_a}\chi_{3}}(k',l')$ and ${_{P_a}\chi_{3}}(2k,2l)$ on the semisimple
classes coincide. Then, the character values of ${_{P_a}\chi_{3}}(k,l)$
(see the character table A.5 in \cite{HimstedtHuang2F4MaxParab}) imply
that the values of ${_{P_a}\chi_{3}}(k',l')$ and ${_{P_a}\chi_{3}}(2k,2l)$ 
coincide on all classes, hence ${_{P_a}\chi_{3}}(k',l') = {_{P_a}\chi_{3}}(2k,2l)$ and
therefore ${_{P_a}\chi_{3}}(k,l)^\alpha = {_{P_a}\chi_{3}}(2k,2l)$. Similarly, 
${_{P_a}\chi_{4}}(k)^\alpha = {_{P_a}\chi_{4}}(2k)$. Hence, $x^\alpha = 2x$
for all $x\in I$. 
\smallskip

\noindent Let $I = {_GI_{27}} \cup {_GI_{33}}$. If $x = \{k ,-k\} \in I$, 
then $x \in C_I(H)$ if and only if $(2^t-1)k \equiv 0$ or 
$(2^t+1)k \equiv 0$. Let 
$
C_\pm := \left\{ \{k, -k\} \in C_I(H) \, | \, (2^t\pm1)k \equiv 0 \right\},
$
so that $C_I(H) = C_- \cup C_+$ and $C_- \cap C_+ = \emptyset$.
We claim 
\[
C_- = \left\{ \{k,-k\} \in {_GI_{27}} \, | \, k \text{   is a
  multiple of   } \frac{q^2-1}{2^{t}-1}\right\}.
\]
The inclusion $\supseteq$ is clear. Let $x = \{k, -k\} \in C_-$. If 
$x \in {_GI_{33}}$, then $(2^t-1) k \equiv 0$ mod $q^2+1$ and
Lemma~\ref{la:gcds}(ii) implies $k \equiv 0$, which is
impossible. Hence $x \in {_GI_{27}}$ and $(2^t-1)k \equiv 0$ mod $q^2-1$.
By Lemma~\ref{la:gcds}(i), $k$ is a multiple of $(q^2-1)/(2^t-1)$, 
proving the claim. Next, we claim 
$$C_+ = \left\{ \{k,-k\} \in {_GI_{33}} \,
  | \, k \text{   is a multiple of   } \frac{q^2+1}{2^t+1}\right\}.$$ 
The inclusion $\supseteq$ is clear.  
Let $x = \{k, -k\} \in C_+$. If $x \in {_GI_{27}}$, then $(2^t+1) k \equiv 0$
mod $q^2-1$ and Lemma~\ref{la:gcds}(iii) implies $k \equiv 0$, which is
impossible. Hence $x \in {_GI_{33}}$ and $(2^t+1)k \equiv 0$ mod $q^2+1$.
By Lemma~\ref{la:gcds}(iv), $k$ is a multiple of $(q^2+1)/(2^t+1)$ and
the claim holds.

\smallskip

\noindent Thus in all cases, $|C_I(H)| = |C_-| + |C_+| =
\frac{2^t-2}{2} + \frac{2^t-2}{2} = 2^t-2$. 

\medskip

\noindent Let $I = {_BI_{1}}$.  If $(k,l) \in I$, then $(k,l) \in C_I(H)$ if and only if $(2^t-1)k \equiv 0$ and $(2^t-1)l \equiv 0$ mod
$q^2-1$. By Lemma~\ref{la:gcds}(i), this is equivalent with $k, l$ are multiples of $\frac{q^2-1}{2^t-1}$. Thus, $|C_I(H)| = (2^t - 1)^2$.

\medskip

Now, let $I$ be one of the parameter sets ${_BI_{2}}$, ${_BI_{3}}$,
${_BI_{4}}$, ${_BI_{5}}$, ${_BI_{9}}$, ${_BI_{11}}$, ${_BI_{12}}$,
${_BI_{23}}$, ${_BI_{26}}$, ${_BI_{27}}$, ${_BI_{39}}$, ${_BI_{42}}$,
${_BI_{44}}$, ${_BI_{45}}$, ${_BI_{51}}$, ${_BI_{52}}$,
${_{P_a}I_{1}}$, ${_{P_a}I_{2}}$, ${_{P_a}I_{5}}$, ${_{P_a}I_{7}}$,
${_{P_a}I_{8}}$, ${_{P_a}I_{12}}$, ${_{P_a}I_{26}}$, ${_{P_a}I_{34}}$,
${_{P_a}I_{35}}$, ${_{P_b}I_{1}}$, ${_{P_b}I_{2}}$, ${_{P_b}I_{3}}$,
${_{P_b}I_{4}}$, ${_{P_b}I_{8}}$, ${_{P_b}I_{12}}$, ${_{P_b}I_{27}}$,
${_{P_b}I_{28}}$, ${_{P_b}I_{40}}$.  

\smallskip

If $k \in I$, then $k \in C_I(H)$ if and only if $(2^t-1)k \equiv 0$ mod
$q^2-1$. By Lemma~\ref{la:gcds}(i), we get 
$C_I(H) = \{k \in I \, | \, k \text{   is a multiple of   }
(q^2-1)/(2^t - 1)\}
$
and $|C_I(H)| = 2^t - 1$.

\medskip

\noindent Let $I = {_{P_a}I_{3}} \cup {_{P_a}I_{4}}$. 
First, we compute $|C_{_{P_a}I_{3}}(H)|$. Let
\[
U_i := \begin{cases} \left\{ \{(k, l) , (k', l')\} \in C_{_{P_a}I_{3}}(H) \, |
  \, 2^t k \equiv k , \, 2^t l \equiv l\right\} & \text{ if   } i = 1, \\ 
\left\{ \{(k, l) , (k', l')\} \in C_{_{P_a}I_{3}}(H) \, | \, 2^t k \equiv k', \, 2^t l \equiv l' \right\} & \text{ if   } i = 2, \\ 
\end{cases}
\]
where $k'=\frac{q}{\sqrt{2}}(k+l)$ and $l'=\frac{q}{\sqrt{2}}(k-l)$.

If $x = \{(k, l) , (\frac{q}{\sqrt{2}}(k+l), \frac{q}{\sqrt{2}}(k-l))\} \in 
{_{P_a}I_{3}}$, then $x \in U_1$ if and only if
$(2^t-1)k \equiv 0$ and $(2^t-1)l \equiv 0$ mod $q^2-1$. By
Lemma~\ref{la:gcds}(i), this is equivalent with $k$, $l$ are multiples of $\frac{q^2-1}{2^t-1}$. Thus, $|U_1| = (2^t-2)(2^t - 1)/2$.

\noindent If $x = \{(k, l) , (\frac{q}{\sqrt{2}}(k+l),
\frac{q}{\sqrt{2}}(k-l))\} \in {_{P_a}I_{3}}$, then $x \in U_2$ if and
only if $2^tk \equiv \frac{q}{\sqrt{2}}(k+l)$  and $2^tl \equiv
\frac{q}{\sqrt{2}}(k-l)$ mod $q^2-1$. Multiplying both congruences by
$\sqrt{2}q$, we get $\sqrt{2}q \cdot 2^t k\equiv k+l$ and $\sqrt{2}q
\cdot 2^t l\equiv k-l$ mod $q^2-1$. Then, $l \equiv (\sqrt{2}q \cdot
2^t - 1) k$ and $k \equiv (\sqrt{2}q \cdot 2^t + 1) l$ mod
$q^2-1$. Combining these two congruences, we get $k \equiv (\sqrt{2}q
\cdot 2^t - 1) (\sqrt{2}q \cdot 2^t + 1) k = (2q^2 \cdot 2^{2t}-1)k
\equiv (2 \cdot 2^{2t} - 1)k$. Then $2(2^{2t}-1)k \equiv 0$ mod
$q^2-1$. By Lemma~\ref{la:gcds}(i) and (iii), it follows that
$\gcd(2(2^t-1)(2^t+1), q^2-1) = \gcd(2^t-1 \, , \, q^2-1) =
2^t-1$. Thus $(2^t - 1)k \equiv 0$ mod $q^2-1$. Then $k \equiv 2^tk
\equiv \frac{q}{\sqrt{2}}(k+l)$. Multiplying both sides by $\sqrt{2}q$,
we get $\sqrt{2}qk \equiv k+l$. But then $l \equiv (\sqrt{2}q - 1)k$,
a~contradiction to the definition of ${_{P_a}I_{3}}$. Hence, 
$U_2 = \emptyset$. So
$|C_{_{P_a}I_{3}}(H)| = |U_1|+|U_2| = (2^t-2)(2^t-1)/2$.

\smallskip
\noindent Next, we calculate $|C_{_{P_a}I_{4}}(H)|$. If $x = \{k, q^2k\}
 \in {_{P_a}I_{4}}$, then $x \in C_{_{P_a}I_{4}}(H)$ if and only if 
$(2^t-1)k \equiv 0$ or $(2^t-q^2)k \equiv 0$ mod
$(q^2+1)(q^2-1)$. Suppose $(2^t-1)k \equiv~0$. By
Lemma~\ref{la:gcds}(i) and (ii), it follows that $\gcd(2^t-1, (q^2+1)(q^2-1)) = 
\gcd(2^t-1 \, , \, q^2-1) = 2^t-1$. Thus,
$(q^2+1)\cdot \frac{(q^2-1)}{2^t-1} \, | \, k$. But then
$(q^2+1) \, | \, k$, a contradiction to the definition of ${_{P_a}I_{4}}$.
So we have proved that $x \in C_{_{P_a}I_{4}}(H)$ if and only if
$(2^t-q^2)k \equiv 0$ mod $(q^2+1)(q^2-1)$. We claim
\[
C_{_{P_a}I_{4}}(H) = \left\{ \{k,q^2k\} \in {_{P_a}I_{4}} \, | \, k \text{   is a
  multiple of   } \frac{(q^2+1)(q^2-1)}{(2^t+1)(2^t-1)}\right\}.
\]
Let $k = \frac{(q^2+1)(q^2-1)}{(2^t+1)(2^t-1)} \cdot m$ 
for some $m \in \Z$. Because $t | 2n+1$ we have $2t \, | \, 2n+1-t$.  Then we 
get $(2^t+1)(2^t-1) \, | \, 2^{2n+1-t}-1$. Thus
$(2^{2n+1-t} - 1) k = \frac{2^{2n+1-t}-1}{(2^t+1)(2^t-1)} (q^2+1)(q^2-1)
\cdot m \equiv 0$ mod $(q^2+1)(q^2-1)$. So $(2^t-q^2) k \equiv 0$ mod
$(q^2+1)(q^2-1)$ and $\{k,q^2k\} \in C_{_{P_a}I_{4}}(H)$.

Conversely, suppose $\{k,q^2k\} \in C_{_{P_a}I_{4}}(H)$. Then
$(2^t-q^2) k \equiv 0$ mod $(q^2+1)(q^2-1)$. Hence
$(2^t+1)k \equiv 0$ mod $q^2+1$ and $(2^t-1)k \equiv 0$ mod $q^2-1$. By
Lemma~\ref{la:gcds}(i) and (iv), this is equivalent with 
$\frac{q^2+1}{2^t+1} \, | \, k$ and $\frac{q^2-1}{2^t-1} \, | \, k$. 
Since
$\frac{q^2+1}{2^t+1} \, | \, q^2+1$ and $\frac{q^2-1}{2^t-1} \, | \,
q^2-1$ and since $\gcd(q^2+1, q^2-1)=1$,
we have $\gcd\left(\frac{q^2+1}{2^t+1},\frac{q^2-1}{2^t-1}\right) = 1$.
Therefore $\frac{(q^2+1)(q^2-1)}{(2^t+1)(2^t-1)} \, | \, k$, and 
the claim holds. So by the definition of ${_{P_a}I_{4}}$, we get 
$|C_{_{P_a}I_{4}}(H)| = 2^t(2^t-1)/2$.

\smallskip

\noindent So, $|C_I(H)| =
|C_{_{P_a}I_{3}}(H)|+|C_{_{P_a}I_{4}}(H)| = (2^t-1)^2$.

\bigskip

\noindent Let $I = {_{P_a}I_{32}} \cup {_{P_a}I_{33}}$. The computation
of $|C_I(H)|$ is analogous to the case $I = {_{G}I_{27}} \cup
{_{G}I_{33}}$. One gets $|C_I(H)| = 2^t - 1$. 

\medskip

\noindent Let $I = {_{P_b}I_{5}} \cup {_{P_b}I_{6}} \cup {_{P_b}I_{7}}$. 
First, we compute $|C_{_{P_b}I_{5}}(H)|$. Let
\[
U_i := \begin{cases} \left\{ \{(k, l) , (k, -l)\} \in C_{_{P_b}I_{5}}(H) \, |
  \, 2^t k \equiv k , \, 2^t l \equiv l\right\} & \text{ if   } i = 1, \\ 
\left\{ \{(k, l) , (k, -l)\} \in C_{_{P_b}I_{5}}(H) \, | \, 2^t k \equiv k , \, 
  2^t l \equiv -l \right\} & \text{ if   } i = 2. \\ 
\end{cases}
\]
If $x = \{(k, l) , (k, -l)\} \in {_{P_b}I_{5}}$, then $x \in U_1$ 
if and only if $(2^t-1)k \equiv 0$ and $(2^t-1)l \equiv 0$ mod 
$q^2-1$. By Lemma~\ref{la:gcds}(i), this is equivalent with $k$, $l$ are multiples of $\frac{q^2-1}{2^t-1}$. Thus, $|U_1| = (2^t-1)(2^t - 2)/2$.

\smallskip

\noindent If $x = \{(k, l) , (k, -l)\} \in {_{P_b}I_{5}}$, then $x \in U_2$ 
if and only if $(2^t-1)k \equiv 0$ and $(2^t+1)l \equiv 0$ mod 
$q^2-1$. By Lemma~\ref{la:gcds}(iii), the second congruence is equivalent with $l \equiv 0$ mod $q^2-1$, a contradiction to the definition of~${_{P_b}I_{5}}$. Hence
$U_2 = \emptyset$. So
$|C_{{_{P_b}I_{5}}}(H)| = |U_1|+|U_2| = (2^t-1)(2^t-2)/2$.

\smallskip

\noindent Next, we calculate $|C_{_{P_b}I_{6}\, \cup\, {_{P_b}I_{7}}}(H)|$. 
Let $\phi_8' := q^2+\sqrt{2}q+1$ and $\phi_8'' :=
q^2-\sqrt{2}q+1$. The set $J_1 := \left\{\{(k,l), (k,-l), 
(k,q^2l), (k,-q^2l)\} | (k,l) \in \Z_{q^2-1} \times
\Z_{q^2+\sqrt{2}q+1}, l\not= 0\right\}$ becomes an $H$-set via
$x^\alpha := 2 x$ for all $x \in J$, and we have 
$J_1 \simeq {_{P_b}I_{6}}$ as $H$-sets; note that
${_{P_b}\chi_{6}}(k)$ is obtained via inflation from a regular
semisimple character of the Levi complement 
$L_b \cong \Z_{q^2-1} \times \Sz(q^2)$. 
In a similar way, we have an isomorphism of $H$-sets 
$${_{P_b}I_{7}} \simeq \left\{\{(k,l), (k,-l), 
(k,q^2l), (k,-q^2l)\} | (k,l) \in \Z_{q^2-1} \times
\Z_{q^2-\sqrt{2}q+1}, l\not= 0\right\} =: J_2.$$ 
Hence, ${_{P_b}I_{6}} \cup {_{P_b}I_{7}} \simeq J_1 \cup J_2$
(disjoint union) as $H$-sets, so that we can identify ${_{P_b}I_{6}}
\cup {_{P_b}I_{7}} =  J_1 \cup J_2$. 

In a first step, we compute $|C_{J_1}(H)|$. Let
\begin{eqnarray*}
U_{\pm1} & := & \left\{ \{(k,l), (k,-l), (k,q^2l), (k,-q^2l)\} \in C_{J_1}(H) \, |
  \, 2^t k \equiv k\, ,\, 2^t l \equiv \pm l \right\} , \\
U_{\pm2} & := & \left\{ \{(k,l), (k,-l), (k,q^2l), (k,-q^2l)\} \in C_{J_1}(H) \, |
  \, 2^t k \equiv k\, ,\, 2^t l \equiv \pm q^2l \right\} ,
\end{eqnarray*}
where the congruences are mod $q^2-1$ and $q^2+\sqrt{2}q+1$,
respectively. Then $C_{J_1}(H) = U_1 \cup U_{-1} \cup U_{2}
\cup U_{-2}$. We claim that $U_{\pm 1} = \emptyset$. 
Suppose 
$
\{(k,l), \dots\} \in U_1.
$
Then $(2^t-1)k \equiv 0$ mod $q^2-1$ and $(2^t-1)l \equiv 0$ mod
$q^2+\sqrt{2}q+1$. Since $2^t-1\, | \, q^2-1\, | \, q^4-1$ and since
$\gcd(q^4-1, (q^2+\sqrt{2}q+1) \cdot (q^2-\sqrt{2}q+1)) = \gcd(q^4-1,
q^4+1)=1$, the second congruence is equivalent with 
$l\equiv 0$ mod $q^2+\sqrt{2}q+1$, contradicting the definition of 
$J_1$. So $U_1 = \emptyset$. Suppose $\{(k,l), \dots\} \in U_{-1}$. 
Then $(2^t-1)k \equiv 0$ mod $q^2-1$ and $(2^t+1)l \equiv 0$ mod
$q^2+\sqrt{2}q+1$. Since $2^t+1\, | \, q^2+1\, | \, 
q^4-1$ and since $\gcd(q^4-1, q^4+1)=1$,  the second congruence
is equivalent with $l\equiv 0$ mod $q^2+\sqrt{2}q+1$, contradicting the
definition of $J_1$. So $U_{-1} = \emptyset$. 

Suppose $\{(k,l), \dots\} \in U_{2}$. Then $2^tk \equiv k$ mod $q^2-1$ and
$2^tl \equiv q^2l$ mod $q^2+\sqrt{2}q+1$, and hence $(2^t-1)k \equiv 0$ mod
$q^2-1$ and $(2^{2n+1-t}-1)l \equiv 0$ mod $q^2+\sqrt{2}q+1$. By
Lemmas~\ref{la:gcds}(i), \ref{la:gcds1} and the definition of $J_1$,
we get  
\[
|U_2| = \begin{cases} \frac 1 4 (2^t-1)(2^t\pm2^{(t+1)/2})
    & \text{   if   } 4t\, | \, 2n+1-t, \\
0 & \text{   if   } 4t \nmid \, 2n+1-t.
\end{cases}
\]
Suppose $\{(k,l), (k,-l), (k,q^2l), (k,-q^2l)\} \in U_{-2}$. Then
$2^tk \equiv k$ mod $q^2-1$ and $2^tl \equiv -q^2l$ mod
$q^2+\sqrt{2}q+1$, and hence $(2^t-1)k \equiv 0$ mod $q^2-1$
and$(2^{2n+1-t}+1)k \equiv 0$ mod $q^2+\sqrt{2}q+1$. By
Lemma~\ref{la:gcds}(i), \ref{la:gcds1} and the definition of $J_1$, we
get  
\[
|U_{-2}| = \begin{cases} 0
    & \text{   if   } 4t\, | \, 2n+1-t, \\
\frac 1 4 (2^t-1)(2^t\pm2^{(t+1)/2}) & \text{   if   } 4t \nmid\, 2n+1-t.
\end{cases}
\]
So
$$|C_{J_1}(H)| = |U_1| + |U_{-1}| + |U_2| + |U_{-2}|
= \frac 1 4 (2^t-1)(2^t\pm2^{(t+1)/2}).$$
Analogously, we obtain
$$|C_{J_2}(H)| = \frac 1 4 (2^t-1)(2^t\mp2^{(t+1)/2}).$$
So
$$|C_{_{P_b}I_6\, \cup _{P_b}I_7}(H)| = |C_{J_1 \cup J_2}(H)| =
|C_{J_1}(H)|+|C_{J_2}(H)| = 2^t(2^t-1)/2.$$
Hence we get 
$$|C_I(H)| = |C_{_{P_b}I_5}(H)| + |C_{_{P_b}I_6\, \cup _{P_b}I_7}(H)|
= \frac{(2^t-1)(2^t-2)}{2} + \frac{2^t(2^t-1)}{2} = (2^t-1)^2.$$ 

\bigskip

\noindent Let $I \in \{{_{G}I_{23}} \cup {_{G}I_{43}} \cup
          {_{G}I_{47}} \, , {_{G}I_{24}} \cup {_{G}I_{44}} \cup
          {_{G}I_{48}} \, , {_{G}I_{25}} \cup {_{G}I_{45}} \cup
          {_{G}I_{49}} \, , {_{P_b}I_{16}} \cup {_{P_b}I_{21}} \cup
          {_{P_b}I_{26}} \, , {_{P_b}I_{47}} \cup {_{P_b}I_{48}} \cup
          {_{P_b}I_{49}} \, , {_{P_b}I_{54}} \cup {_{P_b}I_{55}} \cup
          {_{P_b}I_{56}}\}$.  Then these (disjoint) unions of
          parameter sets are isomorphic $H$-sets, so that we can
          assume $I = {_{P_b}I_{47}} \cup {_{P_b}I_{48}} \cup
          {_{P_b}I_{49}}$. First, we compute
          $|C_{_{P_b}I_{47}}(H)|$. If $x = \{k ,-k\} \in I$, then $x
          \in C_{_{P_b}I_{47}}(H)$ if and only if $(2^t-1)k \equiv 0$
          or $(2^t+1)k \equiv 0$ mod $q^2-1$. Let 
\[
C_\pm := \left\{ \{k, -k\} \in C_{_{P_b}I_{47}}(H) \, | \, (2^t\pm1)k \equiv 0 \right\},
\]
so that $C_I(H) = C_- \cup C_+$ and $C_- \cap C_+ = \emptyset$. We claim that $C_+ = \emptyset$. Let $x = \{k, -k\} \in C_+$, then $(2^t+1) k \equiv 0$ mod $q^2-1$ and Lemma~\ref{la:gcds}(iii) implies $k \equiv 0$, which is impossible. Hence $x \in C_-$ and $(2^t-1)k \equiv 0$ mod $q^2-1$.
By Lemma~\ref{la:gcds}(i), $k$ is a multiple of $(q^2-1)/(2^t-1)$ and $|C_I(H)| = |C_-| + |C_+| = \frac{2^t-2}{2}$.
  
\medskip

\noindent Next, we compute $|C_{{_{P_b}I_{48}}\, \cup
  {_{P_b}I_{49}}}(H)|$. We start with ${_{P_b}I_{48}}$. Let
\begin{eqnarray*}
U_{\pm1} & := & \left\{ \{k, -k, q^2k, -q^2k\} \in C_{{_{P_b}I_{48}}}(H) \, |
  \, 2^t k \equiv \pm k \right\} , \\
U_{\pm2} & := & \left\{ \{k, -k, q^2k, -q^2k\} \in C_{{_{P_b}I_{48}}}(H) \, |
  \, 2^t k \equiv \pm q^2k \right\} ,
\end{eqnarray*}
where the congruences are mod $q^2+\sqrt{2}q+1$. 
Then $C_{{_{P_b}I_{48}}}(H) = U_1 \cup U_{-1} \cup U_{2} \cup U_{-2}$. We claim
that $U_{\pm 1} = \emptyset$. Suppose 
$\{k, -k, q^2k, -q^2k\} \in U_1$. Then $(2^t-1)k \equiv 0$. 
Since $2^t-1\, | \, q^2-1\, | \, q^4-1$ and
since $\gcd(q^4-1, (q^2+\sqrt{2}q+1) \cdot (q^2-\sqrt{2}q+1)) =
\gcd(q^4-1, q^4+1)=1$, we have $k\equiv 0$ mod $q^2+\sqrt{2}q+1$,
contradicting the definition of ${_{P_b}I_{48}}$. So $U_1 = \emptyset$. 
Suppose $\{k, -k, q^2k, -q^2k\} \in U_{-1}$. Then $(2^t+1)k \equiv 0$
mod $q^2+\sqrt{2}q+1$. Since $2^t+1\, | \, q^2+1\, | \, q^4-1$ and
since $\gcd(q^4-1, q^4+1)=1$, it follows that $k\equiv 0$ mod 
$q^2+\sqrt{2}q+1$, contradicting the definition of ${_{P_b}I_{48}}$. 
So $U_{-1} = \emptyset$. Suppose $\{k, -k, q^2k, -q^2k\} \in U_{2}$.
Then $2^tk \equiv q^2k$ mod $q^2+\sqrt{2}q+1$ and hence 
$(2^{2n+1-t}-1)k \equiv 0$ mod $q^2+\sqrt{2}q+1$. By 
Lemma~\ref{la:gcds1} and the definition of ${_{P_b}I_{48}}$, we get 
\[
|U_2| = \begin{cases} \frac 1 4 (2^t\pm2^{(t+1)/2})
    & \text{   if   } 4t\, | \, 2n+1-t , \\
0 & \text{   if   } 4t \nmid \, 2n+1-t.
\end{cases}
\]
Suppose $\{k, -k, q^2k, -q^2k\} \in U_{-2}$. Then $2^tk \equiv -q^2k$
mod $q^2+\sqrt{2}q+1$ and hence $(2^{2n+1-t}+1)k \equiv 0$ mod 
$q^2+\sqrt{2}q+1$. By Lemma~\ref{la:gcds1} and the definition of
${_{P_b}I_{48}}$, we get 
\[
|U_{-2}| = \begin{cases} 0
    & \text{   if   } 4t\, | \, 2n+1-t, \\
\frac 1 4 (2^t\pm2^{(t+1)/2}) & \text{   if   } 4t \nmid \, 2n+1-t.
\end{cases}
\]
So
$$|C_{{_{P_b}I_{48}}}(H)| = |U_1| + |U_{-1}| + |U_2| + |U_{-2}| =
\frac 1 4 (2^t\pm2^{(t+1)/2}).$$
Analogously, we obtain 
$$|C_{{_{P_b}I_{49}}}(H)| = \frac 1 4 (2^t\mp2^{(t+1)/2}).$$
So we have $|C_{{_{P_b}I_{48}}\, \cup {_{P_b}I_{49}}}(H)| =
|C_{{_{P_b}I_{48}}}| + |C_{{_{P_b}I_{49}}}| = 2^t/2$. Hence, we get 
$|C_I(H)| = |C_{{_{P_b}I}_{47}}(H)| + |C_{{_{P_b}I}_{48}\, \cup
  {_{P_b}I}_{49}}(H)| = \frac{2^t-2}{2} + \frac{2^t}{2} = 2^t-1.$
\hfill $\Box$


\begin{proposition} \label{prop:idfixnr}
Let $t \, | \, 2n+1$, $H = \la \alpha^t \ra \le O$ and let
$(I,J) \in \{
({_GI_{15}} \cup {_GI_{16}}, {_{P_a}I_{18}} \cup {_{P_a}I_{19}}), \, 
({_BI_{24}} \cup {_BI_{25}}, {_{P_a}I_{27}} \cup {_{P_a}I_{28}}), \, 
({_BI_{30}} \cup {_BI_{31}} \cup {_BI_{32}} \cup {_BI_{33}} \cup
{_BI_{35}} \cup {_BI_{36}}, {_{P_b}I_{31}} \cup {_{P_b}I_{32}} \cup
{_{P_b}I_{33}} \cup {_{P_b}I_{34}} \cup {_{P_b}I_{36}} \cup
{_{P_b}I_{37}}), \, 
({_BI_{34}} \cup {_BI_{37}}, {_{P_b}I_{35}} \cup {_{P_b}I_{38}}), \, 
({_BI_{38}} \cup {_BI_{40}} \cup {_BI_{41}}, {_{P_b}I_{39}} \cup
{_{P_b}I_{41}} \cup {_{P_b}I_{42}}), \,  
({_BI_{46}} \cup {_BI_{47}} \cup {_BI_{48}} \cup {_BI_{49}},
{_{P_a}I_{36}} \cup {_{P_a}I_{37}} \cup {_{P_a}I_{38}} \cup
{_{P_a}I_{39}}), \, 
({_BI_{50}}, {_{P_a}I_{40}})\}$. Then $|C_I(H)| = |C_J(H)|$.
\end{proposition}

\noindent {\bf Proof:} Suppose $(I, J) = ({_GI_{15}} \cup {_GI_{16}}, 
{_{P_a}I_{18}} \cup {_{P_a}I_{19}})$. Considering degrees we see that
the sets $\{{_G\chi_{15}}, {_G\chi_{16}}\}$ and 
$\{{_{P_a}\chi_{18}}, {_{P_a}\chi_{19}}\}$ are invariant under the
action of $\alpha$ (note that ${_{P_a}\chi_{20}}$ is real-valued and 
${_{P_a}\chi_{18}}$, ${_{P_a}\chi_{19}}$ are not). The characters
${_G\chi_{15}}$, ${_G\chi_{16}}$ coincide on all classes of $G$ except
for the classes $c_{5,3}$ and $c_{5,4}$. Similarly,
${_{P_a}\chi_{18}}, {_{P_a}\chi_{19}}$ coincide on all conjugacy
classes of $P_a$ except for $c_{8,3}$ and $c_{8,4}$. From the class
representatives in \cite[Table~A.2]{HimstedtHuang2F4Borel} and
\cite[Table~A.3]{HimstedtHuang2F4MaxParab} we get that $\alpha$ swaps 
${_G\chi_{15}}$ and ${_G\chi_{16}}$ if and only if it swaps 
${_{P_a}\chi_{18}}$ and ${_{P_a}\chi_{19}}$. So, $|C_I(H)| =
|C_J(H)|$.

Next, suppose $(I, J) = ({_BI_{50}}, {_{P_a}I_{40}})$. 
By construction, ${_{P_a}\chi}_{40}(k) =
{_B\chi}_{50}(k)^{P_a}$ is induced from the $\alpha$-stable
Borel subgroup $B$ for all $k \in {_BI_{50}} = {_{P_a}I_{40}}$ 
(see~\cite[Table~3]{HimstedtHuang2F4MaxParab}) and 
induction of characters induces a bijection from 
$\{{_B\chi}_{50}(k) \, | \, k \in {_BI_{50}}\}$ onto
$\{{_{P_a}\chi}_{40}(k) \, | \, k \in {_{P_a}I_{40}}\}$ mapping
${_B\chi}_{50}(k)$ to ${_{P_a}\chi}_{40}(k)$.
We have 
\[
{_{P_a}\chi}_{40}(k^\alpha) = {_{P_a}\chi}_{40}(k)^\alpha =
\left({_B\chi}_{50}(k)^{P_a}\right)^\alpha =
\left({_B\chi}_{50}(k)^\alpha\right)^{P_a} =
\left({_B\chi}_{50}(k^\alpha)\right)^{P_a}. 
\]
So the above-mentioned bijection is an isomorphism of $H$-sets.
Hence, ${_BI_{50}} \simeq {_{P_a}I_{40}}$ as $H$-sets and 
$|C_{_BI_{50}}(H)| = |C_{_{P_a}I_{40}}(H)|$. The remaining pairs of
parameter sets can be treated analogously, see the construction of 
the characters ${_{P_a}\chi}_{28}(k)$, ${_{P_a}\chi}_{38}(k)$,
${_{P_a}\chi}_{39}(k)$, ${_{P_b}\chi}_{36}(k)$,
${_{P_b}\chi}_{37}(k)$, ${_{P_b}\chi}_{38}(k)$ and
${_{P_b}\chi}_{42}(k)$ in
\cite{HimstedtHuang2F4MaxParab}. \hfill $\Box$   



\section{Dade's Invariant Conjecture for \texorpdfstring{${^2F}_4(q^2)$}{2F4}}
\label{sec:proof3d4}

In this section, we prove Dade's invariant conjecture for $G =
{^2F}_4(2^{2n+1})$, $n>0$, in the defining characteristic.
As in the previous section, let $O = \Out(G) = \langle \alpha \rangle$, where
$\alpha$ is a field automorphism of (odd) order~$2n+1$. We fix a Borel
subgroup $B$ and maximal parabolic subgroups $P_a$ and $P_b$ of~$G$
containing $B$ as in~\cite{HimstedtHuang2F4Borel} and
\cite{HimstedtHuang2F4MaxParab}. In particular, we may assume that
$\alpha$ stabilizes $B$, $P_a$ and $P_b$.

By the remarks on p.~152 in \cite{HumphreysDefGrps}, $G$ has only two
$2$-blocks, the principal block $B_0$ and one defect-$0$-block (corresponding
to the Steinberg character). Hence we have to verify Dade's conjecture only
for the principal block $B_0$.

According to a corollary of the Borel-Tits theorem~\cite{BW}, the
normalizers of radical $2$-subgroups are parabolic subgroups. The
radical $2$-chains of $G$ (up to $G$-conjugacy) are given in Table
\ref{tab:chains}. 


\begin{table}[ht]
\caption[]{Radical $2$-chains of $G$.} \label{tab:chains}
\begin{center}
\begin{tabular}{|l|l|l|l|} \hline
\rule{0cm}{0.4cm}
\hspace{-0.12cm}C & & $N_G(C)$ & $N_A(C)$
\rule[-0.1cm]{0cm}{0.3cm}\\
\hline \cline{1-3} \hline
\rule{0cm}{0.4cm}
\hspace{-0.12cm}$C_1$ & $\{1\}$ & $G$ & $A$
\rule[-0.2cm]{0cm}{0.1cm}\\
\rule{-0.12cm}{0.1cm}
$C_2$ & $\{1\} < O_p(P_a)$ & $P_a$ & $P_a \rtimes \langle \alpha \rangle$
\rule[-0.2cm]{0cm}{0.1cm}\\
\rule{-0.12cm}{0.1cm}
$C_3$ & $\{1\} < O_p(P_a) < O_p(B)$ & $B$ & $B \rtimes \langle \alpha \rangle$
\rule[-0.2cm]{0cm}{0.1cm}\\
\rule{-0.12cm}{0.1cm}
$C_4$ & $\{1\} < O_p(P_b)$ & $P_b$ & $P_b \rtimes \langle \alpha \rangle$
\rule[-0.2cm]{0cm}{0.1cm}\\
\rule{-0.12cm}{0.1cm}
$C_5$ & $\{1\} < O_p(P_b) < O_p(B)$ & $B$ & $B \rtimes \langle \alpha \rangle$
\rule[-0.2cm]{0cm}{0.1cm}\\
\rule{-0.12cm}{0.1cm}
$C_6$ & $\{1\} < O_p(B)$ & $B$ & $B \rtimes \langle \alpha \rangle$
\rule[-0.2cm]{0cm}{0.3cm}\\
\hline
\end{tabular}
\end{center}
\end{table}


\noindent Since $C_5$ and $C_6$ have the same normalizers
$N_G(C_5) = N_G(C_6)$ and $N_A(C_5) = N_A(C_6)$, it follows that
\[
k(N_G(C_5), B_0, d, u) = k(N_G(C_6), B_0, d, u)
\]
for all $d \in \N$ and $u \mid 2n+1$. Thus the contribution of $C_5$ and 
$C_6$ in the alternating sum of Dade's invariant conjecture is zero. So, Dade's
invariant conjecture for $G$ is equivalent with
\begin{equation} \label{eq:altsum}
k(G, B_0, d, u) + k(B, B_0, d, u)= k(P_a, B_0, d, u) + k(P_b, B_0, d, u)
\end{equation}
for all $d \in \N$ and $u\mid 2n+1$.

\begin{theorem}
\label{thmp}
Let $\widetilde{B}$ be a $2$-block of $G ={^2F}_4(2^{2n+1})$, $n>0$,
with a positive defect. Then $\widetilde{B}$ satisfies Dade's
invariant conjecture. 
\end{theorem}

\noindent {\bf Proof:} By the proceeding remarks, we can assume
$\widetilde{B} = B_0$. Suppose $u \mid 2n+1$ and set $t:=\frac{2n+1}{u}$ and 
$H :=\langle \alpha^t \rangle$. Let $S \in \{G, B, P_a, P_b\}$. 
By the character tables in \cite{HimstedtHuang2F4Borel},
\cite{HimstedtHuang2F4MaxParab} and CHEVIE, we 
have $k(S, B_0, d, u) = 0$ when $d \not \in \{11n+6, 14n+7, 14n+8, 15n+8, 16n+8, 17n+9, 18n+9, 20n+10, 20n+11, 20n+12, 21n+11, 22n+11, 23n+12, 24n+12\}$. 

\smallskip

(a) If $d = 11n+6$, then we have $k(G, B_0, d, u) =
\sum_{j \in \{19,20\}}|C_{_GI_{j}}(H)| = 2$ and $k(P_b, B_0, d, u) = \sum_{j \in \{46,53\}}|C_{_{P_b}I_{j}}(H)| = 2$ as well as
$k(B, B_0, d, u) = k(P_a, B_0, d, u) = 0$ by Tables~A.2 and~A.16. Thus
(\ref{eq:altsum}) holds in this case.  

\smallskip

(b) If $d = 14n+7$, then we have 
\[
k(G, B_0, d, u) = |C_{_GI_{18}}(H)| = 1 \text{   and   } k(P_a, B_0,
d, u) = |C_{_{P_a}I_{31}}(H)| = 1
\]
by Tables~A.3 and A.16. Thus~(\ref{eq:altsum}) holds in this case.   

\smallskip

(c) If $d = 14n+8$, then we have $k(G, B_0, d, u) =
\sum_{j \in \{55,56,57,58\}}|C_{_BI_{j}}(H)| = 4$ and $k(P_b, B_0, d, u) = \sum_{j \in \{44,45,51,52\}}|C_{_{P_b}I_{j}}(H)| = 4$ by Tables~A.4 and~A.16. Thus~(\ref{eq:altsum}) holds in this case.  

\smallskip
(d) If $d = 15n+8$, then we have 
$k(B, B_0, d, u) =\sum_{j \in \{51,52,53,54\}} |C_{_BI_{j}}(H)| = 2\cdot 2^t$ and
$k(P_b, B_0, d, u) =\sum_{j \in \{43,47,48,49,50,54,55,56\}}
|C_{_{P_b}I_{j}}(H)| = 2\cdot 2^t$ by Tables~A.5 and~A.16. Thus
(\ref{eq:altsum}) holds in this case. 

\smallskip

(e) If $d = 16n+8$, then Table~A.6 and Proposition~\ref{prop:idfixnr} imply, that (\ref{eq:altsum}) is equivalent to
\[
\sum_{j \in \{42,43\}} |C_{_BI_{j}}(H)| =
\sum_{j \in \{30,32,33\}} |C_{_{P_a}I_{j}}(H)|. 
\]
By Table~A.16, the sums on both sides of the above
equation are equal. Thus (\ref{eq:altsum}) also holds in this case. 

\smallskip

(f) If $d = 17n+9$, then Table~A.7 and Proposition~\ref{prop:idfixnr}
imply, that (\ref{eq:altsum}) is equivalent to 
\[
\sum_{j \in J_B} |C_{_BI_{j}}(H)| =
\sum_{j \in J_{P_a}} |C_{_{P_a}I_{j}}(H)| 
\]
with $J_B := \{44,45\}$ and $J_{P_a} := \{34,35\}$. By Table~A.16, the
sums on both sides of the above equation are equal. Thus
(\ref{eq:altsum}) also holds in this case. 

\smallskip

(g) If $d = 18n+9$, then Table~A.8 and
Proposition~\ref{prop:idfixnr} imply, that (\ref{eq:altsum}) is equivalent to
\[
|C_{_GI_{31}}(H)| + |C_{_BI_{39}}(H)| =
|C_{_{P_a}I_{29}}(H)| + |C_{_{P_b}I_{40}}(H)|.
\]
By Table~A.16, the sums on both sides of the above equation are
equal. Thus (\ref{eq:altsum}) also holds in this case. 

\smallskip

(h) If $d = 20n+10$, then Table~A.9 and Proposition~\ref{prop:idfixnr}
imply, that (\ref{eq:altsum}) is equivalent to
\[
\sum_{j \in J_G} |C_{_GI_{j}}(H)| + \sum_{j \in J_B} |C_{_BI_{j}}(H)| =
\sum_{j \in J_{P_a}} |C_{_{P_a}I_{j}}(H)| + \sum_{j \in J_{P_b}} |C_{_{P_b}I_{j}}(H)|
\]
where the sums are over $J_G = \{17,25,45,49\}$,
$J_B = \{11,12,23\}$, $J_{P_a} = \{12,20,26\}$,
$J_{P_b} = \{4,16,21,26\}$, respectively.
By Table~A.16, we have
\[
k(G, B_0, d, u) + k(B, B_0, d, u) = \sum_{j \in J_G} |C_{_GI_{j}}(H)| + \sum_{j \in J_B} |C_{_BI
_{j}}(H)| = 4\cdot 2^t - 3
\]
and
\[
k(P_a, B_0, d, u) + k(P_b, B_0, d, u) = \sum_{j \in J_{P_a}}
|C_{_{P_a}I_{j}}(H)| + \sum_{j \in J_{P_b}} |C_{_{P_b}I_{j}}(H)| =
4\cdot 2^t - 3. 
\]
Thus (\ref{eq:altsum}) also holds in this case.

\smallskip

(i) If $d = 20n+11$, then Table~A.10 and
Proposition~\ref{prop:idfixnr} imply, that (\ref{eq:altsum}) is equivalent to
\[
\sum_{j \in J_G} |C_{_GI_{j}}(H)| + \sum_{j \in J_B} |C_{_BI_{j}}(H)| =
\sum_{j \in J_{P_a}} |C_{_{P_a}I_{j}}(H)| + \sum_{j \in J_{P_b}} |C_{_{P_b}I_{j}}(H)|
\]
where the sums are over $J_G := \{7,10\}$, $J_B :=
\{13,14\}$, $J_{P_a} := \{13,25\}$ and $J_{P_b} :=
\{14,15\}$, respectively. By Table~A.16, we have 
\[
k(G, B_0, d, u) + k(B, B_0, d, u) = \sum_{j \in J_G} |C_{_GI_{j}}(H)|
+ \sum_{j \in J_B} |C_{_BI_{j}}(H)| = 4 
\]
and
\[
k(P_a, B_0, d, u) + k(P_b, B_0, d, u) = \sum_{j \in J_{P_a}}
|C_{_{P_a}I_{j}}(H)| + \sum_{j \in J_{P_b}} |C_{_{P_b}I_{j}}(H)| = 4. 
\]
Thus (\ref{eq:altsum}) also holds in this case.

\smallskip

(j) If $d = 20n+12$, then Table~A.11 implies, that (\ref{eq:altsum}) is
equivalent to
\[
\sum_{j \in J_G} |C_{_GI_{j}}(H)| + \sum_{j \in J_B} |C_{_BI_{j}}(H)| =
\sum_{j \in J_{P_a}} |C_{_{P_a}I_{j}}(H)| + \sum_{j \in J_{P_b}} |C_{_{P_b}I_{j}}(H)|
\]
with $J_G := \{5,6,8,9,11,12,13,14\}$, $J_B :=
\{15,16,17,18,19,20,21,22\}$, $J_{P_a} :=
\{14,15,16,17,21,22,23,24\}$, $J_{P_b} :=
\{17,18,19,20,22,23,24,25\}$. By Table~A.16, we have
\[
k(G, B_0, d, u) + k(B, B_0, d, u) = \sum_{j \in J_G} |C_{_GI_{j}}(H)| + \sum_{j \in J_B} |C_{_BI_{j}}(H)| = 16
\]
and
\[
k(P_a, B_0, d, u) + k(P_b, B_0, d, u) = \sum_{j \in J_{P_a}} |C_{_{P_a}I_{j}}(H)| + \sum_{j \in J_{P_b}} |C_{_{P_b}I_{j}}(H)| = 16.
\]
Thus (\ref{eq:altsum}) also holds in this case.

\smallskip

(k) If $d = 21n+11$, then Table~A.12 implies, that (\ref{eq:altsum}) is equivalent to
\[
\sum_{j \in J_B} |C_{_BI_{j}}(H)| =
\sum_{j \in J_{P_b}} |C_{_{P_b}I_{j}}(H)| 
\]
with $J_B := \{26,27,28,29\}$ and $J_{P_b} := \{27,28,29,30\}$. By Table~A.16, the sums on both sides of the above equation are equal. Thus (\ref{eq:altsum}) also holds in this case.

\smallskip

(l) If $d = 22n+11$, then Table~A.13 implies, that (\ref{eq:altsum}) is
equivalent to
\[
\sum_{j \in J_G} |C_{_GI_{j}}(H)| + \sum_{j \in J_B} |C_{_BI_{j}}(H)| = \sum_{j \in J_{P_a}} |C_{_{P_a}I_{j}}(H)| + \sum_{j \in J_{P_b}} |C_{_{P_b}I_{j}}(H)|
\]
with $J_G := \{4,27,30,33\}$, $J_B := \{9,10\}$, $J_{P_a} := \{2,11\}$ and 
$J_{P_b} := \{12,13\}$.
By Table~A.16, we have
\[
k(G, B_0, d, u) + k(B, B_0, d, u) = \sum_{j \in J_G} |C_{_GI_{j}}(H)| + \sum_{j \in J_B} |C_{_BI_{j}}(H)| = 2\cdot 2^t
\]
and
\[
k(P_a, B_0, d, u) + k(P_b, B_0, d, u) = \sum_{j \in J_{P_a}} |C_{_{P_a}I_{j}}(H)| + \sum_{j \in J_{P_b}} |C_{_{P_b}I_{j}}(H)| = 2\cdot 2^t.
\]
Thus (\ref{eq:altsum}) also holds in this case.

\smallskip

(m) If $d = 23n+12$, then Table~A.14 implies, that (\ref{eq:altsum}) is
equivalent to
\[
\sum_{j \in J_G} |C_{_GI_{j}}(H)| + \sum_{j \in J_B} |C_{_BI_{j}}(H)| = \sum_{j \in J_{P_a}} |C_{_{P_a}I_{j}}(H)| + \sum_{j \in J_{P_b}} |C_{_{P_b}I_{j}}(H)|
\]
with $J_G := \{2,3,23,24,43,44,47,48\}$, $J_B := \{3,4,6,7\}$, $J_{P_a} := \{7,8,9,10\}$ and 
$J_{P_b} := \{2,3,9,10\}$.
By Table~A.16, we have
\[
k(G, B_0, d, u) + k(B, B_0, d, u) = \sum_{j \in J_G} |C_{_GI_{j}}(H)| + \sum_{j \in J_B} |C_{_BI_{j}}(H)| = 4\cdot 2^t
\]
and
\[
k(P_a, B_0, d, u) + k(P_b, B_0, d, u) = \sum_{j \in J_{P_a}} |C_{_{P_a}I_{j}}(H)| + \sum_{j \in J_{P_b}} |C_{_{P_b}I_{j}}(H)| = 4\cdot 2^t.\]
Thus (\ref{eq:altsum}) also holds in this case.

\smallskip

(n) If $d = 24n+12$, then Table~A.15 implies, that (\ref{eq:altsum}) is
equivalent to
\[
|C_{_GI^{ss}}(H)| + \sum_{j \in J_B} |C_{_BI_{j}}(H)| = \sum_{j \in J_{P_a}} |C_{_{P_a}I_{j}}(H)| + \sum_{j \in J_{P_b}} |C_{_{P_b}I_{j}}(H)|
\]
with $J_B := \{1,2,5,8\}$, $J_{P_a} := \{1,3,4,5,6\}$ and 
$J_{P_b} := \{1,5,6,7,8,11\}$.
By \cite[Lemma 6.2]{HimstedtHuang2F4Borel} and Table~A.16, we have
\[
k(G, B_0, d, u) + k(B, B_0, d, u) = |C_{_GI^{ss}}(H)| + \sum_{j \in J_B} |C_{_BI_{j}}(H)| = 2\cdot 2^{2t}
\]
and
\[
k(P_a, B_0, d, u) + k(P_b, B_0, d, u) = \sum_{j \in J_{P_a}} |C_{_{P_a}I_{j}}(H)| + \sum_{j \in J_{P_b}} |C_{_{P_b}I_{j}}(H)| = 2\cdot 2^{2t}.
\]
Thus (\ref{eq:altsum}) also holds in this case. \hfill $\Box$

\newpage


\section*{Appendix A}

\bigskip
\noindent Table~A.1 \textit{Parameter sets for the irreducible
  characters of $G$, $B$, $P_a$, $P_b$.
 (For the parameter sets ${_GI_{34}}$, ${_GI_{35}}$, ${_GI_{36}}$,
  \dots, ${_GI_{41}}$, ${_GI_{42}}$, ${_GI_{46}}$,
  ${_GI_{50}}$, ${_GI_{51}}$ see the remarks in
  Section~\ref{sec:action}. In this table, we use the abbreviations  
  $k':=\frac{q}{\sqrt{2}}(k+l)$, $l':=\frac{q}{\sqrt{2}}(k-l)$, 
  $a' := -\sqrt{2}q^3+\sqrt{2}q+1$, 
  $b' := -\sqrt{2}q^3-q^2+\sqrt{2}q+2$,
  $a'' := \sqrt{2}q^3-\sqrt{2}q+1$ and
  $b'' := \sqrt{2}q^3-q^2-\sqrt{2}q+2$.)}  

\medskip

\begin{center}
\begin{small}
\begin{tabular}{|l|l|l|l|l|} \hline
\rule{0cm}{0.3cm}
\hspace{-0.25cm} Parameter Set & \hspace{-0.25cm} Characters \hspace{-0.25cm} &
                 Parameters & Equivalence & Number of\hspace{-0.2cm} \\ 
                        &               &      &  Relation  & Characters\hspace{-0.2cm} 
\rule[-0.1cm]{0cm}{0.3cm}\\
\hline \cline{1-3} \hline
\rule{0cm}{0.3cm}
\hspace{-0.1cm}${_GI_{1}} = .. = {_GI_{21}}$ & $\chi_1, \dots, \chi_{21}$ &&& $1$
\rule[-0.2cm]{0cm}{0.1cm}\\
\hline
\rule{-0.12cm}{0.1cm}
${_GI_{22}} = .. ={_GI_{27}}$ & $\chi_{22}(k), .. ,
\chi_{27}(k)$\hspace{-0.2cm} & $k=0, \dots, q^2-2$ & $\{k
\equiv -k\}$ & $\frac{q^2-2}{2}$\\
&& \, $k \neq 0$ &&
\rule[-0.2cm]{0cm}{0.1cm}\\
\hline
\rule{-0.12cm}{0.2cm}
${_GI_{28}}$ & $\chi_{28}(k,l)$ & \multicolumn{2}{c|}{see the remarks in
  Section~\ref{sec:action}} & 
 $\frac{q^4-10q^2+16}{16}$\hspace{-0.2cm} 
\rule[-0.2cm]{0cm}{0.1cm}\\
\hline
\rule{-0.12cm}{0.1cm}
${_GI_{29}} = .. = {_GI_{31}}$ & $\chi_{29},\chi_{30},\chi_{31}$\hspace{-0.2cm} & & & $1$ 
\rule[-0.2cm]{0cm}{0.1cm}\\
\hline
\rule{-0.12cm}{0.1cm}
${_GI_{32}} = {_GI_{33}}$ & $\chi_{32}(k),\chi_{33}(k)$ & $k=0, \dots, q^2$ &
$\{k \equiv -k\}$ & $\frac{q^2-2}{2}$\\
\hline
\end{tabular}
\end{small}
\end{center}

\newpage

\bigskip
\noindent Table~A.1 (cont.) \textit{Parameter sets for irreducible
  characters.}

\vspace{0.3cm}

\begin{center}
\begin{small}
\begin{tabular}{|l|l|l|l|l|} \hline
\rule{0cm}{0.3cm}
\hspace{-0.25cm} Parameter Set & \hspace{-0.25cm} Characters \hspace{-0.25cm} &
                 Parameters & \hspace{-0.2cm} Equivalence \hspace{-0.2cm}
                 & \hspace{-0.2cm} Number of\hspace{-0.2cm} \\ 
                        &               &      & \hspace{-0.2cm} Relation
                 & \hspace{-0.2cm} Characters\hspace{-0.2cm} 
\rule[-0.1cm]{0cm}{0.3cm}\\
\hline \cline{1-3} \hline
\rule{0cm}{0.1cm}
\hspace*{-0.25cm} ${_BI_{11}} = {_BI_{12}} $ & ${_B\chi}_{11}(k), {_B\chi}_{12}(k)$ & $k=0, \dots, q^2-2$\hspace{-0.15cm} && $q^2-1$
\rule[-0.2cm]{0cm}{0.1cm}\\
\hline
\rule{-0.12cm}{0.1cm}
${_BI_{13}} = .. = {_BI_{22}}$ & ${_B\chi}_{ 13 }, \dots, {_B\chi}_{22}$ & && $1$
\rule[-0.2cm]{0cm}{0.1cm}\\
\hline
\rule{-0.12cm}{0.1cm}
${_BI_{23}}$ & ${_B\chi}_{ 23 }(k)$ & $k=0, \dots, q^2-2$ & & $q^2-1$
\rule[-0.2cm]{0cm}{0.1cm}\\
\hline
\rule{-0.12cm}{0.1cm}
${_BI_{24}}$ & ${_B\chi}_{ 24 }$ & && $1$\hspace{-0.15cm} 
\rule[-0.2cm]{0cm}{0.1cm}\\
\hline
\rule{-0.12cm}{0.1cm}
${_BI_{25}}$ & ${_B\chi}_{ 25 }(k)$ & $k=1, \dots, q^2$ & & $q^2$
\rule[-0.2cm]{0cm}{0.1cm}\\
\hline
${_BI_{26}} = {_BI_{27}}$ & ${_B\chi}_{ 26 }(k), {_B\chi}_{ 27 }(k)$ & $k=0, \dots, q^2-2$\hspace{-0.15cm} & & $q^2-1$
\rule[-0.2cm]{0cm}{0.1cm}\\
\hline
\end{tabular}
\end{small}
\end{center}

\newpage

\noindent Table~A.1 (cont.) \textit{Parameter sets for
  irreducible characters.}


\begin{center}
\begin{small}
\begin{tabular}{|l|l|l|l|l|} \hline
\rule{0cm}{0.3cm}
\hspace{-0.25cm} Parameter Set & Characters & Parameters & \hspace{-0.15cm} Equivalence &
Number of\hspace{-0.2cm} \\
&            &            & \hspace{-0.15cm} Relation  & Characters\hspace{-0.2cm} 
\rule[-0.1cm]{0cm}{0.3cm}\\
\hline \cline{1-3} \hline
\rule{0cm}{0.3cm}
\hspace{-0.1cm}${_{P_a}I_{34}} = {_{P_a}I_{35}}$ & ${_{P_a}\chi}_{34}(k), {_{P_a}\chi}_{35}(k)$ \hspace{-0.25cm} & $k=0, \dots, q^2-2$ & & $q^2-1$
\rule[-0.2cm]{0cm}{0.2cm}\\
\hline
\rule{-0.12cm}{0.2cm}
${_{P_a}I_{36}} = {_{P_a}I_{37}}$ & ${_{P_a}\chi}_{36}, {_{P_a}\chi}_{37}$& & & $1$
\rule[-0.2cm]{0cm}{0.2cm}\\
\hline
\rule{-0.12cm}{0.2cm}
${_{P_a}I_{38}} = {_{P_a}I_{39}}$ & ${_{P_a}\chi}_{38}(k), {_{P_a}\chi}_{39}(k),$ \hspace{-0.25cm} & $k=1, \dots, q^2$ & & $q^2$\\
\hspace{0.73cm} $= {_{P_a}I_{40}}$ & ${_{P_a}\chi}_{40}(k)$ & & & 
\rule[-0.2cm]{0cm}{0.2cm}\\
\hline \cline{1-3} \hline
\rule{-0.12cm}{0.2cm}
${_{P_b}I_{1}} = .. = {_{P_b}I_{4}}$ & ${_{P_b}\chi}_{1}(k), .., {_{P_b}\chi}_{4}(k)$ \hspace{-0.25cm} & $k=0, \dots, q^2-2$ &
 & $q^2-1$
\rule[-0.2cm]{0cm}{0.2cm}\\
\hline
\rule{-0.12cm}{0.2cm}
${_{P_b}I_{5}}$ & ${_{P_b}\chi}_{5}(k,l)$& $k,l=0, \dots, q^2-2$ \hspace{-0.25cm} &
 $\{(k,l) \equiv$ & $\frac{q^4-3q^2+2}{2}$\\
\hline
\end{tabular}
\end{small}
\end{center}

\newpage

\noindent Table~A.1 (cont.) \textit{Parameter sets for irreducible characters.} 

\medskip

\begin{center}
\begin{small}
\begin{tabular}{|l|l|l|l|l|} \hline
\rule{0cm}{0.3cm}
\hspace{-0.25cm} Parameter Set & Characters & Parameters & \hspace{-0.15cm} Equivalence &
Number of\hspace{-0.2cm} \\
                        &            &            & \hspace{-0.15cm} Relation  &
                 Characters\hspace{-0.2cm} 
\rule[-0.1cm]{0cm}{0.3cm}\\
\hline \cline{1-3} \hline
\rule{-0.12cm}{0.2cm}
${_{P_b}I_{49}}$ & ${_{P_b}\chi}_{49}(k)$ & $k=0, \dots,
q^2-\sqrt{2}q$\hspace{-0.15cm} & $\{k\equiv -k\equiv$ &
$\frac{q^2-\sqrt{2}q}{4}$\\ 
 & & \, $k\neq 0$ & \hspace{0.12cm} $q^2k \equiv -q^2k\}$ &
\rule[-0.2cm]{0cm}{0.1cm}\\
\hline
\rule{-0.12cm}{0.1cm}
${_{P_b}I_{50}} = .. = {_{P_b}I_{53}}$ & ${_{P_b}\chi}_{50}, .. , {_{P_b}\chi}_{53}$ & && $1$
\rule[-0.2cm]{0cm}{0.1cm}\\
\hline
\rule{-0.12cm}{0.1cm}
${_{P_b}I_{54}}$ & ${_{P_b}\chi}_{54}(k)$ & $k=0, \dots, q^2-2$ & $\{k\equiv -k\}$ & $\frac{q^2-2}{2}$\\
 & & \, $k\neq 0$ & &
\rule[-0.2cm]{0cm}{0.1cm}\\
\hline
\rule{-0.12cm}{0.1cm}
${_{P_b}I_{55}}$ & ${_{P_b}\chi}_{55}(k)$ & $k=0, \dots,
q^2+\sqrt{2}q$ & $\{k\equiv -k\equiv$ & $\frac{q^2+\sqrt{2}q}{4}$\\ 
 & & \, $k\neq 0$ & \hspace{0.12cm} $q^2k \equiv -q^2k\}$ &
\rule[-0.2cm]{0cm}{0.1cm}\\
\hline
\rule{-0.12cm}{0.1cm}
${_{P_b}I_{56}}$ & ${_{P_b}\chi}_{56}(k)$ & $k=0, \dots,
q^2-\sqrt{2}q$ & $\{k\equiv -k\equiv$ & $\frac{q^2-\sqrt{2}q}{4}$\\ 
 & & \, $k\neq 0$ & \hspace{0.12cm} $q^2k \equiv -q^2k\}$ &
\rule[-0.2cm]{0cm}{0.1cm}\\
\hline
\end{tabular}
\end{small}
\end{center}

\bigskip
\bigskip

\noindent Table~A.2 \textit{The irreducible characters of the
  chain normalizers of defect $11n+6$.} 

\medskip

\begin{center}
\begin{small}
\begin{tabular}{|l|l|l|c|l|l|} \hline
\rule{0cm}{0.3cm}
Group & Character & Degree \hspace{1.6cm} & Parameter & Number \hspace{0.42cm} & Class \hspace{1.2cm} 
\rule[-0.1cm]{0cm}{0.3cm}\\
\hline \cline{1-3} \hline
\rule{0cm}{0.4cm}
$G$ & $\chi_{19}$ & $\frac{q^{13}}{\sqrt{2}}\phi_1 \phi_2 \phi_4^2 \phi_{12}$ & ${_GI_{19}}$ & $1$ & $c_{1,3}$
\rule[-0.2cm]{0cm}{0.5cm}\\
\rule{-0.12cm}{0.1cm}
    & $\chi_{20}$ & $\frac{q^{13}}{\sqrt{2}}\phi_1 \phi_2 \phi_4^2 \phi_{12}$ & ${_GI_{20}}$ & $1$& $c_{1,3}$ 
\rule[-0.2cm]{0cm}{0.5cm}\\
\hline
\rule{0cm}{0.4cm}
$P_b$ & ${_{P_b}\chi}_{46}$ &   $\frac{q^{13}}{\sqrt{2}}\phi_1 \phi_2$ & ${_{P_b}I_{46}}$ & $1$ & $c_{1,9}$
\rule[-0.2cm]{0cm}{0.3cm}\\
\rule{0cm}{0.4cm}
    & ${_{P_b}\chi}_{53}$ &   $\frac{q^{13}}{\sqrt{2}}\phi_1 \phi_2$ & ${_{P_b}I_{53}}$ & $1$ & $c_{1,9}$
\rule[-0.2cm]{0cm}{0.3cm}\\
\hline
\end{tabular}
\end{small}
\end{center}

\bigskip
\bigskip

\noindent Table~A.3 \textit{The irreducible characters of the
  chain normalizers of defect $14n+7$.} 

\medskip

\begin{center}
\begin{small}
\begin{tabular}{|l|l|l|c|l|l|} \hline
\rule{0cm}{0.3cm}
Group & Character & Degree \hspace{1.6cm} & Parameter & Number \hspace{0.42cm} & Class \hspace{1.2cm}
\rule[-0.1cm]{0cm}{0.3cm}\\
\hline \cline{1-3} \hline
\rule{0cm}{0.4cm}
$G$ & $\chi_{18}$ & $q^{10}\phi_{12} \phi_{24}$ & ${_GI_{18}}$ & $1$ & 
\rule[-0.2cm]{0cm}{0.5cm}\\
\hline
\rule{0cm}{0.4cm}
$P_a$ & ${_{P_a}\chi}_{31}$ & $q^{10}\phi_1 \phi_2$ & ${_{P_a}I_{31}}$ & $1$ &
\rule[-0.2cm]{0cm}{0.3cm}\\
\hline
\end{tabular}
\end{small}
\end{center}

\bigskip
\bigskip

\noindent Table~A.4 \textit{The irreducible characters of the
  chain normalizers of defect $14n+8$.} 

\medskip

\begin{center}
\begin{small}
\begin{tabular}{|l|l|l|c|l|l|} \hline
\rule{0cm}{0.3cm}
Group & Character & Degree \hspace{1.6cm} & Parameter & Number \hspace{0.42cm} & Class \hspace{1.2cm}
\rule[-0.1cm]{0cm}{0.3cm}\\
\hline \cline{1-3} \hline
\rule{0cm}{0.4cm}
$B$ & ${_B\chi}_{55}$ & $\frac{1}{2}q^{10}\phi_1^2 \phi_2^2$ & ${_BI_{55}}$ & $1$ & $c_{1,57}, c_{1,62}$
\rule[-0.2cm]{0cm}{0.5cm}\\
\rule{0cm}{0.4cm}
    & ${_B\chi}_{56}$ & $\frac{1}{2}q^{10}\phi_1^2 \phi_2^2$ & ${_BI_{56}}$ & $1$ & $c_{1,57}, c_{1,62}$
\rule[-0.2cm]{0cm}{0.5cm}\\
\rule{0cm}{0.4cm}
    & ${_B\chi}_{57}$ & $\frac{1}{2}q^{10}\phi_1^2 \phi_2^2$ & ${_BI_{57}}$ & $1$ & $c_{1,57}, c_{1,62}$
\rule[-0.2cm]{0cm}{0.5cm}\\
\rule{0cm}{0.4cm}
    & ${_B\chi}_{58}$ & $\frac{1}{2}q^{10}\phi_1^2 \phi_2^2$ & ${_BI_{58}}$ & $1$ & $c_{1,57}, c_{1,62}$
\rule[-0.2cm]{0cm}{0.5cm}\\
\hline
\rule{0cm}{0.4cm}
$P_b$ & ${_{P_b}\chi}_{44}$ & $\frac{1}{2}q^{10}\phi_1^2 \phi_2^2$ & ${_{P_b}I_{44}}$ & $1$ & $c_{1,31}, c_{1,33}$ \rule[-0.2cm]{0cm}{0.3cm}\\
\rule{0cm}{0.4cm}
      & ${_{P_b}\chi}_{45}$ & $\frac{1}{2}q^{10}\phi_1^2 \phi_2^2$ & ${_{P_b}I_{45}}$ & $1$ & $c_{1,31}, c_{1,33}$ \rule[-0.2cm]{0cm}{0.3cm}\\
\rule{0cm}{0.4cm}
      & ${_{P_b}\chi}_{51}$ & $\frac{1}{2}q^{10}\phi_1^2 \phi_2^2$ & ${_{P_b}I_{51}}$ & $1$ & $c_{1,31}, c_{1,33}$ \rule[-0.2cm]{0cm}{0.3cm}\\
\rule{0cm}{0.4cm}
      & ${_{P_b}\chi}_{52}$ & $\frac{1}{2}q^{10}\phi_1^2 \phi_2^2$ & ${_{P_b}I_{52}}$ & $1$ & $c_{1,31}, c_{1,33}$ \rule[-0.2cm]{0cm}{0.3cm}\\
\hline
\end{tabular}
\end{small}
\end{center}

\newpage

\noindent Table~A.5 \textit{The irreducible characters of the
  chain normalizers of defect $15n+8$.} 

\begin{center}
\begin{small}
\begin{tabular}{|l|l|l|c|l|l|} \hline
\rule{0cm}{0.3cm}
Group & Character & Degree & Parameter & Number & Class
\rule[-0.1cm]{0cm}{0.3cm}\\
\hline \cline{1-3} \hline
\rule{0cm}{0.4cm}
$B$ & ${_B\chi}_{51}(k)$ & $\frac{q^9}{\sqrt{2}}\phi_1 \phi_2$ & ${_BI_{51}}$ & $q^2-1$ & $c_{5,0}(i), c_{5,2}(i)$ 
\rule[-0.2cm]{0cm}{0.1cm}\\
\rule{-0.12cm}{0.1cm}
 & ${_B\chi}_{52}(k)$ & $\frac{q^9}{\sqrt{2}}\phi_1 \phi_2$ & ${_BI_{52}}$ & $q^2-1$& 
$c_{5,0}(i), c_{5,2}(i)$ \rule[-0.2cm]{0cm}{0.5cm}\\
\rule{-0.12cm}{0.1cm}
 & ${_B\chi}_{53}$ & $\frac{q^9}{\sqrt{2}}\phi_1^2 \phi_2^2$ & ${_BI_{53}}$ & $1$& $c_{1,12}$ 
\rule[-0.2cm]{0cm}{0.5cm}\\
\rule{-0.12cm}{0.1cm}
 & ${_B\chi}_{54}$ & $\frac{q^9}{\sqrt{2}}\phi_1^2 \phi_2^2$ & ${_BI_{54}}$ & $1$& $c_{1,12}$ 
\rule[-0.2cm]{0cm}{0.5cm}\\
\hline
\rule{0cm}{0.4cm}
$P_b$ & ${_{P_b}\chi}_{43}$ & $\frac{q^9}{\sqrt{2}}\phi_1 \phi_2$ & ${_{P_b}I_{43}}$ & $1$ & $c_{1,9}$ 
\rule[-0.2cm]{0cm}{0.1cm}\\
\rule{-0.12cm}{0.1cm}
    & ${_{P_b}\chi}_{47}(k)$ & $\frac{q^9}{\sqrt{2}}\phi_1 \phi_2 \phi_8$ & ${_{P_b}I_{47}}$ & $\frac{1}{2}(q^2-2)$ & $c_{5,0}(i), c_{5,2}(i)$ 
\rule[-0.2cm]{0cm}{0.1cm}\\
\rule{-0.12cm}{0.1cm}
    & ${_{P_b}\chi}_{48}(k)$ & $\frac{q^9}{\sqrt{2}}\phi_1^2 \phi_2^2 \phi_8''$ & ${_{P_b}I_{48}}$ & $\frac{1}{4}(q^2+\sqrt{2}q)$ & $c_{10,0}(i), c_{10,2}(i)$ 
\rule[-0.2cm]{0cm}{0.1cm}\\
\hline
\end{tabular}
\end{small}
\end{center}

\bigskip

\noindent Table~A.6 \textit{The irreducible characters of the
  chain normalizers of defect $16n+8$.}

\begin{center}
\begin{small}
\begin{tabular}{|l|l|l|c|l|l|} \hline
\rule{0cm}{0.3cm}
Group & Character & Degree & Parameter & Number \hspace{0.6cm} & Class \hspace{1.3cm}
\rule[-0.1cm]{0cm}{0.1cm}\\
\hline \cline{1-3} \hline
\rule{0cm}{0.1cm}
$B$ & ${_B\chi}_{42}(k)$ & $q^{8}\phi_1 \phi_2$ & ${_BI_{42}}$ & $q^2-1$ & $c_{9,0}(i)$ 
\rule[-0.2cm]{0cm}{0.1cm}\\
\rule{0cm}{0.1cm}
    & ${_B\chi}_{43}$ & $q^{8}\phi_1^2 \phi_2^2$ & ${_BI_{43}}$ & $1$ & $c_{1,2}$ 
\rule[-0.2cm]{0cm}{0.1cm}\\
\rule{0cm}{0.1cm}
    & ${_B\chi}_{50}(k)$ & $q^{8}\phi_1^2 \phi_2^2$ & ${_BI_{50}}$ & $q^2$ & 
\rule[-0.2cm]{0cm}{0.1cm}\\
\hline
\rule{0cm}{0.1cm}
$P_a$ & ${_{P_a}\chi}_{30}$ & $q^{8}\phi_1 \phi_2$ & ${_{P_a}I_{30}}$ & $1$ & 
\rule[-0.2cm]{0cm}{0.1cm}\\
\rule{0cm}{0.1cm}
    & ${_{P_a}\chi}_{32}(k)$ & $q^{8}\phi_1 \phi_2 \phi_4$ & ${_{P_a}I_{32}}$ & $\frac{q^2-2}{2}$ & $c_{6,0}(i)$
\rule[-0.2cm]{0cm}{0.1cm}\\
\rule{0cm}{0.1cm}
    & ${_{P_a}\chi}_{33}(k)$ & $q^{8}\phi_1^2 \phi_2^2$ & ${_{P_a}I_{33}}$ & $\frac{q^2}{2}$ & $c_{9,0}(i)$ 
\rule[-0.2cm]{0cm}{0.1cm}\\
\rule{0cm}{0.1cm}
    & ${_{P_a}\chi}_{40}(k)$ & $q^{8}\phi_1^2 \phi_2^2 \phi_4$ & ${_{P_a}I_{40}}$ & $q^2$ & 
\rule[-0.2cm]{0cm}{0.1cm}\\
\hline
\end{tabular}
\end{small}
\end{center}

\bigskip

\noindent Table~A.7 \textit{The irreducible characters of the
  chain normalizers of defect $17n+9$.} 

\begin{center}
\begin{small}
\begin{tabular}{|l|l|l|c|l|l|} \hline
\rule{0cm}{0.3cm}
Group & Character & Degree & Parameter & Number \hspace{0.6cm} & Class \hspace{1.3cm}
\rule[-0.1cm]{0cm}{0.3cm}\\
\hline \cline{1-3} \hline
\rule{0cm}{0.4cm}
$B$ & ${_B\chi}_{44}(k)$ & $\frac{q^7}{\sqrt{2}}\phi_1 \phi_2$ &
${_BI_{44}}$ & $q^2-1$ & $c_{6,0}(i), c_{6,2}(i)$  
\rule[-0.2cm]{0cm}{0.1cm}\\
\rule{0cm}{0.4cm}
    & ${_B\chi}_{45}(k)$ & $\frac{q^7}{\sqrt{2}}\phi_1 \phi_2$ &
${_BI_{45}}$ & $q^2-1$ & $c_{6,0}(i), c_{6,2}(i)$  
\rule[-0.2cm]{0cm}{0.1cm}\\
\rule{0cm}{0.1cm}
    & ${_B\chi}_{46}$ & $\frac{q^7}{\sqrt{2}}\phi_1^2 \phi_2^2$ & ${_BI_{46}}$ & $1$ & 
$c_{1,9}$ \rule[-0.2cm]{0cm}{0.1cm}\\
\rule{0cm}{0.1cm}
    & ${_B\chi}_{47}$ & $\frac{q^7}{\sqrt{2}}\phi_1^2 \phi_2^2$ & ${_BI_{47}}$ & $1$ & 
$c_{1,9}$ \rule[-0.2cm]{0cm}{0.1cm}\\
\rule{0cm}{0.1cm}
    & ${_B\chi}_{48}(k)$ & $\frac{q^7}{\sqrt{2}}\phi_1^2 \phi_2^2$ & ${_BI_{48}}$ & $q^2$ & 
\rule[-0.2cm]{0cm}{0.1cm}\\
\rule{0cm}{0.1cm}
    & ${_B\chi}_{49}(k)$ & $\frac{q^7}{\sqrt{2}}\phi_1^2 \phi_2^2$ & ${_BI_{49}}$ & $q^2$ & 
\rule[-0.2cm]{0cm}{0.1cm}\\
\hline
\rule{0cm}{0.1cm}
$P_a$ & ${_{P_a}\chi_{34}(k)}$ & $\frac{q^7}{\sqrt{2}}\phi_1 \phi_2
\phi_4$ & ${_{P_a}I_{34}}$ & $q^2-1$ & $c_{4,0}(i), c_{4,2}(i)$  
\rule[-0.2cm]{0cm}{0.1cm}\\
\hline
\end{tabular}
\end{small}
\end{center}

\vspace{0.8cm}

\noindent Table~A.8 \textit{The irreducible characters of the
  chain normalizers of defect $18n+9$.} 

\smallskip

\begin{center}
\begin{small}
\begin{tabular}{|l|l|l|c|l|l|} \hline
\rule{0cm}{0.3cm}
Group & Character & Degree \hspace{1.4cm} & Parameter &
Number \hspace{0.6cm} & Class \hspace{0.8cm}
\rule[-0.1cm]{0cm}{0.3cm}\\
\hline \cline{1-3} \hline
\rule{0cm}{0.4cm}
$G$ & ${\chi}_{31}$ & $q^{6}\phi_1 \phi_2 \phi_8^2 \phi_{24}$ & ${_GI_{31}}$ & $1$ & 
\rule[-0.2cm]{0cm}{0.5cm}\\
\hline
\rule{0cm}{0.4cm}
$B$ & ${_B\chi}_{38}$ & $q^{6}\phi_1^2 \phi_2^2$ & ${_BI_{38}}$ & $1$ & $c_{1,4}$
\rule[-0.2cm]{0cm}{0.5cm}\\
\rule{0cm}{0.4cm}
    & ${_B\chi}_{39}(k)$ & $q^{6}\phi_1 \phi_2$ & ${_BI_{39}}$ & $q^2-1$ & $c_{8,0}(i)$
\rule[-0.2cm]{0cm}{0.5cm}\\
\rule{0cm}{0.4cm}
    & ${_B\chi}_{40}$ & $q^{6}\phi_1^2 \phi_2^2$ & ${_BI_{40}}$ & $1$ & $c_{1,4}$
\rule[-0.2cm]{0cm}{0.5cm}\\
\rule{0cm}{0.4cm}
    & ${_B\chi}_{41}(k)$ & $q^{6}\phi_1^2 \phi_2^2$ & ${_BI_{41}}$ & $q^2$ & 
\rule[-0.2cm]{0cm}{0.5cm}\\
\hline
\rule{0cm}{0.4cm}
$P_a$ & ${_{P_a}\chi_{29}}$ & $q^{6}\phi_1^2 \phi_2^2 \phi_4$ & & $1$ & 
\rule[-0.2cm]{0cm}{0.3cm}\\
\hline
\rule{0cm}{0.4cm}
$P_b$ & ${_{P_b}\chi}_{39}$ & $q^{6}\phi_1^2 \phi_2^2 \phi_8$ & ${_{P_b}I_{39}}$ & $1$ & $c_{1,3}$
\rule[-0.2cm]{0cm}{0.5cm}\\
\rule{0cm}{0.4cm}
    & ${_{P_b}\chi}_{40}(k)$ & $q^{6}\phi_1 \phi_2 \phi_8$ & ${_{P_b}I_{40}}$ & $q^2-1$ & $c_{6,0}(i)$
\rule[-0.2cm]{0cm}{0.5cm}\\
\rule{0cm}{0.4cm}
    & ${_{P_b}\chi}_{41}$ & $q^{6}\phi_1^2 \phi_2^2 \phi_8$ & ${_{P_b}I_{41}}$ & $1$ & $c_{1,3}$
\rule[-0.2cm]{0cm}{0.5cm}\\
\rule{0cm}{0.4cm}
    & ${_{P_b}\chi}_{42}(k)$ & $q^{6}\phi_1^2 \phi_2^2 \phi_8$ & ${_{P_b}I_{42}}$ & $q^2$ & 
\rule[-0.2cm]{0cm}{0.5cm}\\
\hline
\end{tabular}
\end{small}
\end{center}

\bigskip

\noindent Table~A.9 \textit{The irreducible characters of the
  chain normalizers of defect $20n+10$.} 

\smallskip

\begin{center}
\begin{small}
\begin{tabular}{|l|l|l|c|l|l|} \hline
\rule{0cm}{0.3cm}
Group & Character & Degree \hspace{1.4cm} & Parameter & Number \hspace{0.6cm} & Class \hspace{0.8cm}
\rule[-0.1cm]{0cm}{0.3cm}\\
\hline \cline{1-3} \hline
\rule{0cm}{0.4cm}
$G$ & ${\chi}_{15}$ & $\frac{1}{3}q^{4}\phi_1^2 \phi_2^2 \phi_4^2 \phi_8^2$ & ${_GI_{15}}$ & $1$ & $c_{5,3}$
\rule[-0.2cm]{0cm}{0.5cm}\\
\rule{0cm}{0.4cm}
    & ${\chi}_{16}$ & $\frac{1}{3}q^{4}\phi_1^2 \phi_2^2 \phi_4^2\phi_8^2$ & ${_GI_{16}}$ & $1$ & $c_{5,3}$
\rule[-0.2cm]{0cm}{0.5cm}\\
\rule{0cm}{0.4cm}
    & ${\chi}_{17}$ & $\frac{1}{3}q^{4}\phi_1^2 \phi_2^2 \phi_{12} \phi_{24}$ & ${_GI_{17}}$ & $1$ & 
\rule[-0.2cm]{0cm}{0.5cm}\\
\rule{0cm}{0.4cm}
    & ${\chi}_{25}(k)$ & $q^{4}\phi_4^2 \phi_8 \phi_{12} \phi_{24}$ & ${_GI_{25}}$ & $\frac{1}{2}(q^2-2)$ & $c_{2,0}(i)$
\rule[-0.2cm]{0cm}{0.5cm}\\
\rule{0cm}{0.4cm}
    & ${\chi}_{45}(k)$ & $q^{4}\phi_1\phi_2\phi_4^2 \phi_8' \phi_{12} \phi_{24}$ & ${_GI_{45}}$ & $\frac{1}{4}(q^2-\sqrt{2}q)$ & $c_{8,2}(i)$
\rule[-0.2cm]{0cm}{0.5cm}\\
\rule{0cm}{0.4cm}
    & ${\chi}_{49}(k)$ & $q^{4}\phi_1\phi_2\phi_4^2 \phi_8'' \phi_{12} \phi_{24}$ & ${_GI_{49}}$ & $\frac{1}{4}(q^2+\sqrt{2}q)$ & $c_{10,2}(i)$
\rule[-0.2cm]{0cm}{0.5cm}\\
\hline
\rule{0cm}{0.4cm}
$B$ & ${_B\chi}_{11}(k)$ & $q^{4}\phi_1 \phi_2$ & ${_BI_{11}}$ & $q^2-1$ & $c_{5,0}(i)$ 
\rule[-0.2cm]{0cm}{0.5cm}\\
\hline
\end{tabular}
\end{small}
\end{center}

\newpage
\noindent Table~A.9 (cont.) \textit{Irreducible characters of the
  chain normalizers of defect $20n+10$.} 

\bigskip

\begin{center}
\begin{small}
\begin{tabular}{|l|l|l|c|l|l|} \hline
\rule{0cm}{0.3cm}
Group & Character & Degree \hspace{1.4cm} & Parameter & Number \hspace{0.6cm} & Class \hspace{0.8cm}
\rule[-0.1cm]{0cm}{0.3cm}\\
\hline \cline{1-3} \hline
\rule{0cm}{0.4cm}
$P_b$ & ${_{P_b}\chi}_{4}(k)$ & $q^{4}$ & ${_{P_b}I_{4}}$ & $q^2-1$ & $c_{2,0}(i)$
\rule[-0.2cm]{0cm}{0.5cm}\\
\rule{0cm}{0.4cm}
    & ${_{P_b}\chi}_{16}(k)$ & $q^{4}\phi_1 \phi_2 \phi_8$ & ${_{P_b}I_{16}}$ & $\frac{1}{2}(q^2-2)$ & $c_{5,0}(i)$
\rule[-0.2cm]{0cm}{0.5cm}\\
\rule{0cm}{0.4cm}
    & ${_{P_b}\chi}_{21}(k)$ & $q^{4}\phi_1^2 \phi_2^2 \phi_8'$ & ${_{P_b}I_{21}}$ & $\frac{1}{4}(q^2-\sqrt{2}q)$ & $c_{8,0}(i)$
\rule[-0.2cm]{0cm}{0.5cm}\\
\rule{0cm}{0.4cm}
    & ${_{P_b}\chi}_{26}(k)$ & $q^{4}\phi_1^2 \phi_2^2 \phi_8''$ & ${_{P_b}I_{26}}$ & $\frac{1}{4}(q^2+\sqrt{2}q)$ & $c_{10,0}(i)$
\rule[-0.2cm]{0cm}{0.5cm}\\
\rule{0cm}{0.4cm}
    & ${_{P_b}\chi}_{35}$ & $q^{4}\phi_1^2 \phi_2^2 \phi_8$ & ${_{P_b}I_{35}}$ & $1$ & 
\rule[-0.2cm]{0cm}{0.5cm}\\
\rule{0cm}{0.4cm}
    & ${_{P_b}\chi}_{38}(k)$ & $q^{4}\phi_1^2 \phi_2^2 \phi_8$ & ${_{P_b}I_{38}}$ & $\frac{1}{3}(q^2+1)$ & 
\rule[-0.2cm]{0cm}{0.5cm}\\
\hline
\end{tabular}
\end{small}
\end{center}

\bigskip
\bigskip

\noindent Table~A.10 \textit{The irreducible characters of the
  chain normalizers of defect $20n+11$.} 

\medskip

\begin{center}
\begin{small}
\begin{tabular}{|l|l|l|c|l|l|} \hline
\rule{0cm}{0.3cm}
Group & Character & Degree\hspace{1.4cm} & Parameter &
Number\hspace{0.6cm} & Class \hspace{0.8cm} 
\rule[-0.1cm]{0cm}{0.3cm}\\
\hline \cline{1-3} \hline
\rule{0cm}{0.4cm}
$G$ & ${\chi}_{7}$ & $\frac{q^{4}}{2}\phi_8^2 \phi_{24}$ & ${_GI_{7}}$ & $1$ & 
\rule[-0.2cm]{0cm}{0.5cm}\\
\rule{0cm}{0.4cm}
    & ${\chi}_{10}$ & $\frac{q^{4}}{6}\phi_1^2 \phi_2^2 \phi_4^2 \phi_{24}$ & ${_GI_{10}}$ & $1$ & 
\rule[-0.2cm]{0cm}{0.5cm}\\
\hline
\rule{0cm}{0.4cm}
$B$ & ${_B\chi}_{13}$ & $\frac{q^{4}}{2}\phi_1^2 \phi_2^2$ &
${_BI_{13}}$ & $1$ & $c_{1,5}, c_{1,38}$
\rule[-0.2cm]{0cm}{0.5cm}\\
\rule{0cm}{0.4cm}
    & ${_B\chi}_{14}$ & $\frac{q^{4}}{2}\phi_1^2 \phi_2^2$ &
${_BI_{14}}$ & $1$ & $c_{1,5}, c_{1,38}$ \rule[-0.2cm]{0cm}{0.5cm}\\ 
\rule{0cm}{0.4cm}
    & ${_B\chi}_{30}$ & $\frac{q^{4}}{2}\phi_1^2 \phi_2^2$ &
${_BI_{30}}$ & $1$ & $c_{1,17}, c_{1,54}$ \rule[-0.2cm]{0cm}{0.5cm}\\
\rule{0cm}{0.4cm}
    & ${_B\chi}_{31}$ & $\frac{q^{4}}{2}\phi_1^2 \phi_2^2$ &
${_BI_{31}}$ & $1$ & $c_{1,17}, c_{1,54}$ \rule[-0.2cm]{0cm}{0.5cm}\\ 
\rule{0cm}{0.4cm}
    & ${_B\chi}_{32}$ & $\frac{q^{4}}{2}\phi_1^2 \phi_2^2$ &
${_BI_{32}}$ & $1$ & $c_{1,17}, c_{1,54}$ \rule[-0.2cm]{0cm}{0.5cm}\\
\hline
\end{tabular}
\end{small}
\end{center}

\newpage

\noindent Table~A.11 \textit{The irreducible characters of the
  chain normalizers of defect $20n+12$.} 

\bigskip

\begin{center}
\begin{small}
\begin{tabular}{|l|l|l|c|l|l|} \hline
\rule{0cm}{0.3cm}
Group & Character & Degree \hspace{1.4cm} & Parameter &
Number \hspace{0.6cm} & Class \hspace{0.8cm} 
\rule[-0.1cm]{0cm}{0.3cm}\\
\hline \cline{1-3} \hline
\rule{0cm}{0.3cm}
$G$ & ${\chi}_{5}$ & $\frac{q^{4}}{4}\phi_4^{2}\phi_8''^2\phi_{12}\phi_{24}'$ & ${_GI_{5}}$ & $1$ & 
\rule[-0.2cm]{0cm}{0.5cm}\\
\rule{0cm}{0.3cm}
    & ${\chi}_{6}$ & $\frac{q^{4}}{4}\phi_4^{2}\phi_8'^2\phi_{12}\phi_{24}''$ & ${_GI_{6}}$ & $1$ & 
\rule[-0.2cm]{0cm}{0.5cm}\\
\rule{0cm}{0.3cm}
    & ${\chi}_{8}$ & $\frac{q^4}{12}\phi_1^{2}\phi_2^{2}\phi_8'^2\phi_{12}\phi_{24}'$ & ${_GI_{8}}$ & $1$ & 
\rule[-0.2cm]{0cm}{0.5cm}\\
\rule{0cm}{0.3cm}
    & ${\chi}_{9}$ & $\frac{q^4}{12}\phi_1^{2}\phi_2^{2}\phi_8''^2\phi_{12}\phi_{24}''$ & ${_GI_{9}}$ & $1$ & 
\rule[-0.2cm]{0cm}{0.5cm}\\
\rule{0cm}{0.3cm}
    & ${\chi}_{11}$ & $\frac{q^{4}}{4}\phi_1^{2}\phi_2^{2}\phi_4^{2}\phi_{12}\phi_{24}''$ & ${_GI_{11}}$ & $1$ & $c_{1,3}$ 
\rule[-0.2cm]{0cm}{0.5cm}\\
\rule{0cm}{0.3cm}
    & ${\chi}_{12}$ & $\frac{q^{4}}{4}\phi_1^{2}\phi_2^{2}\phi_4^{2}\phi_{12}\phi_{24}''$ & ${_GI_{12}}$ & $1$ & $c_{1,3}$ 
\rule[-0.2cm]{0cm}{0.5cm}\\
\rule{0cm}{0.3cm}
    & ${\chi}_{13}$ & $\frac{q^{4}}{4}\phi_1^{2}\phi_2^{2}\phi_4^{2}\phi_{12}\phi_{24}'$ & ${_GI_{13}}$ & $1$ & $c_{1,3}$ 
\rule[-0.2cm]{0cm}{0.5cm}\\
\hline
\end{tabular}
\end{small}
\end{center}

\newpage

\noindent Table~A.12 \textit{The irreducible characters of the
  chain normalizers of defect $21n+11$.} 

\begin{center}
\begin{small}
\begin{tabular}{|l|l|l|c|l|l|} \hline
\rule{0cm}{0.3cm}
Group & Character & Degree \hspace{1.6cm} & Parameter &
Number \hspace{0.6cm} & Class \hspace{1.1cm} 
\rule[-0.1cm]{0cm}{0.3cm}\\
\hline \cline{1-3} \hline
\rule{0cm}{0.4cm}
$B$ & ${_B\chi}_{26}(k)$ & $\frac{q^3}{\sqrt{2}}\phi_1 \phi_2$ &
${_BI_{26}}$ & $q^2-1$ & $c_{4,0}(i), c_{4,2}(i)$  
\rule[-0.2cm]{0cm}{0.5cm}\\
\rule{0cm}{0.4cm}
    & ${_B\chi}_{27}(k)$ & $\frac{q^3}{\sqrt{2}}\phi_1 \phi_2$ &
${_BI_{27}}$ & $q^2-1$ & $c_{4,0}(i), c_{4,2}(i)$  
\rule[-0.2cm]{0cm}{0.5cm}\\
\rule{0cm}{0.4cm}
    & ${_B\chi}_{28}$ & $\frac{q^3}{\sqrt{2}}\phi_1^2 \phi_2^2$ &
${_BI_{28}}$ & $1$ &  $c_{1,17}$ \rule[-0.2cm]{0cm}{0.5cm}\\ 
\rule{0cm}{0.4cm}
    & ${_B\chi}_{29}$ & $\frac{q^3}{\sqrt{2}}\phi_1^2 \phi_2^2$ & ${_BI_{29}}$ & $1$ & 
$c_{1,17}$ \rule[-0.2cm]{0cm}{0.5cm}\\
\hline
\rule{0cm}{0.4cm}
$P_b$ & ${_{P_b}\chi}_{27}(k)$ & $\frac{q^3}{\sqrt{2}}\phi_1 \phi_2
\phi_8$ & ${_{P_b}I_{27}}$ & $q^2-1$ & $c_{4,0}(i), c_{4,2}(i)$  
\rule[-0.2cm]{0cm}{0.5cm}\\
\rule{0cm}{0.4cm}
    & ${_{P_b}\chi}_{28}(k)$ & $\frac{q^3}{\sqrt{2}}\phi_1 \phi_2
\phi_8$ & ${_{P_b}I_{28}}$ & $q^2-1$ & $c_{4,0}(i), c_{4,2}(i)$  
\rule[-0.2cm]{0cm}{0.5cm}\\
\rule{0cm}{0.4cm}
    & ${_{P_b}\chi}_{29}$ & $\frac{q^3}{\sqrt{2}}\phi_1^2 \phi_2^2
\phi_8$ & ${_{P_b}I_{29}}$ & $1$ & $c_{1,6}$  
\rule[-0.2cm]{0cm}{0.5cm}\\
\rule{0cm}{0.4cm}
    & ${_{P_b}\chi}_{30}$ & $\frac{q^3}{\sqrt{2}}\phi_1^2 \phi_2^2
\phi_8$ & ${_{P_b}I_{30}}$ & $1$ & $c_{1,6}$  
\rule[-0.2cm]{0cm}{0.5cm}\\
\hline
\end{tabular}
\end{small}
\end{center}

\bigskip
\bigskip

\noindent Table~A.13 \textit{The irreducible characters of the
  chain normalizers of defect $22n+11$.}

\begin{center}
\begin{small}
\begin{tabular}{|l|l|l|c|l|l|} \hline
\rule{0cm}{0.3cm}
Group & Character & Degree \hspace{1.6cm} & Parameter &
Number \hspace{0.6cm} & Class \hspace{1.1cm} 
\rule[-0.1cm]{0cm}{0.3cm}\\
\hline \cline{1-3} \hline
\rule{0cm}{0.3cm}
$G$ & $\chi_{4}$ & $q^2\phi_{12} \phi_{24}$ & ${_GI_{4}}$ & $1$ &
\rule[-0.2cm]{0cm}{0.1cm}\\
\rule{-0.12cm}{0.1cm}
    & $\chi_{27}(k)$ & $q^2\phi_4 \phi_8^2 \phi_{12} \phi_{24}$ &
${_GI_{27}}$ & $\frac{1}{2}(q^2-2)$ & $c_{2,0}(i)$ 
\rule[-0.2cm]{0cm}{0.1cm}\\
\rule{-0.12cm}{0.1cm}
    & $\chi_{30}$ & $q^2\phi_1^2 \phi_2^2 \phi_8^2 \phi_{24}$ & ${_GI_{30}}$ & $1$ & 
\rule[-0.2cm]{0cm}{0.1cm}\\
\rule{-0.12cm}{0.1cm}
    & $\chi_{33}(k)$ & $q^2\phi_1 \phi_2 \phi_8^2 \phi_{12} \phi_{24}$
& ${_GI_{33}}$ & $\frac{1}{2}(q^2-2)$ & $c_{6,0}(i)$ \rule[-0.2cm]{0cm}{0.5cm}\\
\hline
\rule{0cm}{0.4cm}
$B$ & ${_B\chi}_{9}(k)$ & $q^2\phi_1 \phi_2$ & ${_BI_{9}}$ & $q^2-1$ & $c_{4,0}(i)$ 
\rule[-0.2cm]{0cm}{0.1cm}\\
\rule{-0.12cm}{0.1cm}
    & ${_B\chi}_{10}$ & $q^2\phi_1^2 \phi_2^2$ & ${_BI_{10}}$ & $1$ &
\rule[-0.2cm]{0cm}{0.1cm}\\
\hline
\rule{0cm}{0.4cm}
$P_a$ & ${_{P_a}\chi}_{2}(k)$  & $q^2$ & ${_{P_a}I_{2}}$ & $q^2-1$ & 
$c_{2,0}(i)$
\rule[-0.2cm]{0cm}{0.3cm}\\
\rule{0cm}{0.4cm}
      & ${_{P_a}\chi}_{11}$  & $q^2\phi_1^2 \phi_2^2 \phi_4$ & ${_{P_a}I_{11}}$ & $1$ & 
\rule[-0.2cm]{0cm}{0.3cm}\\
\hline
\rule{0cm}{0.3cm}
$P_b$ & ${_{P_b}\chi}_{12}(k)$  & $q^2\phi_1 \phi_2 \phi_8$ & ${_{P_b}I_{12}}$ & $q^2-1$ &
$c_{4,0}(i)$ \rule[-0.2cm]{0cm}{0.1cm}\\
\rule{-0.12cm}{0.1cm}
    & ${_{P_b}\chi}_{13}$  & $q^2\phi_1^2 \phi_2^2 \phi_8$ & ${_{P_b}I_{13}}$ & $1$ & 
\rule[-0.2cm]{0cm}{0.1cm}\\
\hline
\end{tabular}
\end{small}
\end{center}

\bigskip
\bigskip

\noindent Table~A.14 \textit{The irreducible characters of the
  chain normalizers of defect~$23n+12$.} 

\begin{center}
\begin{small}
\begin{tabular}{|l|l|l|c|l|l|} \hline
\rule{0cm}{0.3cm}
Group & Character & Degree & Parameter & Number & Class
\rule[-0.1cm]{0cm}{0.3cm}\\
\hline \cline{1-3} \hline
\rule{0cm}{0.3cm}
$G$ & $\chi_{2}$ & $\frac{q}{\sqrt{2}}\phi_1 \phi_2 \phi_4^2
\phi_{12}$ & ${_GI_{2}}$ & $1$ & $c_{1,4}$  
\rule[-0.2cm]{0cm}{0.1cm}\\
\rule{-0.12cm}{0.1cm}
    & $\chi_{3}$ & $\frac{q}{\sqrt{2}}\phi_1 \phi_2 \phi_4^2
\phi_{12}$ & ${_GI_{3}}$ & $1$ & $c_{1,4}$ 
\rule[-0.2cm]{0cm}{0.1cm}\\
\rule{-0.12cm}{0.1cm}
    & $\chi_{23}(k)$ & $\frac{q}{\sqrt{2}}\phi_1 \phi_2 \phi_4^2
\phi_8 \phi_{12} \phi_{24}$ & ${_GI_{23}}$ & $\frac{1}{2}(q^2-2)$ &
$c_{2,0}(i), c_{2,2}(i)$ 
\rule[-0.2cm]{0cm}{0.1cm}\\
\rule{-0.12cm}{0.1cm}
    & $\chi_{24}(k)$ & $\frac{q}{\sqrt{2}}\phi_1 \phi_2 \phi_4^2
\phi_8 \phi_{12} \phi_{24}$ & ${_GI_{24}}$ & $\frac{1}{2}(q^2-2)$ &
$c_{2,0}(i), c_{2,2}(i)$ 
\rule[-0.2cm]{0cm}{0.1cm}\\
\rule{-0.12cm}{0.1cm}
    & $\chi_{43}(k)$ & $\frac{q}{\sqrt{2}}\phi_1^2 \phi_2^2 \phi_4^2
\phi_8' \phi_{12} \phi_{24}$ & ${_GI_{43}}$ &
$\frac{1}{4}(q^2-\sqrt{2}q)$ & $c_{1,11}, c_{8,2}(i)$ 
\rule[-0.2cm]{0cm}{0.1cm}\\
\rule{-0.12cm}{0.1cm}
    & $\chi_{44}(k)$ & $\frac{q}{\sqrt{2}}\phi_1^2 \phi_2^2 \phi_4^2
\phi_8' \phi_{12} \phi_{24}$ & ${_GI_{44}}$ &
$\frac{1}{4}(q^2-\sqrt{2}q)$ & $c_{1,11}, c_{8,2}(i)$ 
\rule[-0.2cm]{0cm}{0.1cm}\\
\rule{-0.12cm}{0.1cm}
    & $\chi_{47}(k)$ & $\frac{q}{\sqrt{2}}\phi_1^2 \phi_2^2 \phi_4^2
\phi_8'' \phi_{12} \phi_{24}$ & ${_GI_{47}}$ &
$\frac{1}{4}(q^2+\sqrt{2}q)$ & $c_{1,11}, c_{10,2}(i)$ \hspace{-0.25cm}
\rule[-0.2cm]{0cm}{0.1cm}\\
\hline
\end{tabular}
\end{small}
\end{center}
\newpage

\noindent Table~A.14 (cont.) \textit{Irreducible characters of the
  chain normalizers of defect~$23n+12$.} 

\bigskip

\begin{center}
\begin{small}
\begin{tabular}{|l|l|l|c|l|l|} \hline
\rule{0cm}{0.3cm}
Group & Character & Degree & Parameter & Number \hspace{1.4cm} & Class
\rule[-0.1cm]{0cm}{0.3cm}\\
\hline \cline{1-3} \hline
\rule{0cm}{0.3cm}
$P_a$ & ${_{P_a}\chi}_{7}(k)$  & $\frac{q}{\sqrt{2}}\phi_1 \phi_2 \phi_4$ & ${_{P_a}I_{7}}$ & $q^2-1$ & $c_{3,0}(i), c_{3,2}(i)$ 
\rule[-0.2cm]{0cm}{0.3cm}\\
\rule{0cm}{0.4cm}
      & ${_{P_a}\chi}_{8}(k)$  & $\frac{q}{\sqrt{2}}\phi_1 \phi_2 \phi_4$ & ${_{P_a}I_{8}}$ & $q^2-1$ & $c_{3,0}(i), c_{3,2}(i)$
\rule[-0.2cm]{0cm}{0.3cm}\\
\rule{0cm}{0.4cm}
      & ${_{P_a}\chi}_{9}$  & $\frac{q}{\sqrt{2}}\phi_1^2 \phi_2^2 \phi_4$ & ${_{P_a}I_{9}}$ & $1$ & $c_{1,12}$
\rule[-0.2cm]{0cm}{0.3cm}\\
\rule{0cm}{0.4cm}
      & ${_{P_a}\chi}_{10}$  & $\frac{q}{\sqrt{2}}\phi_1^2 \phi_2^2 \phi_4$ & ${_{P_a}I_{10}}$ & $1$ & $c_{1,12}$
\rule[-0.2cm]{0cm}{0.3cm}\\
\hline
\rule{0cm}{0.4cm}
$P_b$ & ${_{P_b}\chi}_{2}(k)$ & $\frac{q}{\sqrt{2}}\phi_1 \phi_2$ & ${_{P_b}I_{2}}$ & $q^2-1$ & $c_{2,0}(i), c_{2,2}(i)$
\rule[-0.2cm]{0cm}{0.1cm}\\
\rule{-0.12cm}{0.1cm}
    & ${_{P_b}\chi}_{3}(k)$  & $\frac{q}{\sqrt{2}}\phi_1 \phi_2$ & ${_{P_b}I_{3}}$ & $q^2-1$ & $c_{2,0}(i), c_{2,2}(i)$
\rule[-0.2cm]{0cm}{0.1cm}\\
\rule{-0.12cm}{0.1cm}
    & ${_{P_b}\chi}_{9}$  & $\frac{q}{\sqrt{2}}\phi_1^2 \phi_2^2 \phi_8$ & ${_{P_b}I_{9}}$ & $1$ & $c_{1,30}$
\rule[-0.2cm]{0cm}{0.1cm}\\
\rule{-0.12cm}{0.1cm}
    & ${_{P_b}\chi}_{10}$  & $\frac{q}{\sqrt{2}}\phi_1^2 \phi_2^2 \phi_8$ & ${_{P_b}I_{10}}$ & $1$ & $c_{1,30}$
\rule[-0.2cm]{0cm}{0.5cm}\\
\hline
\end{tabular}
\end{small}
\end{center}

\vspace{0.8cm}


\noindent Table~A.15 \textit{The irreducible characters of the
  chain normalizers of defect $24n+12$. (The parameter set
  ${_GI^{ss}}$ parameterizes the semisimple irreducible characters of
  $G$.)} 

\bigskip

\begin{center}
\begin{small}
\begin{tabular}{|l|l|l|c|l|l|} \hline
\rule{0cm}{0.3cm}
Group & Character & Degree & Parameter &
Number & Class
\rule[-0.1cm]{0cm}{0.3cm}\\
\hline \cline{1-3} \hline
\rule{0cm}{0.3cm}
$G$ &  & & ${_GI^{ss}}$ & $q^4$ &  
\rule[-0.2cm]{0cm}{0.5cm}\\
\hline
\rule{0cm}{0.4cm}
$B$ & ${_B\chi}_{1}(k,l)$ & $1$ & ${_BI_{1}}$ & $(q^2-1)^2$ & $c_{10,0}(i,j)$ 
\rule[-0.2cm]{0cm}{0.1cm}\\
\rule{-0.12cm}{0.1cm}
    & ${_B\chi}_{2}(k)$ & $\phi_1 \phi_2$ & ${_BI_{2}}$ & $q^2-1$ & $c_{2,0}(i)$
\rule[-0.2cm]{0cm}{0.1cm}\\
\rule{-0.12cm}{0.1cm}
    & ${_B\chi}_{5}(k)$ & $\phi_1 \phi_2$ & ${_BI_{5}}$ & $q^2-1$ & $c_{3,0}(i)$
\rule[-0.2cm]{0cm}{0.1cm}\\
\rule{-0.12cm}{0.1cm}
    & ${_B\chi}_{8}$ & $\phi_1^2 \phi_2^2$ & ${_BI_{8}}$ & $1$ & 
\rule[-0.2cm]{0cm}{0.5cm}\\
\hline
\rule{0cm}{0.4cm}
$P_a$ & ${_{P_a}\chi}_{1}(k)$ & $1$ & ${_{P_a}I_{1}}$ & $q^2-1$ & 
$c_{2,0}(i)$
\rule[-0.2cm]{0cm}{0.3cm}\\
\rule{0cm}{0.4cm}
      & ${_{P_a}\chi}_{3}(k,l)$ & $\phi_4$ & ${_{P_a}I_{3}}$ & $\frac{q^4-3q^2+2}{2}$ & $c_{7,0}(i,j)$
\rule[-0.2cm]{0cm}{0.3cm}\\
\hline
\end{tabular}
\end{small}
\end{center}

\newpage

\noindent Table~A.16 \textit{Number of fixed points of $H =
  \la \alpha^t \ra$ on parameter sets of the irreducible
  characters. The unions of parameter sets in this table are disjoint unions.}

\bigskip

\begin{center}
\begin{small}
\begin{tabular}{|l|c|} \hline
\rule{0cm}{0.3cm}
Parameter set $I$ \hspace{6cm} & Number of fixed \\
& points $|C_I(H)|$ 
\rule[-0.1cm]{0cm}{0.3cm}\\
\hline \cline{1-2} \hline
\rule{0cm}{0.3cm}
\hspace*{-0.22cm} ${_GI^{ss}}$ & $2^{2t}$ \\
${_GI_{2}} \cup {_GI_{3}}$ & $2$ \\
${_GI_{4}} \cup {_GI_{30}}$ & $2$ \\
${_GI_{5}} \cup {_GI_{6}} \cup {_GI_{8}} \cup {_GI_{9}} \cup {_GI_{11}} \cup {_GI_{12}} \cup {_GI_{13}} \cup {_GI_{14}}$ & $8$ \\
${_GI_{7}} \cup {_GI_{10}}$ & $2$ \\
${_GI_{17}}$ & $1$ \\
${_GI_{18}}$ & $1$ \\
${_GI_{19}} \cup {_GI_{20}}$ & $2$ \\
${_GI_{23}} \cup {_GI_{43}} \cup {_GI_{47}}$ & $2^t-1$ \\
${_GI_{24}} \cup {_GI_{44}} \cup {_GI_{48}}$ & $2^t-1$ \\
${_GI_{25}} \cup {_GI_{45}} \cup {_GI_{49}}$ & $2^t-1$ \\
${_GI_{27}} \cup {_GI_{33}}$ & $2^t-2$ \\
${_GI_{31}}$ & $1$ \\
\hline
${_BI_{1}}$ & $(2^t-1)^2$ \\
${_BI_{2}}$ & $2^t-1$ \\
${_BI_{3}}$ & $2^t-1$ \\
${_BI_{4}}$ & $2^t-1$ \\
${_BI_{5}}$ & $2^t-1$ \\
${_BI_{6}} \cup {_BI_{7}}$ & $2$ \\
${_BI_{8}}$ & $1$ \\
${_BI_{9}}$ & $2^t-1$ \\
${_BI_{10}}$ & $1$ \\
${_BI_{11}}$ & $2^t-1$ \\
${_BI_{12}}$ & $2^t-1$ \\
${_BI_{13}} \cup {_BI_{14}}$ & $2$ \\
${_BI_{15}} \cup {_BI_{16}} \cup {_BI_{17}} \cup {_BI_{18}} \cup
{_BI_{19}} \cup {_BI_{20}} \cup {_BI_{21}} \cup {_BI_{22}}$ & $8$ \\ 
${_BI_{23}}$ & $2^t-1$ \\
${_BI_{26}}$ & $2^t-1$ \\
${_BI_{27}}$ & $2^t-1$ \\
${_BI_{28}} \cup {_BI_{29}}$ & $2$ \\
${_BI_{39}}$ & $2^t-1$ \\
${_BI_{42}}$ & $2^t-1$ \\
${_BI_{43}}$ & $1$ \\
${_BI_{44}}$ & $2^t-1$ \\
${_BI_{45}}$ & $2^t-1$ \\
${_BI_{51}}$ & $2^t-1$ \\
${_BI_{52}}$ & $2^t-1$ \\
${_BI_{53}} \cup {_BI_{54}}$ & $2$ \\
${_BI_{55}} \cup {_BI_{56}} \cup {_BI_{57}} \cup {_BI_{58}}$ & $4$ \\
\hline
\end{tabular}
\end{small}
\end{center}

\newpage

\noindent Table~A.16 (cont.) \textit{Number of fixed points of $H =
  \la \alpha^t \ra$ on parameter sets of the irreducible
  characters.}

\bigskip

\begin{center}
\begin{small}
\begin{tabular}{|l|c|} \hline
\rule{0cm}{0.3cm}
Parameter set $I$ & Number of fixed \\
& points $|C_I(H)|$ 
\rule[-0.1cm]{0cm}{0.3cm}\\
\hline \cline{1-2} \hline
${_{P_a}I_{1}}$ & $2^t-1$ \\
${_{P_a}I_{2}}$ & $2^t-1$ \\
${_{P_a}I_{3}} \cup {_{P_a}I_{4}}$ & $(2^t-1)^2$ \\
${_{P_a}I_{5}}$ & $2^t-1$ \\
${_{P_a}I_{6}}$ & $1$ \\
${_{P_a}I_{7}}$ & $2^t-1$ \\
\hline
\end{tabular}
\end{small}
\end{center}

\newpage

\section*{Appendix B}

The irreducible characters of the simple Ree groups ${^2F}_4(q^2)$
were computed by G.~Malle, see \cite{MalleUni2F4} and the CHEVIE
library \cite{CHEVIE}. However, it seems that $10$ families of
irreducible characters of ${^2F}_4(q^2)$ have never been 
published. We give a brief sketch of the construction of these
irreducible characters using techniques from Deligne-Lusztig theory.
We proceed very much along the same line as F.~L\"ubeck~\cite{Luebeck}. 
For many of our calculations we used CHEVIE and GAP~\cite{GAP4}.
We have implemented the completed character table of ${^2F}_4(q^2)$ as
a generic character table in CHEVIE, 
and whenever we write $_G\chi_i(k,\dots)$ we mean the character of
${^2F}_4(q^2)$ which is labeled by $\chi_i(k, \dots)$ in this CHEVIE
table. The tables in this Appendix might also be of independent
interest; Tables~B.9 and B.10 are also used in~\cite{Himstedt2F4Decomp}.


\subsection*{Notation and setup}
\label{subsec:notaset}

We fix an integer $n > 0$ and set $\theta = 2^n$ and $q :=
\sqrt{2^{2n+1}} \in \R_{>0}$. Let $\F_{q^2}$ be a finite field with
$q^2$ elements, $\F$ its algebraic closure. To describe our
calculations, we fix a simple algebraic group $\G$ of adjoint type
$F_4$ defined over $\F$. Up to isomorphism of algebraic groups, $\G$
is determined by a root datum $(X, \Phi, Y, \Phi^\vee)$, see
Section~1.9 in \cite{Carter2}:

Let $X = Y = \Z^4$ be free abelian groups of rank $4$ and let 
$\la \cdot, \cdot \ra: X \times Y \rightarrow \Z$, 
$(x,y) \mapsto xy^{tr}$ be the natural pairing. Let $\{\hat{e}_1,
\hat{e}_2, \hat{e}_3, \hat{e}_4 \}$ and $\{e_1, e_2, e_3, e_4\}$ be
the standard bases of $X$ and $Y$, respectively. We can describe the
root system $\Phi$ and the coroot system $\Phi^\vee$ by giving a set
$\Delta = \{r_1, r_2, r_3, r_4\} \subseteq X$ of simple roots and a
set $\Delta^\vee = \{r_1^\vee, r_2^\vee, r_3^\vee, r_4^\vee\}
\subseteq Y$ of simple coroots. We write $r_1, r_2, r_3, r_4$ as rows
of a matrix $A$ and $r_1^\vee, r_2^\vee, r_3^\vee, r_4^\vee$ as rows
of a matrix $A^\vee$:
\[
A = \left(\begin{array}{rrrr}
1 & 0 & 0 & 0\\
0 & 1 & 0 & 0\\
0 & 0 & 1 & 0\\
0 & 0 & 0 & 1
\end{array}\right) \quad \text{and} \quad 
A^\vee = \left(\begin{array}{rrrr}
 2 & -1 &  0 &  0\\
-1 &  2 & -1 &  0\\
 0 & -2 &  2 & -1\\
 0 &  0 & -1 &  2
\end{array}\right).
\]
So, the Cartan matrix $\la r_i, r_j^\vee \ra$ is equal to $A^\vee$ and
we see that the Dynkin diagram of $\Phi$ is

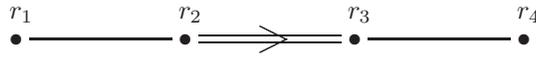
\setlength{\unitlength}{0.9mm}
\begin{center}
\begin{picture}(85,13)
\thinlines
\put(5,5){\circle*{1.5}}
\put(30,5){\circle*{1.5}}
\put(55,5){\circle*{1.5}}
\put(80,5){\circle*{1.5}}
\put(7,5){\line(1,0){21}}
\put(32,5.5){\line(1,0){21}}
\put(32,4.5){\line(1,0){21}}
\put(57,5){\line(1,0){21}}
\put(41,7){\line(2,-1){4}}
\put(41,3){\line(2,1){4}}
\put(4,8){$r_1$}
\put(29,8){$r_2$}
\put(54,8){$r_3$}
\put(79,8){$r_4$}
\end{picture}
\end{center}

To define a Frobenius morphism $F: \G \rightarrow \G$ it is enough to
describe the induced action of $F$ on $X$ in the form $q \cdot F_0$
where $q$ is as above and $F_0$ is the automorphism of the vector space 
$V = X \otimes_\Z \R$ given by
$
(x_1, x_2, x_3, x_4) \mapsto (\frac{x_4}{\sqrt{2}},
\frac{x_3}{\sqrt{2}}, \frac{2x_2}{\sqrt{2}}, \frac{2x_1}{\sqrt{2}}). 
$
This determines ${^2F}_4(q^2) := \G^F$ up to isomorphism.

A realization of the group $\G$ with Frobenius map $F$ is described in
\cite{HimstedtHuang2F4Borel} and \cite{ShinodaClasses2F4} (and
originally in \cite{Ree2F4}):
Let $V$ be a Euclidean vector space with scalar product $( \cdot ,
\cdot )$, let $\{\epsilon_1,\epsilon_2,\epsilon_3,\epsilon_4\}$ be an
orthonormal basis of~$V$ and let $\Phi$ be the set consisting of the
$48$ vectors  
\[
\epsilon_i \, , \, \epsilon_i+\epsilon_j \, , \,\,
\frac{1}{2}(\epsilon_i+\epsilon_j+\epsilon_k+\epsilon_l), 
\]
where $i,j,k,l \in \{\pm 1, \pm 2, \pm 3, \pm 4\}$, $|i|,|j|,|k|,|l|$
are different and $\epsilon_{-i} = -\epsilon_i$ for all~$i$. The set $\Phi$ is a
root system of type $F_4$ and the set $\Delta := \{r_1, r_2, r_3,
r_4\}$ with the simple roots 
$
r_1 := \epsilon_2 - \epsilon_3 \, , \, r_2 := \epsilon_3 - \epsilon_4
\, , \, r_3 := \epsilon_4 \, , \, r_4 := \frac{1}{2}(\epsilon_1 -
\epsilon_2 - \epsilon_3 - \epsilon_4) 
$
is a basis of $\Phi$. For each $r \in \Phi$, we define the
corresponding coroot $r^\vee := \frac{2r}{(r,r)}$. Then, 
$\Phi^\vee := \{r^\vee \, | \, r \in \Phi\}$ is a root system of type
  $F_4$ with basis $\{r_1^\vee, r_2^\vee, r_3^\vee, r_4^\vee\}$.
The Weyl group $\W$ is the group generated by $w_{r_1}, w_{r_2},
w_{r_3}, w_{r_4}$, the reflections at the hyperplane orthogonal to the
simple roots. 

Let $\G$ be the simple algebraic group over $\F$ of adjoint type with
Frobenius morphism $F: \G \rightarrow \G$, maximal $F$-stable
torus~$\T$ and root system $\Phi$ which is described in
\cite{HimstedtHuang2F4Borel} and~\cite{ShinodaClasses2F4}. In
\cite[Section 1.9]{Carter2} it is shown how the character group $X$ and
the cocharacter group $Y$ of the torus $\T$ are endowed with a natural
pairing $\la \cdot, \cdot \ra$ and how $\Phi$ and $\Phi^\vee$ can be
naturally embedded into $X$ and $Y$, respectively. This gives a root
datum $(X, \Phi, Y, \Phi^\vee)$ with the properties described at the
beginning of this section.

The elements of the maximal torus $\T$ can be parameterized as follows:
There is a natural isomorphism between the abelian groups $\T$ and 
$\Hom(X, \F^\times)$ (see \cite[Section~1.11 and 
Proposition~3.1.2 (i)]{Carter2}) and we write $h(z_1, z_2, z_3, z_4)$ for 
the element of $\T$ corresponding to $\chi \in \Hom(X, \F^\times)$ with 
$\chi(\epsilon_i) = z_i$ ($i=1,2,3,4$). 

We mention that there is an alternative
parameterization of the elements of $\T$ (which in fact we use for our
GAP-programs): By \cite[Proposition~3.1.2 (ii)]{Carter2}, we have  
$\T \cong Y \otimes \F^\times$ as abelian groups. Every element of 
$Y \otimes \F^\times$ can be written uniquely as 
$\sum_{i=1}^4 e_i \otimes \lambda_i$ with $\lambda_i \in \F^\times$ and we 
write $(\lambda_1, \lambda_2, \lambda_3, \lambda_4)$ for the corresponding
element of $\T$. The maps transforming one parameterization into the
other are given by: 
\begin{eqnarray*}
(\lambda_1, \lambda_2, \lambda_3, \lambda_4) & \mapsto & h(\lambda_1
\lambda_2^2 \lambda_3^3 \lambda_4^2, \lambda_1 \lambda_2 \lambda_3, \lambda_2
\lambda_3, \lambda_3) \quad \text{and} \\
h(z_1, z_2, z_3, z_4) & \mapsto & (z_2 z_3^{-1}, z_3 z_4^{-1}, z_4, (z_1 z_2^{-1}
z_3^{-1} z_4^{-1})^{1/2}). 
\end{eqnarray*}

The Weyl group $\W$ is isomorphic with $N_\G(\T)/\T$ and so acts on
$\T$ by conjugation which also induces an action of $\W$ on $X$ and
on $Y$, see \cite{Luebeck}. The action of~$\W$ on $X$ and $Y$ can be
determined by the formulas in \cite[p. 19]{Carter2}. Once we know the
action of $\W$ on $X$, we can also describe the action on $\T$ via 
$\T \cong \Hom(X, \F^\times)$. These actions are given in Table~B.1.

Finally, we mention that \cite[Section 2]{HimstedtHuang2F4Borel}
provides generators and relations for the group $\G$. Using GAP
programs written by Christoph K\"ohler and the first author in the
language of the CHEVIE package, this enables us to carry out explicit
computations in the groups $\G$ and $\G^F = {^2F}_4(q^2)$.  

\subsection*{Maximal tori}
\label{subsec:tori}

Every semisimple element $s \in \G^F$ is contained in an $F$-stable
maximal torus $\tilde{\T}$ of $\G$, see \cite[Corollary 3.16]{DigneMichel}.
Suppose $\tilde{\T}$ is such an $F$-stable maximal torus of
$\G$. We choose some $g \in \G$ such that $\tilde{\T}^g = \T$ and we
have $n := g^{-1}F(g) \in N_G(\T)$. The action of $F$ on $\tilde{\T}$
corresponds to the action of $(Fn^{-1}): \T \rightarrow \T, t \mapsto
{^nF}(t)$ (in the sense $F(\tilde{t})^g = {^nF}(\tilde{t}^g)$). Since
$\T$ is abelian we also write $(Fw^{-1}): \T \rightarrow \T, t \mapsto
{^wF}(t)$ where $w := n\T \in \W$. The map $\tilde{\T} \mapsto w :=
n\T$ induces a bijection between the $\G^F$-conjugacy classes of the
$F$-stable maximal tori of $\G$ and the $F$-conjugacy classes of the
Weyl group $\W$, for details see \cite[Section 3.3]{Carter2}. A set of
representatives $w_i$ for the $F$-conjugacy classes of $\W$ and the
corresponding subgroups $T^{Fw^{-1}}$ of fixed points were computed in
\cite[\S 3]{ShinodaClasses2F4}. A list of these representatives, the
order of the corresponding $F$-centralizers and the maximal tori
$T^{(Fw_i^{-1})}$ are given in Table~B.2. Parameters for the elements
of $T^{(Fw_i^{-1})}$ can be found in Table~B.3, where we describe
field elements using the notation in \cite[Table 5]{HimstedtHuang2F4Borel}.

Let $(X, \Phi, Y, \Phi^\vee)$ be the root datum of $\G$ described at
the beginning of this appendix. So the dual group~$\G^*$ is a
connected reductive group with root datum $(Y, \Phi^\vee, X, \Phi)$. 
The action of $F$ on $X$ and $Y$ gives rise to a Frobenius morphism
$F^*: \G^* \rightarrow \G^*$ which determines the finite group
$\G^{*F^*}$ up to isomorphism.  

A maximal $F^*$-stable torus of $\G^*$ is isomorphic with 
$X \otimes_\Z \F^\times$ and, via the isomorphism $\varphi_1$ from
\cite[Section~2]{HimstedtHuang2F4Borel}, also isomorphic with
$X \otimes_\Z \Q_{p'}/\Z$. Using the $\Z$-basis $\{\hat{e}_1,
\hat{e}_2, \hat{e}_3, \hat{e}_4 \}$ of $X$, every element of 
$X \otimes_\Z \Q_{p'}/\Z$ can be written as
\[
s(\mu_1, \mu_2, \mu_3, \mu_4) := \sum_{i=1}^4 \hat{e}_i \otimes \mu_i
\text{   with   } \mu_i \in \Q_{p'}/\Z.
\]
The action of $F^*$ on $\T^*$ is induced by the action of $F$ on $X$
\[
F^*(s(\mu_1, \mu_2, \mu_3, \mu_4)) = s(\theta \mu_4, \theta \mu_3,
2\theta \mu_2, 2\theta \mu_1).
\]
Similarly, the action of $\W$ on $X$ induces an action of $\W$ on
$\T^*$. This enables us to compute the fixed point groups 
$\T_i^{*F^*} := (X \otimes_\Z \Q_{p'}/\Z)^{(w_i F)}$, where $w_i$ runs
through the representatives $w_i$ for the $F$-conjugacy classes of $\W$. 
In this way, we get the parameterization of the maximal tori of
$\G^{*F^*}$ which is described in Table~B.4.

For each $i$, the maximal torus $\T_i^{*F^*}$ of the dual group is
isomorphic with the group of linear characters of the torus $\T_i^F$
of $\G^F$ and this isomorphism can be calculated explicitly as 
in \cite[p. 23]{Luebeck}. Applying this isomorphism
we get the results in Table~B.5. These linear characters will be used 
to parameterize the Lusztig series of irreducible characters of $\G^F$.

\subsection*{Conjugacy classes}
\label{subsec:conjcl}

We say that two semisimple elements $s_1, s_2 \in \G^F$ are in the
same \textsl{class type} if their centralizers $C_\G(s_1)$ and
$C_\G(s_2)$ are $\G^F$-conjugate. A~subset $\Psi \subseteq \Phi$ is
called a \textsl{closed subsystem} if $\Z \Psi \cap \Phi = \Psi$.
Let $M$ be the set of all pairs $(\Pi, w)$ such that $\Pi$ is a set of
simple roots for a closed subsystem of $\Phi$ and $w$ is an
element of $\W$ such that $\Pi$ is $(Fw^{-1})$-invariant. We call 
$(\Pi_1, w_1), (\Pi_2, w_2) \in M$ equivalent if and only if there is
some $v \in \W$ such that $\Pi_1 v = \Pi_2$ and 
$v^{-1} w_1 F(v) = w_2$ and let $\bar{M}$ be the set of
equivalence classes. The semisimple class types of $\G^F$ can be
parameterized by a subset of $\bar{M}$, see 
\cite[Section III.4]{Koehler}. 

A set of representatives for the semisimple conjugacy classes of
$\G^F$ was computed by K.~Shinoda, see \cite[Table IV]{ShinodaClasses2F4}, and 
is given in Table~B.7. The semisimple class types of $\G^F$ are listed
in Table~B.6. The data in these tables can be obtained by the methods
described in \cite[Section III.4]{Koehler} and 
\cite[Section 4.1]{Luebeck}. The notation for the positive roots $r_i$
in Table~B.6 is defined in \cite[Section~2]{HimstedtHuang2F4Borel}.
Analogously, we can compute the semisimple class types and
representatives for the semisimple conjugacy classes of $\G^{*F^*}$
given in Tables~B.9 and B.10. We just have to use the dual torus $\T^*$
instead of $\T$ and the coroots instead of the roots.

Representatives for the unipotent and mixed conjugacy classes were
calculated by K.~Shinoda, see \cite[Tables II and V]{ShinodaClasses2F4}
and are given in Table~B.8. They can be obtained by explicit
computations in $\G^F$ using the relations in 
\cite[Tables 1-4]{HimstedtHuang2F4Borel}. 

\subsection*{Lusztig series and Deligne-Lusztig characters}
\label{subsec:dellin}

Let $\T_i$ be one of the $F$-stable maximal tori of $\G$ in
Table~B.3 and $\theta$ a linear character of $\T_i^F$. Associated with
the pair $(\T_i, \theta)$ is a generalized character 
$R_{\T_i}^\G \theta$ of $\G^F$, called \textsl{Deligne-Lusztig character}.
In fact, $R_{\T_i}^\G \theta$ only depends on the $\G^F$-conjugacy
class of $(\T_i, \theta)$.
Every irreducible character of $\G^F$ occurs as a constituent of some 
$R_{\T_i}^\G \theta$. Every pair $(\T_i, \theta)$ can be identified
with a pair $(\T_i^*, s)$ where $\T_i^*$ is one of the $F^*$-stable
maximal tori of the dual group $\G^*$ in Table~B.4 and 
$s \in \T_i^{*F^*}$. Two Deligne-Lusztig \pagebreak 
characters $R_{\T_i}^\G \theta_1$
and $R_{\T_j}^\G \theta_2$ can have a common constituent only if the
corresponding semisimple elements $s_1, s_2 \in \G^{*F^*}$ are
$\G^{*F^*}$-conjugate (here we use that the centralizers
$C_{\G^*}(s_1)$, $C_{\G^*}(s_2)$ are connected). This induces an
equivalence relation on $\Irr(\G^F)$ in a natural way. The equivalence
classes $\mathcal{E}(\G^F, (s))$ are parameterized by the semisimple
conjugacy classes $(s)$ of $\G^{*F^*}$ and are called 
\textsl{Lusztig series}. We say, that a class function is contained in
the Lusztig series $\mathcal{E}(\G^F, (s))$ if and only if it is
in the $\C$-span of the irreducible characters in $\mathcal{E}(\G^F, (s))$.

We say that two Lusztig series $\mathcal{E}(\G^F, (s_1))$, 
$\mathcal{E}(\G^F, (s_2))$ are \textsl{of the same type} if and only
if $s_1, s_2$ belong to the same semisimple class type of $\G^{*F^*}$.
So, Tables~B.9 and B.10 give an overview over the Lusztig series 
of $\G^F$. The irreducible characters of $\G^F$ we want to construct
are those in the Lusztig series of types $g_8$, $g_{10}$, $g_{13}$, 
$g_{14}$.

The values of the Deligne-Lusztig characters can be computed using the
formula given in \cite[Satz 2.1 (b)]{Luebeck}: Fix an $F$-stable
maximal torus $\T_i$ of $\G$ and a linear character $\theta$ of $\T_i^F$.
Let $g$ be an element of $\G^F$ and $g = su = us$ with $s \in \G^F$
semisimple and $u \in \G^F$ unipotent be the Jordan decomposition of
$g$. We write $\mathbf{C} := C_\G(s)$ for the centralizer of $g$ in
$\G$. Note that in our situation, $\mathbf{C}$ is connected (because
$\G$ is simply connected). Let $\T_{i1}, \dots, \T_{ir}$ be a set of
representatives for the $\mathbf{C}^F$-conjugacy classes of $F$-stable
maximal tori in $\mathbf{C}$ which are contained in the
$\G^F$-conjugacy class of $\T_i$ and let $s_{ir1}, \dots, s_{irl_r}$
for $1 \le r \le k$ be the different elements of the form $s^x$ with
$x \in \G^F$ such that ${^x\T}_i$ is in the $\mathbf{C}^F$-conjugacy
class of $\T_{ir}$. Then 
\begin{equation}\label{eq:RTtheta}
(R_{\T_i}^\G \theta)(g) = \sum_{r=1}^k \left(\sum_{t=1}^{l_r}
\theta(s_{irt})\right) Q_{T_{ir}}^\mathbf{C}(u)
\end{equation}
where $Q_{T_{ir}}^\mathbf{C}$ is the Green function of $\mathbf{C}$
with respect to $\T_{ir}$. We will only need $R_{\T_i}^\G \theta$ for
$i=3,4,5,6,7$ and arbitrary $\theta$. The values of 
$R_{\T_i}^\G \theta$ for $i=3,4,5$ can be read off from the character
table of ${^2F}_4(q^2)$ in CHEVIE, see the characters
$-{_G\chi}_{35}(k)$, $-{_G\chi}_{36}(k)$, ${_G\chi}_{37}(k,l)$ in the
CHEVIE table. 

So we only have to compute the values of $R_{\T_6}^\G \theta$ and
$R_{\T_7}^\G \theta$. For these Deligne-Lusztig characters, the right
hand side of (\ref{eq:RTtheta}) can only be non-zero if $s$ is of
class type $h_1$, $h_8$, $h_{10}$, $h_{13}$ or $h_{14}$. For these
classes, $\mathbf{C}$ has Dynkin type ${^2F}_4$, ${^2B}_2$ or $A_0$. So
the corresponding Green functions are known, see the CHEVIE library of
Green functions. 

The elements $s_{irt}$ can be computed from the root datum of $\G$ as
it is described in \cite[p.~29, (3)]{Luebeck}. It turns out
that in the computations for $R_{\T_6}^\G \theta$ and $R_{\T_7}^\G \theta$ 
we always have $r=k=1$. This is all information which is needed to
compute the values of $R_{\T_6}^\G \theta$ and $R_{\T_7}^\G \theta$
via formula~(\ref{eq:RTtheta}). 

The Lusztig series $g_{13}$ and $g_{14}$ correspond to regular
semisimple elements of $\G^*$. This implies that $R_{\T_6}^\G \theta$
and $R_{\T_7}^\G \theta$ are (up to possibly a sign) irreducible
characters of $\G^F$. Considering the value on the identity element,
we see that in our situation this sign is $+1$. This gives us the
irreducible characters $_G\chi_{50}(k,l)$ and $_G\chi_{51}(k,l)$ of
$\G^F$ and completes the construction of all irreducible characters in 
the Lusztig series of types $g_{13}$, $g_{14}$. However, since some
values of these irreducible characters involve sums of up to $96$
roots of unity, we cannot print their values here in this article. 

\subsection*{Construction of irreducible characters of type
  \texorpdfstring{$g_8$}{g8} and \texorpdfstring{$g_{10}$}{g10}} 
\label{subsec:irrG}

Our main tool to construct the irreducible characters in the Lusztig
series of types $g_8$ and $g_{10}$ is Lusztig's Jordan decomposition
of characters: Let $\T_i$ be an $F$-stable maximal torus of $\G$ and
$\theta$ a linear character of $\T_i^F$ and let $(\T_i^*, s)$ be a
pair corresponding to $(\T_i, \theta)$. Then, there  is a bijection 
\begin{equation*} \label{eq:jord}
\mathcal{L}': \mathcal{E}(\G^F, (s)) \rightarrow
\mathcal{E}(C_{\G^*}(s)^{F^*}, (1)).
\end{equation*}
and for $\chi \in \mathcal{E}(\G^F, (s))$ we have the following
identity of scalar products:
\begin{equation} \label{eq:jordanscal}
(\chi, R_{\T_i}^\G \theta)_{\G^F} = \epsilon_\G \cdot \epsilon_{C_{\G^*}(s)}
\cdot (\mathcal{L}'(\chi), R_{\T_i^*}^{C_{\G^*}(s)^*} \mathbf{1})_{C_{\G^*}(s)^{F^*}}
\end{equation}
where $\epsilon_\G, \epsilon_{C_{\G^*}(s)} \in \{1, -1\}$ and
$\mathbf{1}$ is the trivial character of $\T_i^{*F^*}$. The signs 
$\epsilon_\G, \epsilon_{C_{\G^*}(s)}$ are defined in \cite[Section
  6.5]{Carter2} and can be computed directly from their
definition. Alternatively, since ${^2F}_4(q^2)$ is isomorphic with its
dual group, they can also be read off from the values of the Steinberg
character, see \cite[Theorem 6.5.9]{Carter2}.

We describe the construction of the irreducible characters of $\G^F$
in the Lusztig series of type $g_8$. The calculations for the series
of type $g_{10}$ are analogous.
Let~$g_8(k)$ be a representative for the semisimple conjugacy classes
of $\G^{F^*}$ as in Table~B.10 and let $C^* := C_{\G^*}(g_8(k))$.
By Table~B.9 we know that $C^{*F^*}$ has Dynkin type~${^2B}_2$. So the
Jordan decomposition of characters implies that the Lusztig series 
$\mathcal{E}(\G^F, (g_8(k)))$ contains exactly $4$ irreducible
characters, one of them corresponding to the trivial character, two
corresponding to the cuspidal characters ${^2B}_2[a]$, ${^2B}_2[b]$
and one corresponding to the Steinberg character for groups of type
${^2B}_2$. We denote these characters, in this order, by
$_G\chi_{8,1}(k)$, $_G\chi_{8,2}(k)$, $_G\chi_{8,3}(k)$, $_G\chi_{8,4}(k)$.

As described in \cite[Section 4.1 (2), (3) and Section 6]{Luebeck} 
we can compute representatives for the $\G^{*F^*}$-conjugacy classes 
$g_8(k)$ in the various tori $\T_i^*$ and see that the only maximal
tori containing such representatives are $\T_3^*$, $\T_5^*$ and
$\T_6^*$. We obtain that, for fixed $k$, the Deligne-Lusztig
characters $-{_G\chi}_{35}((2\theta^3-\theta) k)$, 
${_G\chi}_{37}(0, k)$, ${_G\chi}_{50}(0, k)$ all correspond to the  
$\G^{*F^*}$-conjugacy class $g_8(k)$. Using CHEVIE, we
can easily compute the norm $(-{_G\chi}_{35}((2\theta^3-\theta) k),
-{_G\chi}_{35}((2\theta^3-\theta) k) = 2$ and, using the scalar product property
(\ref{eq:jordanscal}), we can compute the scalar products 
$(-{_G\chi}_{35}((2\theta^3-\theta) k), {_G\chi}_{8,1}(k)) = 
(-{_G\chi}_{35}((2\theta^3-\theta) k), {_G\chi}_{8,4}(k)) = -1$. 
So that we get
\[
-{_G\chi}_{35}((2\theta^3-\theta) k) = -{_G\chi}_{8,1}(k) - {_G\chi}_{8,4}(k).
\]
Analogously, we obtain
\begin{eqnarray*}
{_G\chi}_{37}(0, k) & = & -{_G\chi}_{8,1}(k)-{_G\chi}_{8,2}(k)-{_G\chi}_{8,3}(k)+{_G\chi}_{8,4}(k),\\
{_G\chi}_{50}(0, k) & = & -{_G\chi}_{8,1}(k)+{_G\chi}_{8,2}(k)+{_G\chi}_{8,3}(k)+{_G\chi}_{8,4}(k).
\end{eqnarray*}
So, we see that the class function
\[
f_8(k) := {_G\chi}_{8,2}(k)-{_G\chi}_{8,3}(k)
\]
is orthogonal to the $\C$-vector space spanned by the Deligne-Lusztig
characters. By inverting the coefficient matrix of the above linear
equations, we can now write ${_G\chi}_{8,1}(k), \dots,
{_G\chi}_{8,4}(k)$ as linear combinations of the Deligne-Lusztig
characters and $f_8(k)$ as follows:
\[
\left(\begin{array}{c}
{_G\chi}_{8,1}(k)\\
{_G\chi}_{8,2}(k)\\
{_G\chi}_{8,3}(k)\\
{_G\chi}_{8,4}(k)
\end{array}\right)
=
\frac{1}{4}
\left(\begin{array}{rrrr}
-2 & -1 & -1 &  0\\
 0 & -1 &  1 &  2\\
 0 & -1 &  1 & -2\\
-2 &  1 &  1 &  0
\end{array}\right)
\cdot
\left(\begin{array}{c}
-{_G\chi}_{35}((2\theta^3-\theta) k)\\
{_G\chi}_{37}(0, k)\\
{_G\chi}_{50}(0, k)\\
f_8(k)
\end{array}\right)
\]
This already gives us all values of the irreducible characters 
${_G\chi}_{42}(k) := {_G\chi}_{8,1}(k)$ and
${_G\chi}_{45}(k) := {_G\chi}_{8,4}(k)$. To determine the values of
the non-uniform characters ${_G\chi}_{43}(k)$, ${_G\chi}_{44}(k)$ we
have to compute the values of the class functions $f_8(k)$.

\subsection*{The class functions \texorpdfstring{$f_8(k)$}{f8(k)}}

To determine the values of the class function $f_8(k)$ we use
arguments similar to those in \cite[Section 8]{Luebeck}. We proceed in
several steps.

\smallskip

\noindent Step 1: Class functions constructed from Levi subgroups.

\noindent Fix a representative $h \in \G^F$ for the semisimple
conjugacy classes of type $h_8$ as in Table~B.7 and let $\LL$ be its
centralizer in $\G$. Since the order of the representative is not
divisible by a bad prime for $F_4$, the centralizer $\LL$ is an
$F$-stable Levi subgroup of~$\G$. 

Since $\LL^F$ has Dynkin type ${^2B}_2$ it has two cuspidal unipotent
irreducible characters ${^2B}_2[a]$, ${^2B}_2[b]$ and their difference
$f := {^2B}_2[a] - {^2B}_2[b]$ is a unipotent class function of~$\LL^F$
which is orthogonal to all Deligne-Lusztig characters of
$\LL^F$. Using the notation from Table~B.8, the values of $f$ are:
\begin{center}
\begin{tabular}{|l|c|c|} \hline
\rule{0cm}{0.3cm}
Class representative & $h_8(i) x_8(1) x_{16}(1) x_{21}(1)$ & $h_8(i)
x_8(1) x_{16}(1) x_{24}(1)$ 
\rule[-0.1cm]{0cm}{0.3cm}\\
\hline \cline{1-3} \hline
\rule{0cm}{0.3cm}
Value of $f$ & $\sqrt{2}q \epsilon_4$ &$-\sqrt{2}q \epsilon_4$ \\
\hline
\end{tabular}
\end{center}
where $\epsilon_4$ is a complex fourth root of unity as in 
\cite[Table 5]{HimstedtHuang2F4Borel}. The values of $f$ on the
remaining conjugacy classes of $\LL^F$ are zero.

Let $\LL^* \subseteq \G^*$ be the dual group of $\LL$. Then, $\LL^*$
is the centralizer of the semisimple elements of $\G^*$ of type
$g_8$. Let 
\[
g_8(k) = s\left(0, \frac{(\theta-1)k}{q^2-\sqrt{2}q+1},
\frac{-k}{q^2-\sqrt{2}q+1}, 0 \right)
\] 
be one of these elements. Since $g_8(k)$ is contained in the center of
$\LL^{F^*}$, we can identify it with a linear character $\lambda(k)$
of $\LL^F$ in a natural way, see \cite[p.~69]{Luebeck}. Using Lusztig
induction, we can define the following class functions on $\G^F$:
\[
\psi(k) := R_{\LL}^\G (f \cdot \lambda(k)).
\]

\smallskip

\noindent Step~2: Connection between $f_8(k)$ and $\psi(k)$.

\noindent From \cite[Lemma 8.1]{Luebeck} and \cite[Lemma 8.2]{Luebeck}
we know that $\psi(k)$ is a generalized character of~$\G^F$ in the
Lusztig series $g_8(k)$ which is orthogonal to all Deligne-Lusztig
characters of $\G^F$. Since the space of class functions in the
Lusztig series $g_8(k)$ which are orthogonal to all Deligne-Lusztig
characters has dimension one, it follows that $\psi(k)$ is a multiple
of $f_8(k)$. Because $\psi(k)$ is a generalized character and $f_8$ has
constituents with multiplicity $\pm 1$ we get that $\psi(k)$ is an
integer multiple of $f_8(k)$.

\smallskip

\noindent Step~3: Many values of $f_8(k)$, $\psi(k)$ are zero.

\noindent The dimension of the space of all class functions which are
orthogonal to all Deligne-Lusztig characters of $\G^F$ has dimension
$q^2+9$ since we get $10$ linearly independent such class functions
from the unipotent characters and $(q^2-2)/2$, $(q^2-\sqrt{2}q)/4$, 
$(q^2+\sqrt{2}q)/4$ such class functions from the Lusztig series of
types $g_2$, $g_8$ and $g_{10}$ respectively. The same argument as in
\cite[p.~71]{Luebeck} shows that $\psi(k)$ vanishes on all conjugacy
classes of $\G^F$ except for possibly the following tuples of classes:
$(c_{1,3}, c_{1,4})$, $(c_{1,7}, c_{1,8}, c_{1,9})$, 
$(c_{1,10}, c_{1,11}, c_{1,12})$,  
$(c_{1,13}, c_{1,14})$, $(c_{1,15}, c_{1,16}, c_{1,17}, c_{1,18})$, 
$(c_{2,2}(i), c_{2,3}(i))$, $(c_{5,2}, c_{5,3}, c_{5,4})$, 
$(c_{8,2}(i), c_{8,3}(i))$, $(c_{10,2}(i), c_{10,3}(i))$.
Additionally, the argument shows that the values of $\psi(k)$
satisfy the following relations:
\begin{eqnarray*}
\psi(k)(c_{1,4}) & = & -\psi(k)(c_{1,3})\\
\psi(k)(c_{1,9}) & = & -\psi(k)(c_{1,8})+\psi(k)(c_{1,7})\\
\psi(k)(c_{1,12}) & = & -2\psi(k)(c_{1,10})-\psi(k)(c_{1,11})\\
\psi(k)(c_{1,14}) & = & -\psi(k)(c_{1,13})\\
\psi(k)(c_{1,18}) & = & -\psi(k)(c_{1,15})-\psi(k)(c_{1,16})-\psi(k)(c_{1,17})\\
\psi(k)(c_{8,3}(i)) & = & -\psi(k)(c_{8,2}(i))
\end{eqnarray*}

\medskip

\noindent Step~4: Values of the remaining irreducible characters.

\noindent We use the formula \cite[Proposition 12.2]{DigneMichel} for 
Lusztig induction to derive some information on the non-zero values of
$\psi(k)$. (Unfortunately, the formula does not give us the values of
$\psi(k)$ explicitly, since we do not know the two-parameter Green
functions occurring in the right hand side of the formula.)

Firstly, we can conclude that $\psi(k)$ vanishes on the classes
$c_{2,2}(i)$, $c_{2,3}(i)$, $c_{5,2}$, $c_{5,3}$, $c_{5,4}$,
$c_{10,2}(i)$, $c_{10,3}(i)$ since the semisimple part of the
representatives for these classes cannot be conjugate to an element of
the center $Z(\LL^F)$. Note that by \cite[p.13]{ShinodaClasses2F4}, we have
$\LL^F = Z(\LL^F) \times (\LL^F)'$. Secondly, we see that 
\begin{eqnarray*}
\psi(k)(c_{8,2}(i)) & = & x_{8,2} \epsilon_4
(\varphi_8''^{qik/\sqrt{2}} + \varphi_8''^{-qik/\sqrt{2}} +
\varphi_8''^{(q/\sqrt{2}-1)ik} + \varphi_8''^{-(q/\sqrt{2}-1)ik}),\\ 
\psi(k)(c_{8,3}(i)) & = & x_{8,3} \epsilon_4
(\varphi_8''^{qik/\sqrt{2}} + \varphi_8''^{-qik/\sqrt{2}} +
\varphi_8''^{(q/\sqrt{2}-1)ik} + \varphi_8''^{-(q/\sqrt{2}-1)ik})
\end{eqnarray*}
for some rational numbers $x_{8,2}$, $x_{8,3}$ not depending on $i,k$.
Thirdly, the formula implies that the values of $\psi(k)$ on the
unipotent conjugacy classes are of the form 
$\psi(k)(c_{1,r}) = \epsilon_4 x_{1,r}$ where $x_{1,r}$ are rational
integers. Considering the $x_{1,r}$, $x_{8,r}$ as indeterminates and
using CHEVIE, we can compute the scalar products of $\psi(k)$ with the
unipotent irreducible characters of $\G^F$. This gives us $\psi(k)$
and hence $f_8(k)$ up to a rational multiple. From the condition
$(f_8(k), f_8(k)) = 2$ we get $f_8(k)$ up to a sign, see Table~B.11.
Note, that the sign of $f_8(k)$ only flips the complex-conjugate irreducible
characters $\chi_{8,2}(k)$ and $\chi_{8,3}(k)$. This gives us the
irreducible characters $\chi_{43}(k)$ and $\chi_{44}(k)$ and thereby
completes the determination of the irreducible characters of $\G^F$ in
the Lusztig series of type $g_8$. The values of the irreducible
characters $\chi_{46}(k), \dots, \chi_{49}(k)$ of $\G^F$ in the 
Lusztig series of type $g_{10}$ can be computed analogously.

Again, for space reasons, we cannot print all values of the
irreducible characters $\chi_{42}(k), \dots, \chi_{49}(k)$. Therefore,
Table B.12 presents only the degrees and those character values which
are used in the verification of Dade's conjecture.

\medskip

\noindent \textbf{Remark:} We have implemented the completed character
table of $\G^F = {^2F}_4(q^2)$ as a generic character table in
CHEVIE. 
To test the correctness of the
table, we have verified the row and column orthogonality
relations. Furthermore, using CHEVIE, we computed the tensor products
of all irreducible characters with the unipotent characters and
computed the scalar products (which have to be nonnegative integers)
of these tensor products with the unipotent characters. 

\newpage



\bigskip
\noindent Table~B.1 \textit{Description of the generators 
$w_{r_1}, w_{r_2}, w_{r_3}, w_{r_4}$ of the Weyl group $\W$ of type
  $F_4$. As automorphisms of the character group $X$ and the
  cocharacter group~$Y$, we describe the generators by their images of
  the chosen $\Z$-bases of $X$ and $Y$, respectively. As elements of
  $N_\G(\T)/\T$, we describe them by giving the torus elements 
  ${^{w_{r_i}}h(z_1,z_2,z_3,z_4)}$.}

\bigskip

\begin{center}
\begin{tabular}{|l|l|l|c|} \hline
\rule{0cm}{0.3cm}
$\W$ & $w_{r_i}$ & Action of $w_{r_i}$ & Matrix (Row convention) 
\rule[-0.1cm]{0cm}{0.3cm}\\
\hline \cline{1-3} \hline
\rule{0cm}{0.3cm}
\hspace{-0.2cm} in $\Aut(X)$ & $w_{r_1}$ & $\hat{e}_1 \mapsto -\hat{e}_1, $ &\\
&& $\hat{e}_2 \mapsto \hat{e}_1+\hat{e}_2,$ & \\
&& $\hat{e}_3 \mapsto \hat{e}_3,$ &\\
&& $\hat{e}_4 \mapsto \hat{e}_4$ &
\raisebox{4ex}[-4ex]{$\left(\begin{array}{rrrr}-1 & 0 & 0 & 0\\1 & 1
    & 0 & 0\\0 & 0 & 1 & 0\\0 & 0 & 0 & 1\end{array}\right)$}
\rule[-0.2cm]{0cm}{0.1cm}\\
\rule{-0.12cm}{0.1cm}
& $w_{r_2}$ & $\hat{e}_1 \mapsto \hat{e}_1+\hat{e}_2, $ &\\
&& $\hat{e}_2 \mapsto -\hat{e}_2,$ & \\
&& $\hat{e}_3 \mapsto \hat{e}_2+\hat{e}_3,$ &\\
\hline
\end{tabular}
\end{center}

\newpage

\bigskip
\noindent Table~B.2 (cf.~\cite[Table III]{ShinodaClasses2F4}) \textit{The
  $F$-conjugacy classes of the Weyl group $\W$ of type~$F_4$. For each
  class, we give a representative $w_i$, the order of the
  $F$-centralizer $C_{\W, F}(w_i)$ and the order of $\T^{Fw_i^{-1}}$.}

\bigskip

\begin{center}
\begin{tabular}{|l|l|l|l|} \hline
\rule{0cm}{0.4cm}
 \hspace{-0.2cm} $\T_i$ \hspace{-0.2cm} & $w_i$ & $|C_{\W, F}(w_i)|$ & $|\T_i^F|$
\rule[-0.1cm]{0cm}{0.4cm}\\
\hline \cline{1-3} \hline
\rule{0cm}{0.6cm}
 \hspace{-0.3cm} $T_1$ \hspace{-0.2cm} & $1$ & $16$ & $(q^2-1)^2$
\rule[-0.2cm]{0cm}{0.6cm}\\
\hline
\rule{-0.12cm}{0.6cm}
 \hspace{-0.2cm} $T_2$ \hspace{-0.2cm} & $w_{r_1}$ & $4$ & $q^4-1$
\rule[-0.2cm]{0cm}{0.6cm}\\
\hline
\rule{-0.12cm}{0.6cm}
 \hspace{-0.2cm} $T_3$ \hspace{-0.2cm} & $w_{r_3}$ & $8$ & $(q^2-1)(q^2-\sqrt{2}q+1)$
\rule[-0.2cm]{0cm}{0.6cm}\\
\hline
\rule{-0.12cm}{0.6cm}
 \hspace{-0.2cm} $T_4$ \hspace{-0.2cm} & $w_{r_2} w_{r_3} w_{r_2}$ & $8$ & $(q^2-1)(q^2+\sqrt{2}q+1)$
\rule[-0.2cm]{0cm}{0.6cm}\\
\hline
\rule{-0.12cm}{0.6cm}
 \hspace{-0.2cm} $T_5$ \hspace{-0.2cm} & $w_{r_1} w_{r_2} w_{r_1}
 w_{r_3} w_{r_2} w_{r_1} w_{r_3} w_{r_2}$ & $16$ & $q^4+1$
\rule[-0.2cm]{0cm}{0.6cm}\\
\hline
\rule{-0.12cm}{0.6cm}
 \hspace{-0.2cm} $T_6$ \hspace{-0.2cm} & $w_{r_1} w_{r_2} w_{r_3} w_{r_2} w_{r_1} w_{r_3}$ & $96$ & $(q^2-\sqrt{2}q+1)^2$
\rule[-0.2cm]{0cm}{0.6cm}\\
\hline
\rule{-0.12cm}{0.6cm}
 \hspace{-0.2cm} $T_7$ \hspace{-0.2cm} & $w_{r_2} w_{r_3} w_{r_2} w_{r_4} w_{r_3} w_{r_2} w_{r_1} w_{r_3}
w_{r_2} \cdot$ && $(q^2+\sqrt{2}q+1)^2$\\
\hline
\end{tabular}
\end{center}

\newpage

\bigskip
\noindent Table~B.3 (cf.~Shinoda \cite[(3.1)]{ShinodaClasses2F4})
\textit{Parameters for the elements of the maximal 
  tori $\T_i^F = \T^{(Fw_i^{-1})}$.}

\bigskip

\begin{center}
\begin{small}
\begin{tabular}{|l|l|} \hline
\rule{0cm}{0.4cm}
Torus & $\T^{(Fw_i^{-1})}$
\rule[-0.1cm]{0cm}{0.4cm}\\
\hline \cline{1-2} \hline
\rule{0cm}{0.5cm}
\hspace{-0.25cm} $\T_1^F$ & $\{ t_1(i, j) := h(\zetat^{i}, \zetat^{(2\theta-1)i}, \zetat^{j},
 \zetat^{(2\theta-1)j}) \, | \, 0 \le i,j \le q^2-2\}$
\rule[-0.2cm]{0cm}{0.6cm}\\
\hline
\rule{-0.12cm}{0.6cm}
$\T_2^F$ & $\{ t_2(i) := h(\zetat_4^{(4\theta^3+2\theta^2+1)i},
\zetat_4^{(2\theta+2\theta^2-1)i}, \zetat_4^{(2\theta-2\theta^2+1)i},
\zetat_4^{(2\theta^2-4\theta^3+1)i}) \, |$\\
& \hfill $\, 0 \le i \le q^4-2 \}$ 
\rule[-0.2cm]{0cm}{0.6cm}\\
\hline
\rule{-0.12cm}{0.6cm}
$\T_3^F$ & $\{ t_3(i) := h(\psit_8''^{(2\theta^2-2\theta+1)i},
\psit_8''^{(4\theta^3-6\theta^2+4\theta-1)i},
\psit_8''^{(2\theta^2-1)i}, \psit_8''^{(-4\theta^4+2\theta^2) i}) \, |$\\
& \hfill $\, 0 \le i \le q^4-\sqrt{2}q^3+\sqrt{2}q-2 \}$ 
\rule[-0.2cm]{0cm}{0.6cm}\\
\hline
\rule{-0.12cm}{0.6cm}
$\T_4^F$ & $\{ t_4(i) := h(\psit_8'^{(2\theta^2+2\theta+1)i},
\psit_8'^{(4\theta^3+2\theta^2-1)i},
\psit_8'^{(2\theta^2-1)i}, \psit_8'^{(-4\theta^4+2\theta^2) i}) \, |$\\
& \hfill $\, 0 \le i \le q^4+\sqrt{2}q^3-\sqrt{2}q-2 \}$ 
\rule[-0.2cm]{0cm}{0.6cm}\\
\hline
\end{tabular}
\end{small}
\end{center}

\newpage 

\bigskip

\noindent Table~B.4 \textit{Parameters for the elements of the maximal
  tori $\T_i^{*F^*} = (X \otimes_\Z \Q_{p'}/\Z)^{w_iF}$ of the dual
  group $\G^{*F^*}$. We use the abbreviations
  $\phi_{24}' := q^4+\sqrt{2}q^3+q^2+\sqrt{2}q+1$ and
  $\phi_{24}'' := q^4-\sqrt{2}q^3+q^2-\sqrt{2}q+1$.}

\bigskip

\begin{center}
\begin{scriptsize}
\begin{tabular}{|l|l|} \hline
\rule{0cm}{0.4cm}
\hspace{-0.25cm} Torus \hspace{-0.25cm} & $(X \otimes_\Z \Q_{p'}/\Z)^{w_iF}$
\rule[-0.1cm]{0cm}{0.4cm}\\
\hline \cline{1-2} \hline
\rule{0cm}{0.5cm}
\hspace{-0.25cm}
 $\T_1^{*F^*}$  \hspace{-0.25cm} & $\{ s_1(k,l) := s\left(
\frac{(2\theta+1)k+(2\theta+1)l}{q^2-1},
\frac{(4\theta+3)k+(2\theta+1)l}{q^2-1},
\frac{(6\theta+4)k+(2\theta+2)l}{q^2-1},
\frac{(2\theta+2)k+(2\theta+2)l}{q^2-1} \right) \, |$\\
& \hfill $\, 0 \le k,l \le q^2-2 \}$ 
\rule[-0.2cm]{0cm}{0.6cm}\\
\hline
\rule{-0.12cm}{0.6cm}
$\T_2^{*F^*}$  \hspace{-0.25cm} & $\{ s_2(k) := s\left(
\frac{2\theta^4(2\theta^3+2\theta^2+\theta)k}{q^4-1},
\frac{2\theta^4(4\theta^3+2\theta^2+2\theta+1)k}{q^4-1},
\frac{2\theta^4(4\theta^2+2\theta+4\theta^3+2)k}{q^4-1},
\frac{2\theta^4(2\theta^2+2\theta+1)k}{q^4-1} \right) |$\\
& \hfill $\, 0 \le k \le q^4-2 \}$ 
\rule[-0.2cm]{0cm}{0.6cm}\\
\hline
\rule{-0.12cm}{0.6cm}
$\T_3^{*F^*}$  \hspace{-0.25cm} & $\{ s_3(k) := s\left(
\frac{(6\theta^3-2\theta^2-\theta+2)k}{(q^2-1)(q^2-\sqrt{2}q+1)},
\frac{(10\theta^3-2\theta^2-2\theta+3)k}{(q^2-1)(q^2-\sqrt{2}q+1)},
\frac{(16\theta^3-4\theta^2-4\theta+5)k}{(q^2-1)(q^2-\sqrt{2}q+1)},
\frac{(8\theta^3-2\theta^2-2\theta+3)k}{(q^2-1)(q^2-\sqrt{2}q+1)} \right) |$\\
& \hfill $\, 0 \le k \le q^4-\sqrt{2}q^3+\sqrt{2}q-2 \}$ 
\rule[-0.2cm]{0cm}{0.6cm}\\
\hline
\rule{-0.12cm}{0.6cm}
$\T_4^{*F^*}$  \hspace{-0.25cm} & $\{ s_4(k) :=
s\left(\frac{(2\theta^3+2\theta^2+\theta)k}{(q^2-1)(q^2+\sqrt{2}q+1)},
\frac{(2\theta^3+4\theta^2+2\theta)k}{(q^2-1)(q^2+\sqrt{2}q+1)},
\frac{(4\theta^2+4\theta+1)k}{(q^2-1)(q^2+\sqrt{2}q+1)},
\frac{(2\theta^2+2\theta+1)k}{(q^2-1)(q^2+\sqrt{2}q+1)}\right) |$\\
\hline
\end{tabular}
\end{scriptsize}
\end{center}

\newpage

\noindent Table~B.5 \textit{Interpretation of the elements 
  $s_i(k, \dots)$ of the tori of the dual group as
  linear characters $\theta_i(k, \dots): \T_i^F \rightarrow \Q_{p'} / \Z$. 
  From those we get complex-valued linear characters of $\T_i^F$ via the 
  embedding $\varphi_2: \Q_{p'} \hookrightarrow \C^\times$, 
  $\frac{r}{s} \mapsto e^{2 \pi i r / s}$.}

\bigskip

\begin{center}
\begin{tabular}{|l|l|} \hline
\rule{0cm}{0.4cm}
\hspace{-0.25cm} Torus & $\theta_i(k, \dots): \T_i^F \rightarrow \Q_{p'} / \Z$
\rule[-0.1cm]{0cm}{0.4cm}\\
\hline \cline{1-2} \hline
\rule{0cm}{0.5cm}
\hspace{-0.25cm} 
$\hat{\T}_1^F$ & $\theta_1(k,l): t_1(a,b) \mapsto
\frac{(k+l)a+(k-l)b}{q^2-1}$ 
\rule[-0.2cm]{0cm}{0.6cm}\\
\hline
\rule{-0.12cm}{0.6cm}
$\hat{\T}_2^F$ & $\theta_2(k): t_2(a) \mapsto
\frac{ak}{q^4-1}$ 
\rule[-0.2cm]{0cm}{0.6cm}\\
\hline
\rule{-0.12cm}{0.6cm}
$\hat{\T}_3^F$ & $\theta_3(k): t_3(a) \mapsto
\frac{ak}{(q^2-1)(q^2-\sqrt{2}q+1)}$ 
\rule[-0.2cm]{0cm}{0.6cm}\\
\hline
\rule{-0.12cm}{0.6cm}
$\hat{\T}_4^F$ & $\theta_4(k): t_4(a) \mapsto
\frac{ak}{(q^2-1)(q^2+\sqrt{2}q+1)}$ 
\rule[-0.2cm]{0cm}{0.6cm}\\
\hline
\rule{-0.12cm}{0.6cm}
$\hat{\T}_5^F$ & $\theta_5(k,l): t_5(a,b) \mapsto
\frac{ak}{q^2-\sqrt{2}q+1}+\frac{bl}{q^2+\sqrt{2}q+1}$ 
\rule[-0.2cm]{0cm}{0.6cm}\\
\hline
\rule{-0.12cm}{0.6cm}
$\hat{\T}_6^F$ & $\theta_6(k,l): t_6(a,b) \mapsto
\frac{ak+bl}{q^2-\sqrt{2}q+1}$ 
\rule[-0.2cm]{0cm}{0.6cm}\\
\hline
\rule{-0.12cm}{0.6cm}
$\hat{\T}_7^F$ & $\theta_7(k,l): t_7(a,b) \mapsto
\frac{ak+bl}{q^2+\sqrt{2}q+1}$ 
\rule[-0.2cm]{0cm}{0.6cm}\\
\hline
\rule{-0.12cm}{0.6cm}
$\hat{\T}_8^F$ & $\theta_8(k,l): t_8(a,b) \mapsto
\frac{(k+l)a+(k-l)b}{q^2+1}$ 
\rule[-0.2cm]{0cm}{0.6cm}\\
\hline
\rule{-0.12cm}{0.6cm}
$\hat{\T}_9^F$ & $\theta_9(k): t_9(a) \mapsto
\frac{ak}{q^4-q^2+1}$ 
\rule[-0.2cm]{0cm}{0.6cm}\\
\hline
\rule{-0.12cm}{0.6cm}
$\hat{\T}_{10}^F$ & $\theta_{10}(k): t_{10}(a) \mapsto
\frac{ak}{q^4-\sqrt{2}q^3+q^2-\sqrt{2}q+1}$ 
\rule[-0.2cm]{0cm}{0.6cm}\\
\hline
\rule{-0.12cm}{0.6cm}
$\hat{\T}_{11}^F$ & $\theta_{11}(k): t_{11}(a) \mapsto
\frac{ak}{q^4+\sqrt{2}q^3+q^2+\sqrt{2}q+1}$ \\
\hline
\end{tabular}
\end{center}

\newpage

\noindent Table~B.6 \textit{Parameterization of the semisimple class types
  of ${^2F}_4(q^2)$.}

\bigskip

\begin{center}
\begin{tabular}{|l|c|l|l|} \hline
\rule{0cm}{0.4cm}
\hspace{-0.25cm} No. & Dynkin & $\Pi$ & $w$ \\
& type &&
\rule[-0.1cm]{0cm}{0.4cm}\\
\hline \cline{1-2} \hline
\rule{0cm}{0.5cm}
\hspace{-0.25cm} 
$h_1$ & ${^2F}_4$ & $[r_1,r_2,r_3,r_4]$ & $1$ 
\rule[-0.2cm]{0cm}{0.6cm}\\
\hline
\rule{-0.12cm}{0.6cm}
$h_2$ & ${^2B}_2$ & $[r_8,r_{16}]$ & $1$ 
\rule[-0.2cm]{0cm}{0.6cm}\\
\hline
\rule{-0.12cm}{0.6cm}
$h_3$ & $A_1$ & $[r_{15},r_{20}]$ & $1$ 
\rule[-0.2cm]{0cm}{0.6cm}\\
\hline
\rule{-0.12cm}{0.6cm}
$h_4$ & $A_0$ & $[ \, ]$ & $1$ 
\rule[-0.2cm]{0cm}{0.6cm}\\
\hline
\rule{-0.12cm}{0.4cm}
$h_5$ & ${^2A}_2$ & if $\theta \equiv 1$ mod $3$:& $w_{r_1}$ \\
&& $[r_8,r_9,r_{10},r_{18}]$ &\\
\hline
\end{tabular}
\end{center}

\newpage

\noindent Table~B.7 (cf.~Shinoda \cite[Table IV]{ShinodaClasses2F4} and
\cite[Table A.1]{HimstedtHuang2F4Borel}) \textit{Parameterization of
  the semisimple conjugacy classes of~${^2F_4(q^2)}$.} 

\begin{center}
\begin{tiny}
\begin{tabular}{|l|l|l|} \hline
\rule{0cm}{0.3cm}
{\normalsize Representative} & {\normalsize Parameters} & 
{\normalsize Number of \hspace{-0.25cm}} \\
&& {\normalsize Classes}
\rule[-0.1cm]{0cm}{0.3cm}\\
\hline
\rule{0cm}{0.3cm}
 $h_1 := h(1, 1, 1, 1)$ & & $1$
\rule[-0.2cm]{0cm}{0.3cm}\\
\hline
\rule{0cm}{0.3cm}
 $h_2(i) := h(1, 1, \zetat_2^{i}, \zetat_2^{(2\theta-1)i})$ & $i=0, \dots,
q^2-2$ & $\frac{q^2-2}{2}$ \\ 
& \, $i \neq 0$ &
\rule[-0.2cm]{0cm}{0.3cm}\\
\hline
\rule{0cm}{0.3cm}
 $h_3(i) := h(\zetat_2^i, \zetat_2^{(2\theta-1) i}, \zetat_2^{(2\theta+1) i},
\zetat_2^i)$ & $i=0, \dots, q^2-2$ & $\frac{q^2-2}{2}$ \\
& \, $i \neq 0$ &
\rule[-0.2cm]{0cm}{0.3cm}\\
\hline
\rule{0cm}{0.3cm}
 $h_4(i,j) := h({\zetat_2}^i, {\zetat_2}^{(2\theta-1)i}, {\zetat_2}^{j},
     {\zetat_2}^{(2\theta-1)j})$ & 
 $i=0, \dots, q^2-2$; $j=0, \dots, q^2-2$ & $\frac{q^4-10q^2+16}{16}$ \\
\hline
\end{tabular}
\end{tiny}
\end{center}

\newpage

\noindent Table~B.8 (K.~Shinoda~\cite[ Tables~II,
  V]{ShinodaClasses2F4} and \cite[Table A.2]{HimstedtHuang2F4Borel})
\textit{The conjugacy classes of~${^2F_4(q^2)}$.}

\enlargethispage{1cm}

\begin{center}
\begin{small}
\begin{tabular}{|l|l|l|} \hline
\rule{0cm}{0.3cm}
 {\normalsize \hspace{-0.3cm} Notation} & {\normalsize Representative} & {\normalsize
  $|C_{^2F_4(q^2)}|$} \\
\hline 
$c_{1,0}$ & $1$ & $q^{24} (q^{12}+1) (q^8-1) (q^6+1) (q^2-1)$\\
$c_{1,1}$ & $\alpha_{12}(1)$ & $q^{24} (q^4+1) (q^2-1)$\\
$c_{1,2}$ & $\alpha_{10}(1)$ & $q^{20} (q^4-1)$\\
$c_{1,3}$ & $\alpha_5(1) \alpha_{12}(1)$ & $2q^{14} (q^4+1) (q^2-1)$\\
$c_{1,4}$ & $\alpha_{5}(1)$ & $2q^{14} (q^4+1) (q^2-1)$\\
$c_{1,5}$ & $\alpha_{7}(1) \alpha_8(1)$ & $q^{16}$\\
$c_{1,6}$ & $\alpha_{5}(1) \alpha_7(1)$ & $q^{14}$\\
\hline
\end{tabular}
\end{small}
\end{center}

\newpage

\noindent Table~B.9 \textit{Parameterization of the types of Lusztig
  series of ${^2F}_4(q^2)$.}

\bigskip

\begin{center}
\begin{tabular}{|l|c|l|l|} \hline
\rule{0cm}{0.4cm}
\hspace{-0.25cm} No. & Dynkin & $\Pi$ & $w$ \\
& type &&
\rule[-0.1cm]{0cm}{0.4cm}\\
\hline \cline{1-2} \hline
\rule{0cm}{0.5cm}
\hspace{-0.25cm} 
$g_1$ & ${^2F}_4$ & $[r_1^\vee,r_2^\vee,r_3^\vee,r_4^\vee]$ & $1$ 
\rule[-0.2cm]{0cm}{0.6cm}\\
\hline
\rule{-0.12cm}{0.6cm}
$g_2$ & ${^2B}_2$ & $[r_8^\vee,r_{16}^\vee]$ & $1$ 
\rule[-0.2cm]{0cm}{0.6cm}\\
\hline
\rule{-0.12cm}{0.6cm}
$g_3$ & $A_1$ & $[r_{15}^\vee,r_{20}^\vee]$ & $1$ 
\rule[-0.2cm]{0cm}{0.6cm}\\
\hline
\rule{-0.12cm}{0.6cm}
$g_4$ & $A_0$ & $[ \, ]$ & $1$ 
\rule[-0.2cm]{0cm}{0.6cm}\\
\hline
\rule{-0.12cm}{0.4cm}
$g_5$ & ${^2A}_2$ & if $\theta \equiv 1$ mod $3$:& $w_{r_1}$ \\
&& $[r_6^\vee,r_{11}^\vee,r_{12}^\vee,r_{16}^\vee]$ &\\
\hline
\end{tabular}
\end{center}

\newpage

\noindent Table~B.10 \textit{Representatives for the semisimple
  conjugacy classes of $\G^{*F^*}$.}

\bigskip

\begin{center}
\begin{tiny}
\begin{tabular}{|l|l|l|} \hline
\rule{0cm}{0.3cm}
{\normalsize Representative} & {\normalsize Parameters} & 
{\normalsize Number of \hspace{-0.25cm}} \\
&& {\normalsize Classes}
\rule[-0.1cm]{0cm}{0.3cm}\\
\hline \cline{1-3} \hline
\rule{0cm}{0.3cm}
 $g_1 := s(0,0,0,0)$ & & $1$
\rule[-0.2cm]{0cm}{0.3cm}\\
\hline
\rule{0cm}{0.3cm}
 $g_2(k) := s\left(0, \frac{(2\theta+2)k}{q^2-1},
\frac{(4\theta+2)k}{q^2-1}, 0 \right)$ & $k=0, \dots,
q^2-2$ & $\frac{q^2-2}{2}$\\
& \, $k \neq 0$ &
\rule[-0.2cm]{0cm}{0.3cm}\\
\hline
\rule{0cm}{0.3cm}
 $g_3(k) := s\left(\frac{2\theta k}{q^2-1},
\frac{(4\theta+2)k}{q^2-1}, \frac{(4\theta+4) k}{q^2-1}, \frac{2
  k}{q^2-1}, \right)$ & $k=0, \dots, q^2-2$ & $\frac{q^2-2}{2}$ \\ 
& \, $k \neq 0$ &
\rule[-0.2cm]{0cm}{0.3cm}\\
\hline
\rule{0cm}{0.3cm}
 $g_4(k,l) := s_1(k,l)$ & 
 $k=0, \dots, q^2-2$; $l=0, \dots, q^2-2$ & $\frac{q^4-10q^2+16}{16}$ \\
\hline
\end{tabular}
\end{tiny}
\end{center}

\newpage

\noindent Table~B.11 \textit{Values of the class functions $\pm
  f_8(k)$ and $\pm f_{10}(k)$ in the Lusztig series of types $g_8$,
  $g_{10}$. All values not listed here are zero.} 

\bigskip

\begin{center}
\begin{small}
\begin{tabular}{|c|c|c|} \hline
\rule{0cm}{0.3cm}
& $c_{1,3}$ & $c_{1,4}$ 
\rule[-0.1cm]{0cm}{0.3cm}\\
\hline \cline{1-2} \hline
\rule{0cm}{0.4cm}
$f_8(k)$ & $\sqrt{2}q^5 \epsilon_4 (q^2+\sqrt{2}q+1)(q^2-1)$ &
$-\sqrt{2}q^5 \epsilon_4 (q^2+\sqrt{2}q+1)(q^2-1)$  
\rule[-0.2cm]{0cm}{0.3cm}\\
\hline
\rule{0cm}{0.4cm}
$f_{10}(k)$ & $\sqrt{2}q^5 \epsilon_4 (q^2-\sqrt{2}q+1)(q^2-1)$ & 
$-\sqrt{2}q^5 \epsilon_4 (q^2-\sqrt{2}q+1)(q^2-1)$ 
\rule[-0.2cm]{0cm}{0.3cm}\\
\hline
\end{tabular}

\bigskip

\begin{tabular}{|c|c|c|} \hline
\rule{0cm}{0.3cm}
& $c_{1,11}$ & $c_{1,12}$ 
\rule[-0.1cm]{0cm}{0.3cm}\\
\hline \cline{1-2} \hline
\rule{0cm}{0.4cm}
$f_8(k)$ & \hspace{0.68cm} $-\sqrt{2} q^2 \epsilon_4 (2q+\sqrt{2})$ \hspace{0.68cm} &
 \hspace{0.98cm}$\sqrt{2} q^2 \epsilon_4 (2q+\sqrt{2})$  \hspace{0.98cm}
\rule[-0.2cm]{0cm}{0.3cm}\\
\hline
\rule{0cm}{0.4cm}
$f_{10}(k)$ & $\sqrt{2} q^2 \epsilon_4 (2q-\sqrt{2})$ &
$-\sqrt{2} q^2 \epsilon_4 (2q-\sqrt{2})$
\rule[-0.2cm]{0cm}{0.3cm}\\
\hline
\end{tabular}

\bigskip

\begin{tabular}{|c|c|c|c|c|} \hline
\rule{0cm}{0.3cm}
& $c_{1,15}$ & $c_{1,16}$ & $c_{1,17}$ & $c_{1,18}$ 
\rule[-0.1cm]{0cm}{0.3cm}\\
\hline \cline{1-2} \hline
\rule{0cm}{0.4cm}
$f_8(k)$ & \hspace{0.36cm} $-\sqrt{2} q \epsilon_4$ \hspace{0.36cm}
& \hspace{0.36cm} $\sqrt{2} q \epsilon_4$ \hspace{0.36cm} & \hspace{0.425cm} 
$-\sqrt{2} q \epsilon_4$ \hspace{0.425cm} & \hspace{0.425cm} $\sqrt{2} q
\epsilon_4$ \hspace{0.425cm} 
\rule[-0.2cm]{0cm}{0.3cm}\\
\hline
\rule{0cm}{0.4cm}
$f_{10}(k)$ & $-\sqrt{2} q \epsilon_4$ & $\sqrt{2} q \epsilon_4$ &
$-\sqrt{2} q \epsilon_4$ & $\sqrt{2} q \epsilon_4$
\rule[-0.2cm]{0cm}{0.3cm}\\
\hline
\end{tabular}

\bigskip

\begin{tabular}{|c|c|} \hline
\rule{0cm}{0.3cm}
& $c_{8,2}(i)$ 
\rule[-0.1cm]{0cm}{0.3cm}\\
\hline \cline{1-2} \hline
\rule{0cm}{0.4cm}
$f_8(k)$ & \hspace{0.1cm} $\sqrt{2}q \epsilon_4
(\varphi_8''^{qik/\sqrt{2}} + \varphi_8''^{-qik/\sqrt{2}} +
\varphi_8''^{(q/\sqrt{2}-1)ik} + \varphi_8''^{-(q/\sqrt{2}-1)ik})$ \hspace{0.1cm} 
\rule[-0.2cm]{0cm}{0.3cm}\\
\hline
\rule{0cm}{0.4cm}
$f_{10}(k)$ & $0$
\rule[-0.2cm]{0cm}{0.3cm}\\
\hline
\end{tabular}

\bigskip

\begin{tabular}{|c|c|} \hline
\rule{0cm}{0.3cm}
& $c_{8,3}(i)$ 
\rule[-0.1cm]{0cm}{0.3cm}\\
\hline \cline{1-2} \hline
\rule{0cm}{0.4cm}
$f_8(k)$ & $-\sqrt{2}q \epsilon_4
(\varphi_8''^{qik/\sqrt{2}} + \varphi_8''^{-qik/\sqrt{2}} +
\varphi_8''^{(q/\sqrt{2}-1)ik} +
\varphi_8''^{-(q/\sqrt{2}-1)ik})$ \hspace{0.05cm}
\rule[-0.2cm]{0cm}{0.3cm}\\
\hline
\rule{0cm}{0.4cm}
$f_{10}(k)$ & $0$
\rule[-0.2cm]{0cm}{0.3cm}\\
\hline
\end{tabular}

\bigskip

\begin{tabular}{|c|c|} \hline
\rule{0cm}{0.3cm}
& $c_{10,2}(i)$ 
\rule[-0.1cm]{0cm}{0.3cm}\\
\hline \cline{1-2} \hline
\rule{0cm}{0.4cm}
$f_8(k)$ & $0$
\rule[-0.2cm]{0cm}{0.3cm}\\
\hline
\rule{0cm}{0.4cm}
$f_{10}(k)$ & \hspace{0.25cm} $\sqrt{2}q \epsilon_4
(\varphi_8'^{qik/\sqrt{2}} + \varphi_8'^{-qik/\sqrt{2}} +
\varphi_8'^{(q/\sqrt{2}+1)ik} + \varphi_8'^{-(q/\sqrt{2}+1)ik})$ \hspace{0.25cm} 
\rule[-0.2cm]{0cm}{0.3cm}\\
\hline
\end{tabular}

\bigskip

\begin{tabular}{|c|c|} \hline
\rule{0cm}{0.3cm}
& $c_{10,3}(i)$ 
\rule[-0.1cm]{0cm}{0.3cm}\\
\hline \cline{1-2} \hline
\rule{0cm}{0.4cm}
$f_8(k)$ & $0$
\rule[-0.2cm]{0cm}{0.3cm}\\
\hline
\rule{0cm}{0.4cm}
$f_{10}(k)$ &  \hspace{0.12cm} $-\sqrt{2}q \epsilon_4
(\varphi_8'^{qik/\sqrt{2}} + \varphi_8'^{-qik/\sqrt{2}} +
\varphi_8'^{(q/\sqrt{2}+1)ik} + \varphi_8'^{-(q/\sqrt{2}+1)ik})$ \hspace{0.12cm}
\rule[-0.2cm]{0cm}{0.3cm}\\
\hline
\end{tabular}
\end{small}
\end{center}

\newpage

\noindent Table~B.12 \textit{Values of the irreducible characters
  ${_G\chi}_{42}(k)$, \dots, ${_G\chi}_{49}(k)$ of the group 
  $G := {^2F}_4(q^2)$ on some conjugacy classes.}

\bigskip

\begin{center}
$
\begin{array}{
|c|c|
}
\hline

& c_{1,0}
\rule[-0.1cm]{0cm}{0.3cm}\\
\hline \cline{1-2} \hline
\rule{0cm}{0.4cm}
{_G\chi}_{
42
}
(
 k
)
&
(q-1)(q+1)(q^2+1)^{2}(q^4-q^2+1)(q^{8}-q^{4}+1)(q^{2}+\sqrt{2}q+1)
\\ 
\hline
{_G\chi}_{
43
}
(
 k
)
&
\frac{q}{\sqrt{2}}(q-1)^{2}(q+1)^{2}(q^2+1)^{2}(q^4-q^2+1)(q^{8}-
q^{4}+1)(q^{2}+\sqrt{2}q+1)
\\ 
\hline
{_G\chi}_{
44
}
(
 k
)
&
\frac{q}{\sqrt{2}}(q-1)^{2}(q+1)^{2}(q^2+1)^{2}(q^4-q^2+1)(q^{8}-
q^{4}+1)(q^{2}+\sqrt{2}q+1)
\\ 
\hline
{_G\chi}_{
45
}
(
 k
)
&
q^{4}(q-1)(q+1)(q^2+1)^{2}(q^4-q^2+1)(q^{8}-q^{4}+1)(q^{2}+\sqrt{2}q+1
)
\\ 
\hline
{_G\chi}_{
46
}
(
 k
)
&
(q-1)(q+1)(q^2+1)^{2}(q^4-q^2+1)(q^{8}-q^{4}+1)(q^{2}-\sqrt{2}q+1)
\\ 
\hline
{_G\chi}_{
47
}
(
 k
)
&
\frac{q}{\sqrt{2}}(q-1)^{2}(q+1)^{2}(q^2+1)^{2}(q^4-q^2+1)(q^{8}-
q^{4}+1)(q^{2}-\sqrt{2}q+1)
\\ 
\hline
{_G\chi}_{
48
}
(
 k
)
&
\frac{q}{\sqrt{2}}(q-1)^{2}(q+1)^{2}(q^2+1)^{2}(q^4-q^2+1)(q^{8}-
q^{4}+1)(q^{2}-\sqrt{2}q+1)
\\ 
\hline
{_G\chi}_{
49
}
(
 k
)
&
q^{4}(q-1)(q+1)(q^2+1)^{2}(q^4-q^2+1)(q^{8}-q^{4}+1)(q^{2}-\sqrt{2}q+1
)
\\ 
\hline
\end{array}$ \bigskip
\medskip

$
\begin{array}{
|c|c|
}
\hline

&c_{
1,11
}
\rule[-0.1cm]{0cm}{0.3cm}\\
\hline \cline{1-2} \hline
\rule{0cm}{0.4cm}
{_G\chi}_{
42
}
(
 k
)
&
-(q^{2}+\sqrt{2}q+1)
\\ 
\hline
{_G\chi}_{
43
}
(
 k
)
&
\hspace{3.1cm} q^{2}+\frac{1}{2}\sqrt{2}q-\varepsilon_4\sqrt{2}q^{3}
-\varepsilon_4q^{2} \hspace{3.1cm} 
\\ 
\hline
{_G\chi}_{
44
}
(
 k
)
&
q^{2}+\frac{1}{2}\sqrt{2}q+\varepsilon_4\sqrt{2}q^{3}+
q^{2}\varepsilon_4
\\ 
\hline
{_G\chi}_{
45
}
(
 k
)
&
0
\\ 
\hline
{_G\chi}_{
46
}
(
 k
)
&
-(q^{2}-\sqrt{2}q+1)
\\ 
\hline
{_G\chi}_{
47
}
(
 k
)
&
-q^{2}+\frac{1}{2}\sqrt{2}q+\varepsilon_4\sqrt{2}q^{3}
-\varepsilon_4q^{2}
\\ 
\hline
{_G\chi}_{
48
}
(
 k
)
&
-q^{2}+\frac{1}{2}\sqrt{2}q-\varepsilon_4\sqrt{2}q^{3}+
q^{2}\varepsilon_4
\\ 
\hline
{_G\chi}_{
49
}
(
 k
)
&
0
\\ 
\hline
\end{array}$ \bigskip
\medskip

$
\begin{array}{
|c|c|
}
\hline

&c_{
8,2
}
(
 i
)
\rule[-0.1cm]{0cm}{0.3cm}\\
\hline \cline{1-2} \hline
\rule{0cm}{0.4cm} 
{_G\chi}_{
42
}
(
 k
)
&
\parbox{
90.82882mm
}{\centering 
\mbox{}\\ $
-\varphi_{8}''^{\frac{1}{2} iq\sqrt{2} k- i k}
-\varphi_{8}''^{\sqrt{2}q i k- i k}
-\varphi_{8}''^{-\sqrt{2}q i k+ i k}
-\varphi_{8}''^{-\frac{1}{2} iq\sqrt{2} k+ i k}
$ 
 \\  $ 
-\varphi_{8}''^{\frac{1}{2} iq\sqrt{2} k}-\varphi_{8}''^{-\frac{1}{2}
  iq\sqrt{2} k}-\varphi_{8}''^{- i k}-\varphi_{8}''^{ i k}-1 
$
 \\[1mm] }
\\ 
\hline
{_G\chi}_{
43
}
(
 k
)
&
\parbox{
90.611126mm
}{\centering 
\mbox{}\\ $
\frac{1}{2}q\sqrt{2}\varepsilon_4
\varphi_{8}''^{\frac{1}{2} iq\sqrt{2} k}
+\frac{1}{2}q\sqrt{2}\varepsilon_4
\varphi_{8}''^{-\frac{1}{2} iq\sqrt{2} k+ i k}$ 
 \\  $ 
+\frac{1}{2}q\sqrt{2}\varepsilon_4
\varphi_{8}''^{\frac{1}{2} iq\sqrt{2} k- i k}+
\frac{1}{2}q\sqrt{2}\varepsilon_4
\varphi_{8}''^{-\frac{1}{2} iq\sqrt{2} k} 
$ 
 \\  $ 
+\varphi_{8}''^{\frac{1}{2} iq\sqrt{2} k- i k}
+\varphi_{8}''^{-\frac{1}{2} iq\sqrt{2} k+ i k}
+\varphi_{8}''^{\frac{1}{2} iq\sqrt{2} k}+
\varphi_{8}''^{-\frac{1}{2} iq\sqrt{2} k}+1
$
 \\[1mm] }
\\ 
\hline
{_G\chi}_{
44
}
(
 k
)
&
\parbox{
90.611126mm
}{\centering 
\mbox{}\\ $
-\frac{1}{2}q\sqrt{2}\varepsilon_4
\varphi_{8}''^{\frac{1}{2} iq\sqrt{2} k}
-\frac{1}{2}q\sqrt{2}\varepsilon_4
\varphi_{8}''^{-\frac{1}{2} iq\sqrt{2} k+ i k}$ 
 \\  $ 
-\frac{1}{2}q\sqrt{2}\varepsilon_4
\varphi_{8}''^{\frac{1}{2} iq\sqrt{2} k- i k}
-\frac{1}{2}q\sqrt{2}\varepsilon_4
\varphi_{8}''^{-\frac{1}{2} iq\sqrt{2} k}+
\varphi_{8}''^{\frac{1}{2} iq\sqrt{2} k- i k}
$ 
 \\  $ 
+\varphi_{8}''^{-\frac{1}{2} iq\sqrt{2} k+ i k}
+\varphi_{8}''^{\frac{1}{2} iq\sqrt{2} k}+
\varphi_{8}''^{-\frac{1}{2} iq\sqrt{2} k}+1
$
 \\[1mm] }
\\ 
\hline
{_G\chi}_{
45
}
(
 k
)
&
\hspace{0.83cm} 
1+\varphi_{8}''^{\frac{1}{2} iq\sqrt{2} k- i k}+
\varphi_{8}''^{-\frac{1}{2} iq\sqrt{2} k+ i k}+
\varphi_{8}''^{\frac{1}{2} iq\sqrt{2} k}+
\varphi_{8}''^{-\frac{1}{2} iq\sqrt{2} k} \hspace{0.83cm} 
\\
\hline
\end{array}$ \bigskip
\end{center}

\newpage

\noindent Table~B.12 (cont.) 

\smallskip

\begin{center}
$
\begin{array}{
|c|c|
}
\hline

&c_{
10,2
}
(
 i
)
\rule[-0.1cm]{0cm}{0.3cm}\\
\hline \cline{1-2} \hline
\rule{0cm}{0.4cm}
{_G\chi}_{
42
}
(
 k
)
&
-1
\\ 
\hline
{_G\chi}_{
43
}
(
 k
)
&
-1
\\ 
\hline
{_G\chi}_{
44
}
(
 k
)
&
-1
\\ 
\hline
{_G\chi}_{
45
}
(
 k
)
&
1
\\ 
\hline
{_G\chi}_{
46
}
(
 k
)
&
\parbox{
90.82882mm
}{\centering 
\mbox{}\\ $
-\varphi_{8}'^{-\sqrt{2}q i k- i k}
-\varphi_{8}'^{-\frac{1}{2} iq\sqrt{2} k- i k}
-\varphi_{8}'^{\frac{1}{2} iq\sqrt{2} k+ i k}
-\varphi_{8}'^{\sqrt{2}q i k+ i k}$ 
 \\  $ 
-\varphi_{8}'^{\frac{1}{2} iq\sqrt{2} k}
-\varphi_{8}'^{-\frac{1}{2} iq\sqrt{2} k}
-\varphi_{8}'^{- i k}-\varphi_{8}'^{ i k}-1
$
 \\[1mm] }
\\ 
\hline
{_G\chi}_{
47
}
(
 k
)
&
\parbox{
90.658722mm
}{\centering 
\mbox{}\\ $
\frac{1}{2}q\sqrt{2}\varepsilon_4
\varphi_{8}'^{\frac{1}{2} iq\sqrt{2} k}+
\frac{1}{2}q\sqrt{2}\varepsilon_4
\varphi_{8}'^{-\frac{1}{2} iq\sqrt{2} k- i k}$ 
 \\  $ 
+\frac{1}{2}q\sqrt{2}\varepsilon_4
\varphi_{8}'^{-\frac{1}{2} iq\sqrt{2} k}+
\frac{1}{2}q\sqrt{2}\varepsilon_4
\varphi_{8}'^{\frac{1}{2} iq\sqrt{2} k+ i k}
-\varphi_{8}'^{-\frac{1}{2} iq\sqrt{2} k- i k}
$ 
 \\  $ 
-\varphi_{8}'^{\frac{1}{2} iq\sqrt{2} k+ i k}
-\varphi_{8}'^{-\frac{1}{2} iq\sqrt{2} k}
-\varphi_{8}'^{\frac{1}{2} iq\sqrt{2} k}-1
$
 \\[1mm] }
\\ 
\hline
{_G\chi}_{
48
}
(
 k
)
&
\parbox{
90.658722mm
}{\centering 
\mbox{}\\ $
-\frac{1}{2}q\sqrt{2}\varepsilon_4
\varphi_{8}'^{\frac{1}{2} iq\sqrt{2} k}
-\frac{1}{2}q\sqrt{2}\varepsilon_4
\varphi_{8}'^{-\frac{1}{2} iq\sqrt{2} k- i k}
$ 
 \\  $ 
-\frac{1}{2}q\sqrt{2}\varepsilon_4
\varphi_{8}'^{-\frac{1}{2} iq\sqrt{2} k}
-\frac{1}{2}q\sqrt{2}\varepsilon_4
\varphi_{8}'^{\frac{1}{2} iq\sqrt{2} k+ i k}
-\varphi_{8}'^{-\frac{1}{2} iq\sqrt{2} k- i k}
$ 
 \\  $ 
-\varphi_{8}'^{\frac{1}{2} iq\sqrt{2} k+ i k}
-\varphi_{8}'^{-\frac{1}{2} iq\sqrt{2} k}
-\varphi_{8}'^{\frac{1}{2} iq\sqrt{2} k}-1
$
 \\[1mm] }
\\
\hline
\end{array}$
\end{center}

\end{document}